\documentclass[a4paper,review]{cas-sc}

\usepackage[section]{placeins} 
\usepackage{amssymb}
\usepackage{amsthm}
\usepackage{amsmath}
\usepackage{upgreek}
\usepackage{pdflscape}
\usepackage{listings}
\usepackage{multirow}
\usepackage[percent]{overpic}
\usepackage{multicol}
\usepackage{color}
\usepackage{stmaryrd}
\usepackage{xfrac}
\usepackage[capitalise]{cleveref}
\usepackage{booktabs}
\usepackage[section]{placeins} 
\usepackage{algorithm}
\usepackage{algorithmic}
\usepackage{calc}
\usepackage[titletoc]{appendix}
\usepackage{siunitx}
\usepackage[sort&compress,square,numbers]{natbib}

\usepackage{tikz-dimline}
\usepackage{graphicx}
\usepackage{graphics}
\usepackage{wrapfig}
\usepackage{float}
\usepackage{subfig}
\usepackage[percent]{overpic}
\usepackage{varwidth}
\usepackage{tikz}
\usepackage{pgfplots}
\usepackage{adjustbox}
\usetikzlibrary{arrows,matrix,positioning,fit}
\usetikzlibrary{shapes,positioning}
\usetikzlibrary{backgrounds}
\usepackage{tikz-layers}
\usepgfplotslibrary{fillbetween}
\usetikzlibrary{intersections}
\pgfplotsset{compat=1.14}
\usepgfplotslibrary{colorbrewer}
\usepgfplotslibrary{patchplots}
\usepgfplotslibrary[colorbrewer]
\usetikzlibrary{pgfplots.colorbrewer}
\usetikzlibrary[pgfplots.colorbrewer]
\usepgfplotslibrary{units}
\usetikzlibrary{spy}
\usepackage{pgfplotstable}
\usepackage{arrayjobx}
\graphicspath{ {./figs/} }
\usetikzlibrary{external}

\newlength\myheight
\newlength\mydepth
\settototalheight\myheight{Xygp}
\newcommand*\inlinegraphics[1]{%
  \settototalheight\myheight{Xygp}%
  \settodepth\mydepth{Xygp}%
  \raisebox{-\mydepth}{\includegraphics[height=\myheight]{#1}}%
}
\newcommand\orcid[1]{\href{https://orcid.org/#1}{\inlinegraphics{orcid_16x16.png}}}

\makeatletter
\def\BState{\State\hskip-\ALG@thistlm}
\makeatother

\newtheorem{theorem}{Theorem}[section]
\newtheorem{corollary}{Corollary}[theorem]
\newtheorem{lemma}[theorem]{Lemma}
\newdefinition{definition}{Definition}[section]

\newcommand\pxvar[2]{\partial_{#2} #1}

\newcommand{\gu}{g(\mathbf{u})}
\newcommand{\gux}{g(\mathbf{u}(\mathbf{x}))}
\newcommand{\gbar}{g(\overline{\mathbf{u}})}
\newcommand{\ux}{\mathbf{u}(\mathbf{x})}
\newcommand{\ubar}{\overline{\mathbf{u}}}

\newcommand{\gstar}{g(\mathbf{u}(\mathbf{x}^*))}
\newcommand{\hu}{h(\mathbf{u})}
\newcommand{\hux}{h(\mathbf{u}(\mathbf{x}))}
\newcommand{\hstar}{h(\mathbf{u}(\mathbf{x}^*))}
\newcommand{\ghat}{g(\widehat{\mathbf{u}}(\mathbf{x}))}
\newcommand{\uhat}{\widehat{\mathbf{u}}(\mathbf{x})}

\begin{document}

\title[mode=title]{Continuously bounds-preserving discontinuous Galerkin methods for hyperbolic conservation laws}
\shorttitle{Continuously bounds-preserving discontinuous Galerkin methods for hyperbolic conservation laws}
\shortauthors{T. Dzanic}

\author[1]{T. Dzanic}[orcid=0000-0003-3791-1134]
\cormark[1]
\cortext[cor1]{Corresponding author}
\ead{dzanic1@llnl.gov}

\address[1]{Lawrence Livermore National Laboratory, Livermore, CA 94551, United States of America}

\begin{abstract}
For finite element approximations of transport phenomena, it is often necessary to apply a form of limiting to ensure that the discrete solution remains well-behaved and satisfies physical constraints. However, these limiting procedures are typically performed at discrete nodal locations, which is not sufficient to ensure the robustness of the scheme when the solution must be evaluated at arbitrary locations (e.g., for adaptive mesh refinement, remapping in arbitrary Lagrangian--Eulerian solvers, overset meshes, etc.). In this work, a novel limiting approach for discontinuous Galerkin methods is presented which ensures that the solution is \emph{continuously} bounds-preserving (i.e., across the entire solution polynomial) for any arbitrary choice of basis, approximation order, and mesh element type. Through a modified formulation for the constraint functionals, the proposed approach requires only the solution of a single spatial scalar minimization problem per element for which a highly efficient numerical optimization procedure is presented. The efficacy of this approach is shown in numerical experiments by enforcing continuous constraints in high-order unstructured discontinuous Galerkin discretizations of hyperbolic conservation laws, ranging from scalar transport with maximum principle preserving constraints to compressible gas dynamics with positivity-preserving constraints.
\end{abstract}

\begin{keywords}
Discontinuous Galerkin \sep
High-order\sep
Hyperbolic conservation laws\sep
Bounds-preserving\sep
Positivity-preserving\sep
Limiting
\end{keywords}



\maketitle



\section{Introduction}
\label{sec:intro}
The efficient and accurate computation of complex transport phenomena remains a driving force in the development of high-fidelity numerical discretization techniques. Within this broad application area, the use of high-order schemes has grown in prevalence, primarily due to their ability to obtain higher accuracy per degree of freedom which is advantageous in the context of scale-resolving simulations. Various high-order numerical methods have been routinely applied to simulating transport phenomena, with these schemes typically falling within the overarching classes of finite difference, finite volume, and finite element methods. Of these approaches, finite element methods offer notable benefits as they can achieve arbitrarily high-order accuracy and are straightforward to generalize to unstructured meshes. One particular formulation of finite element methods, discontinuous Galerkin (DG)~\citep{Hesthaven2008DG} schemes and variants thereof (e.g., flux reconstruction~\citep{Huynh2007}), offers further advantages as they are typically better suited for approximating transport phenomena and possess a compact data structure that is beneficial for modern massively-parallel high performance computing architectures. 
As such, these characteristics have resulted in the successful application of DG schemes to a wide variety of problems ranging from fluid dynamics and astrophysics to solid mechanics and bioengineering~\citep{Cox2019,Markert2022,Zhu2008,Nguyen2014}. 

However, the accuracy and robustness of DG schemes (and finite element methods in general) can become severely degraded when the physical systems in question must abide by strict physical constraints, which is frequently the case in transport-dominated applications (e.g., positivity of density and internal energy in gas dynamics, subluminal velocity in relativistic magnetohydrodynamics, bounded mass fractions in multi-species flows, etc.). Standard numerical implementations generally do not guarantee that these constraints are satisfied by the discrete solution, which may result in inaccurate and physically inconsistent predictions or, more typically, the failure of the scheme altogether. This problem is further exacerbated in more complex applications such as multi-physics modeling, where the strong coupling between various physical models can be particularly sensitive to ill-behaved numerical solutions. 

A typical remedy for this problem is to enforce these constraints via some form of \textit{a posteriori} limiting on the discrete solution to ensure that the solution remains \emph{bounds-preserving} (i.e., constraint-satisfying). These methods generally rely on utilizing a secondary less accurate (but more robust) numerical approximation of the governing equations which is guaranteed to be bounds-preserving. The solutions as predicted by the high-order baseline scheme and the low-order secondary scheme are then combined to yield a bounds-preserving approximation that ostensibly retains the accuracy of the high-order scheme. Within the overarching framework of limiting methods, the quintessential example is the flux-corrected transport method of \citet{Boris1997} from which a variety of approaches are derived~\citep{Guermond2016,Anderson2017,Pazner2021,Lin2023,RuedaRamrez2022,Peyvan2023,Kuzmin2020,Dumbser2016,Lin2024}. Alternate limiting formulations which do not rely on a secondary low-order scheme have also been developed by using bounds-preserving estimates derived from the high-order solution~\citep{Zhang2010,Zhang2011,Zhang2011b, Dzanic2022,Ching2023,Dzanic2023,Pandare2023,Zhang2023} or through slope-limiting/stencil modification of the underlying solution~\citep{Barth1989,Cockburn1989,Qiu2004}.

However, a critical drawback of these limiting approaches is that they are typically performed \emph{discretely}, whereas the solution for the DG approximation is represented (piecewise) continuously. The consequence is that the schemes are only provably bounds-preserving at the discrete spatial locations where the limiting was performed. For many applications, this is not sufficient -- it is sometimes necessary to have to evaluate the solution at arbitrary locations. For example, adaptive mesh refinement frequently requires evaluating the solution at a new set of quadrature points which do not coincidence with the quadrature points used for the element. This is similarly the case for coupled multi-physics solvers that rely on independent meshes, where the solution must be frequently transferred between these meshes. In some cases, the location that the solution must be evaluated at may not even be known \emph{a priori}, such as in the remap stage of an arbitrary Lagrangian--Eulerian solver or with dynamic overset meshes. In other cases, the solution may not even be queried in real time, such as when coupling of physical models is done as a post-processing step after the end of the simulation, making it impossible to perform \textit{a posteriori} limiting based on the location where the solution will be evaluated. For such applications, the limited solution may still violate the constraints as it is not bounds-preserving at these arbitrary points, causing the failure of the scheme and/or the breakdown of physical models. 

This issue is typically counteracted by changing the basis of the underlying finite element method such that the solution remains continuously bounded (i.e., across the entire solution polynomial). One such \emph{continuously bounds-preserving} approach can be obtained by ensuring that the solution approximation is at most linear, but this effectively limits the scheme to second-order accuracy. Alternatively, one may use bounded shape functions such Bernstein polynomials with appropriate limiting to guarantee that the solution is continuously bounds-preserving even for high-order approximations~\citep{Anderson2017, Glaubitz2019}. However, the use of bounded Bernstein polynomials for the solution basis gives rise to some notable drawbacks, such as numerical difficulties due to the ill-conditioned basis transformation operation and the lack of flexibility of the underlying approximation space which does not easily extend to more complex bases (e.g., divergence-free basis functions, rational basis functions, etc.). Additionally, limiting based on the convex hull property of Bernstein polynomials can introduce excessive numerical dissipation as it generally only provides a \emph{sufficient} degree of limiting, not the \emph{necessary} degree of limiting, even for linear constraints and for solutions which already continuously satisfy constraints. As such, there is potential to significantly decrease the computational cost and increase the accuracy and applicability of DG schemes if the underlying discretization can be modified such as to continuously enforce constraints for general DG discretizations.

In this work, a novel \textit{a posteriori} limiting procedure for DG schemes is presented which ensures that the solution is continuously bounds-preserving for any arbitrary choice of basis, approximation order, and mesh element type. The proposed method extends the limiting approach of \citet{Zhang2010} to continuously enforce general algebraic constraints on the high-order solution. We show that limiting based off of the continuous extrema of the constraint functionals is not sufficient to enforce constraints continuously and propose a modified formulation for constraints which can guarantee a continuously bounds-preserving solution with only a single spatial minimization problem per element. This approach enforces the exact amount of limiting necessary for linear constraints and a sufficient amount of limiting for nonlinear constraints. Furthermore, a highly efficient numerical optimization procedure for solving the required minimization problem is presented, although the non-convexity of the problem does not necessarily mathematically guarantee convergence to the global minimum. The proposed approach is applied to a high-order unstructured DG scheme for hyperbolic conservation laws ranging from scalar transport to compressible gas dynamics. 

The remainder of this manuscript is organized as follows. In \cref{sec:preliminaries}, some preliminaries on hyperbolic conservation laws, convex invariants and constraints, and DG schemes are presented. The proposed numerical approach is then introduced in \cref{sec:methodology}, followed by numerical implementation details in \cref{sec:implementation} and the results of the numerical experiments in \cref{sec:results}. Conclusions are then drawn in \cref{sec:conclusion}. 
\section{Preliminaries}
\label{sec:preliminaries}

\subsection{Hyperbolic conservation laws and convex constraints}
This work pertains to hyperbolic conservation laws of the form 
\begin{equation}\label{eq:gen_hype}
     \begin{cases}
        \pxvar{\mathbf{u}(\mathbf{x}, t)}{t} + \boldsymbol{\nabla}\cdot\mathbf{F}(\mathbf{u}) = 0, \quad \mathrm{for}\ (\mathbf{x}, t)\in\Omega\times\mathbb{R}_+, \\
        \mathbf{u}(\mathbf{x}, 0)= \mathbf{u}_0(\mathbf{x}), \quad \quad  \quad \quad  \  \mathrm{for}\ \mathbf{x}\in \Omega,
    \end{cases}
\end{equation}
where $\Omega\subset\mathbb{R}^d$ is a $d$-dimensional spatial domain, $\mathbf{u}\in\mathbb{R}^m$ is a vector-valued solution of $m$ variables, and $\mathbf{F}(\mathbf{u})\in \mathbb{R}^{m \times d}$ is the associated flux. We assume that there exists a contiguous set of \emph{admissible} solutions to \cref{eq:gen_hype}, denoted by $\mathcal G$, which possess certain properties that are desirable for numerical approximations of these problems. Without giving a precise definition to what an admissible solution may be, which in many cases remains an open problem, we assume that these admissible solutions must abide by some strict constraints, formalized by the positivity of some set of $N_c$ linear or nonlinear constraint functionals of $\mathbf{u}$ as
\begin{equation}\label{eq:constraints}
    g_i(\mathbf{u}) \geq 0 \ \forall \ i \in \lbrace 1, \dotsc, N_c \rbrace . 
\end{equation}
Furthermore, we assume that this set $\mathcal G$ is convex, such that the constraint functionals are quasiconcave (i.e., $\mathcal G$ is the upper contour set of quasiconcave functionals). The convexity of the admissible set of solutions is closely related to the invariant region of many hyperbolic conservation laws (e.g., the set of solutions corresponding to the positivity of density and internal energy in gas dynamics is convex), the details of which have been explored by \citet{Lax1954}, \citet{Hoff1985}, and \citet{Frid2001}. The goal of this work is to ensure that the numerical approximation of these systems satisfies the given constraints stemming from convex invariants of the system. For brevity, only weak inequalities are considered, but this work readily extends to strict inequalities as they are effectively identical from the perspective of finite precision arithmetic. 

\subsection{Discontinuous Galerkin methods}
The numerical method of choice for approximating \cref{eq:gen_hype} is the DG method~\citep{Hesthaven2008DG}. In this approach, the domain $\Omega$ is partitioned into $N$ elements $\Omega_k$ such that $\Omega = \bigcup_N\Omega_k$ and $\Omega_i\cap\Omega_j=\emptyset$ for $i\neq j$. With a slight abuse of notation, the discrete solution $\mathbf{u}(\mathbf{x})$ within each element $\Omega_k$ is approximated via a summation of a set of $n_s$ basis functions as 
\begin{equation}
    \mathbf{u}(\mathbf{x}) = \sum_{i=1}^{n_s} \mathbf{u}_{i}\phi_i(\mathbf{x})\subset V_h,
\end{equation}
where $\phi_i(\mathbf{x})$ are the basis functions, $\mathbf{u}_i$ are their associated coefficients, and $V_h$ is the piece-wise polynomial space spanned by the basis functions. For polynomial approximations, which are the focus of this work, we denote $\mathbb P_p$ as the order of the approximation, where $p$ is the maximal order of $\mathbf{u}(\mathbf{x})$.

The DG formulation involves integrating \cref{eq:gen_hype} with respect to a test function $\mathbf{w}(\mathbf{x}) \subset V_h$, which resides in the same finite element space as $\mathbf{u}(\mathbf{x})$, yielding
\begin{equation}\label{eq:semi-disc}
    \sum_{k=1}^{N} \left \{ \int_{\Omega_k} \partial_t \mathbf{u}{\cdot} \mathbf{w}\ \mathrm{d}V  + \int_{\partial \Omega_k} \hat{\mathbf{F}}(\mathbf{u}^-, \mathbf{u}^+, \mathbf{n}) {\cdot} \mathbf{w}\ \mathrm{d}S  -  \int_{\Omega_k} \mathbf{F}(\mathbf{u}){\cdot} \nabla \mathbf{w}\ \mathrm{d}V\right \} = 0.
\end{equation}
Along the element interfaces $\partial \Omega$, the flux function is replaced by a numerical interface flux $\hat{\mathbf{F}}(\mathbf{u}^-, \mathbf{u}^+, \mathbf{n})$. This interface flux is dependent on the solution within the element of interest (denoted by the $\mathbf{u}^-$), its face-adjacent neighbor (denoted by $\mathbf{u}^+$), and the outward facing normal vector $\mathbf{n}$. Due to the discontinuous nature of the solution across the element interfaces, the flux is generally computed using exact~\citep{Toro1997_4} or approximate Riemann solvers such as the approach of \citet{Rusanov1962}.

The DG scheme possesses a highly useful property in that the evolution of the element-wise mean solution $\overline{\mathbf{u}}$, defined as 
\begin{equation}
    \overline{\mathbf{u}} = \frac{\int_{\Omega_k}\mathbf{u} (\mathbf{x})\ \mathrm{d}\mathbf{x}}{\int_{\Omega_k} \mathrm{d}\mathbf{x}},
\end{equation}
is dependent on only the interface flux, i.e.,
\begin{equation}
    \partial_t \overline{\mathbf{u}} = - \int_{\partial \Omega_k} \hat{\mathbf{F}}(\mathbf{u}^-, \mathbf{u}^+, \mathbf{n}) \ \mathrm{d}S.
\end{equation}
This yields a critical property in which the temporal update of the element-wise mean is equivalent to a convex combination of temporal updates of first-order Godunov schemes~\citep{Godunov1959}, which, under some relatively minor assumptions on the numerical scheme (e.g., CFL condition, appropriate Riemann solver, strong stability preserving time integration, appropriate quadrature, etc.), ensures that the element-wise mean preserves convex invariants of the system~\citep{Zhang2010}. We remark here that while the focus of this work is on explicit DG schemes for hyperbolic conservation laws without source terms, the proposed method may be directly extended to any system for which the DG scheme preserves invariants of on the element-wise mean, which allows for its application to implicit schemes~\citep{Qin2018}, mixed hyperbolic-parabolic systems~\citep{Zhang2017}, and certain hyperbolic systems with source terms which possess these properties on the element-wise mean~\citep{Zhang2011c}. This property of the element-wise mean plays a key role in many numerical stabilization methods for high-order DG schemes~\citep{Zhang2010,Zhang2011b,Zhang2017,Qin2018,Dzanic2022,Ching2023,Zhang2023,Dzanic2023b}, and the reader is referred to the aforementioned series of works stemming from \citet{Zhang2010} for a more in-depth overview on the sufficient conditions for this property. Furthermore, the proposed approach can be combined with any other approach which ensures that the element-wise mean satisfies a given set of physical constraints, including additional limiting/stabilization methods which satisfy the constraints on a set of discrete quadrature points or solution post-processing operations such as coarsening for adaptive mesh refinement. 
\section{Methodology}\label{sec:methodology}
In this section, the proposed limiting formulation for enforcing continuous constraints in DG schemes is first introduced. Afterwards, the numerical optimization procedure for performing the limiting and an overview of various possible constraint functionals are then presented. 

\subsection{Numerical formulation}\label{ssec:formulation}
The underlying limiting method of the proposed approach is the \textit{a posteriori} ``squeeze'' limiter of \citet{Zhang2010} which was introduced for the purposes of \emph{discretely} enforcing convex constraints in DG schemes. In this approach, the limiter acts to uniformly contract the high-order solution $\mathbf{u}(\mathbf{x})$ towards the element-wise mean $\overline{\mathbf{u}}$, given as
\begin{equation}\label{eq:limiter}
    \widehat{\mathbf{u}}(\mathbf{x}) = (1 - \alpha) \ux + \alpha \ubar = \ux + \alpha \left (\ubar - \ux \right),
\end{equation}
where $\widehat{\mathbf{u}}(\mathbf{x})$ is the limited solution and $\alpha \in [0,1]$ is a free parameter dictating the degree of limiting. In the discrete form, the limiting procedure ensures that the limited solution is bounds-preserving at some set of discrete nodal locations $\mathbf{x}_i \ \forall \ i \in S$, i.e.,
\begin{equation}
    g(\widehat{\mathbf{u}}(\mathbf{x}_i)) \geq 0 \ \forall \ i \in S,
\end{equation}
where $S$ may be the set of quadrature nodes and/or solution nodes for the scheme and $g$ refers to the constraint functionals in \cref{eq:constraints}. The limiter relies on the ability of DG schemes to preserve convex invariants on the element-wise mean as discussed in \cref{sec:preliminaries} (i.e., $g(\overline{\mathbf{u}}) \geq 0$), such that there always exists a value of $\alpha$ for which the constraints are satisfied by the limited solution. This limiting approach encompasses a variety of numerical stabilization methods, including slope-limited, maximum principle/positivity-preserving, and entropy-bounded schemes~\citep{Barth1989,Zhang2010,Zhang2011b,Dzanic2022,Ching2023}. Furthermore, it is highly advantageous for DG schemes as it is locally conservative, does not require a secondary low-order numerical scheme, is trivial to implement on general unstructured meshes, and can retain the high-order accuracy of the underlying DG scheme.

\begin{figure}[htbp!]
        \centering
        \adjustbox{width=0.4\linewidth,valign=b}{\input{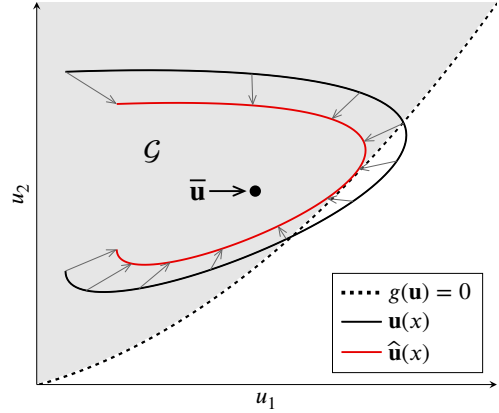}}
        \caption{\label{fig:limiting_schematic} Schematic of the limiting procedure for a vector-valued solution $\mathbf{u}(x) = [u_1(x), u_2(x)]^T$ showcasing the minimal limiting necessary for the limited solution $\widehat{\mathbf{u}}(x)$ to be continuously bounds-preserving. }
    \end{figure}

The goal of the proposed approach is to efficiently find some value of $\alpha$ such that the limited solution is \emph{continuously} bounds-preserving, i.e.,
\begin{equation}
    \ghat \geq 0 \ \forall \ \mathbf{x}.
\end{equation}
Without loss of generality, we consider only one constraint functional $\gu$ for some arbitrary element $\Omega_k$ in the mesh, but the analysis readily extends to sets of functionals and all elements in the mesh. For brevity, the notation $\forall\ \mathbf{x}$ is used to refer to all $\mathbf{x} \in \Omega_k$. A schematic of the goal of the proposed limiting approach is presented in \cref{fig:limiting_schematic}, showing the high-order solution contracting towards the element-wise mean to where the limited solution resides completely within the admissible set $\mathcal G$. In the ideal case, the limiting is performed such as to minimally modify the high-order solution while continuously enforcing constraints, which effectively corresponds to finding the limiting factor for which the limited solution is tangent to the zero contour of the constraint functional.

Since the goal is to now enforce constraints continuously via the limiting approach, it is typically necessary to find the continuous extrema of the solution (or some functional of the solution) to see to what extent the solution violates the constraints. Outside of some trivial cases, this requires the use of optimization techniques to find the continuous extrema in each element. However, this optimization problem must be carefully formulated to avoid incurring a substantial computational cost. A natural choice may be to find the minima of the constraint functional and apply limiting based on that nodal location. If we denote $g^*$ as the spatial minimum of the constraint functional, computed as
\begin{equation}
    g^* = \underset{\mathbf{x}}{\min}\ g(\mathbf{u}(\mathbf{x})), \quad \mathbf{x}^* = \underset{\mathbf{x}}{\arg \min}\ g(\mathbf{u}(\mathbf{x})),
\end{equation}
then it can easily be shown from the work of \citet{Zhang2010} that setting $\alpha$ as 
\begin{equation}\label{eq:discrete_alpha}
    \alpha = \max \left[ 0, \frac{g^*}{g^* - \gbar}\right]
\end{equation}
ensures that $g(\widehat{\mathbf{u}}(\mathbf{x}^*)) \geq 0$ under the standard assumptions. However, we claim here that limiting based on the minimum of the constraint functional \emph{does not ensure that the limited solution will be continuously bounds-preserving}, i.e., $g(\widehat{\mathbf{u}}(\mathbf{x}^*)) \geq 0 \nRightarrow \ghat \geq 0 \ \forall \ \mathbf{x}$.

The reason for this difficulty is that the minimum of the constraint functional is not necessarily the point which requires the most limiting. It can be seen that limiting per \cref{eq:discrete_alpha} is dependent on not only the minimum of the constraint functional but also the value of the constraint functional of the element-wise mean. In fact, $\alpha$ here behaves as the contraction factor necessary to ``squeeze'' the high-order solution into the admissible set, as shown in \cref{fig:limiting_schematic}. For any given spatial location (i.e., a point on the black phase-space curve in \cref{fig:limiting_schematic}), the limiting operation contracts the solution from this point along the line segment between this point and the element-wise mean. In this view, we can define an inner distance and an outer distance, referring to the length of the line segment residing inside and outside of the set, respectively. As the contraction is a convex combination of both the high-order solution point and the element-wise mean, the necessary level of contraction is not dependent on the outer distance but the \emph{ratio} of the outer and inner distances (i.e., the necessary limiting factor depends on not just the distance of the solution point to the boundary of the admissible set but also the distance of the element-wise mean to the boundary of the admissible set). If two arbitrary points are taken along the phase-space curve of the high-order solution (with their two associated line segments towards the element-wise mean), the point which resides closer to the admissible set may in fact require more limiting if the element-wise mean is very close to the boundary of the admissible set along its respective line segment.

\begin{figure}[htbp!]
        \centering
        \adjustbox{width=0.6\linewidth,valign=b}{\begin{tikzpicture}[spy using outlines={rectangle, height=3cm,width=2.3cm, magnification=3, connect spies}]
	\begin{axis}[name=plot1,
		axis line style={latex-latex},
	    axis x line=left,
        axis y line=left,
		xlabel={$u_1$},
    	xmin=-3, xmax=5,
    	ylabel={$u_2$},
    	ymin=-2, ymax=2,
        axis equal image,
        clip mode=individual,
        axis on top,
        label style={font=\footnotesize},
        tick label style={font=\footnotesize},
    	legend style={at={(0.97, 0.97)},anchor=north east,font=\footnotesize},
    	legend cell align={left},
    	style={font=\normalsize}]

        \draw[black!0, fill=gray!20] (0,0) circle (1);
        \node[black] at (0.2, 0.5) { $\mathcal G$};

    	\addplot[black, thick, dotted] coordinates {
            (-1.2246467991473532e-16, -1.0) (-0.06342391965656484, -0.9979866764718843) (-0.12659245357374938, -0.9919548128307953) (-0.18925124436041008, -0.9819286972627067) (-0.2511479871810793, -0.9679487013963562) (-0.31203344569848734, -0.9500711177409453) (-0.3716624556603276, -0.9283679330160726) (-0.4297949120891718, -0.9029265382866211) (-0.4861967361004687, -0.8738493770697849) (-0.5406408174555978, -0.8412535328311811) (-0.5929079290546404, -0.8052702575310586) (-0.6427876096865395, -0.7660444431189779) (-0.6900790114821119, -0.7237340381050702) (-0.7345917086575334, -0.6785094115571321) (-0.776146464291757, -0.6305526670845222) (-0.8145759520503357, -0.580056909571198) (-0.8497254299495146, -0.5272254676105022) (-0.8814533634475821, -0.4722710747726826) (-0.9096319953545184, -0.4154150130018863) (-0.9341478602651068, -0.3568862215918718) (-0.954902241444074, -0.2969203753282746) (-0.9718115683235418, -0.23575893550942695) (-0.984807753012208, -0.1736481776669303) (-0.9938384644612542, -0.11083819990101054) (-0.998867339183008, -0.04758191582374228) (-0.9998741276738751, 0.015865963834808153) (-0.9968547759519423, 0.07924995685678866) (-0.9898214418809327, 0.14231483827328534) (-0.9788024462147786, 0.20480666806519107) (-0.963842158559942, 0.266473813690035) (-0.9450008187146685, 0.3270679633174218) (-0.9223542941045814, 0.3863451256931288) (-0.8959937742913359, 0.44406661260577424) (-0.8660254037844385, 0.5000000000000003) (-0.8325698546347712, 0.5539200638661106) (-0.7957618405308319, 0.6056096871376669) (-0.7557495743542583, 0.6548607339452851) (-0.7126941713788628, 0.7014748877063214) (-0.6667690005162915, 0.7452644496757548) (-0.6181589862206051, 0.7860530947427876) (-0.5670598638627704, 0.8236765814298329) (-0.5136773915734061, 0.8579834132349772) (-0.4582265217274103, 0.8888354486549235) (-0.40093053540661355, 0.9161084574320697) (-0.34202014332566855, 0.9396926207859084) (-0.28173255684142945, 0.9594929736144975) (-0.2203105327865403, 0.9754297868854072) (-0.15800139597334953, 0.9874388886763944) (-0.09505604330418257, 0.9954719225730846) (-0.031727933498067504, 0.9994965423831852) (0.03172793349806785, 0.9994965423831851) (0.09505604330418327, 0.9954719225730846) (0.15800139597335022, 0.9874388886763943) (0.22031053278654067, 0.975429786885407) (0.2817325568414301, 0.9594929736144973) (0.3420201433256689, 0.9396926207859083) (0.40093053540661416, 0.9161084574320694) (0.4582265217274106, 0.8888354486549234) (0.5136773915734063, 0.8579834132349771) (0.5670598638627711, 0.8236765814298325) (0.6181589862206053, 0.7860530947427874) (0.666769000516292, 0.7452644496757543) (0.7126941713788629, 0.7014748877063212) (0.7557495743542587, 0.6548607339452847) (0.7957618405308323, 0.6056096871376664) (0.8325698546347714, 0.5539200638661103) (0.8660254037844388, 0.4999999999999997) (0.8959937742913359, 0.444066612605774) (0.9223542941045816, 0.3863451256931282) (0.9450008187146686, 0.3270679633174214) (0.9638421585599423, 0.2664738136900343) (0.9788024462147787, 0.20480666806519043) (0.9898214418809327, 0.14231483827328512) (0.9968547759519424, 0.079249956856788) (0.9998741276738751, 0.01586596383480771) (0.998867339183008, -0.047581915823742944) (0.9938384644612541, -0.1108381999010112) (0.984807753012208, -0.1736481776669303) (0.9718115683235415, -0.23575893550942761) (0.9549022414440739, -0.29692037532827503) (0.9341478602651065, -0.3568862215918724) (0.9096319953545182, -0.41541501300188677) (0.8814533634475816, -0.4722710747726834) (0.849725429949514, -0.5272254676105029) (0.8145759520503357, -0.580056909571198) (0.7761464642917565, -0.6305526670845228) (0.7345917086575331, -0.6785094115571324) (0.6900790114821116, -0.7237340381050705) (0.6427876096865391, -0.7660444431189782) (0.5929079290546397, -0.8052702575310591) (0.5406408174555974, -0.8412535328311813) (0.4861967361004687, -0.8738493770697849) (0.42979491208917103, -0.9029265382866215) (0.37166245566032713, -0.9283679330160728) (0.3120334456984865, -0.9500711177409455) (0.25114798718107884, -0.9679487013963564) (0.18925124436041008, -0.9819286972627067) (0.12659245357374893, -0.9919548128307953) (0.0634239196565644, -0.9979866764718844) (1.2246467991473532e-16, -1.0) };
        \addlegendentry{$g(\mathbf{u}) = 0$};

    	\addplot[black, thick] coordinates {
(-2.665, -1.0) (-2.659, -0.998) (-2.653, -0.992) (-2.647, -0.982) (-2.641, -0.968) (-2.634, -0.95) (-2.628, -0.928) (-2.62, -0.903) (-2.613, -0.874) (-2.605, -0.841) (-2.597, -0.805) (-2.589, -0.766) (-2.58, -0.724) (-2.571, -0.679) (-2.562, -0.631) (-2.552, -0.58) (-2.542, -0.527) (-2.531, -0.472) (-2.52, -0.415) (-2.508, -0.357) (-2.496, -0.297) (-2.484, -0.236) (-2.471, -0.174) (-2.457, -0.111) (-2.443, -0.048) (-2.428, 0.016) (-2.413, 0.079) (-2.397, 0.142) (-2.38, 0.205) (-2.363, 0.266) (-2.345, 0.327) (-2.326, 0.386) (-2.307, 0.444) (-2.287, 0.5) (-2.265, 0.554) (-2.243, 0.606) (-2.22, 0.655) (-2.197, 0.701) (-2.172, 0.745) (-2.146, 0.786) (-2.119, 0.824) (-2.091, 0.858) (-2.061, 0.889) (-2.031, 0.916) (-1.999, 0.94) (-1.966, 0.959) (-1.932, 0.975) (-1.896, 0.987) (-1.859, 0.995) (-1.82, 0.999) (-1.78, 0.999) (-1.738, 0.995) (-1.694, 0.987) (-1.648, 0.975) (-1.601, 0.959) (-1.551, 0.94) (-1.5, 0.916) (-1.446, 0.889) (-1.39, 0.858) (-1.332, 0.824) (-1.272, 0.786) (-1.209, 0.745) (-1.143, 0.701) (-1.075, 0.655) (-1.003, 0.606) (-0.929, 0.554) (-0.852, 0.5) (-0.772, 0.444) (-0.688, 0.386) (-0.601, 0.327) (-0.511, 0.266) (-0.416, 0.205) (-0.318, 0.142) (-0.216, 0.079) (-0.109, 0.016) (0.002, -0.048) (0.117, -0.111) (0.238, -0.174) (0.363, -0.236) (0.493, -0.297) (0.629, -0.357) (0.771, -0.415) (0.918, -0.472) (1.071, -0.527) (1.231, -0.58) (1.397, -0.631) (1.57, -0.679) (1.75, -0.724) (1.938, -0.766) (2.133, -0.805) (2.336, -0.841) (2.548, -0.874) (2.769, -0.903) (2.998, -0.928) (3.237, -0.95) (3.486, -0.968) (3.746, -0.982) (4.015, -0.992) (4.296, -0.998) (4.589, -1.0) 
        };
        \addlegendentry{$\mathbf{u}(x)$};

    	\addplot[black, thick, red!90!black] coordinates {
   (-1.547, -0.349) (-1.546, -0.348) (-1.544, -0.346) (-1.542, -0.342) (-1.539, -0.338) (-1.537, -0.332) (-1.535, -0.324) (-1.532, -0.315) (-1.53, -0.305) (-1.527, -0.294) (-1.524, -0.282) (-1.522, -0.269) (-1.519, -0.254) (-1.516, -0.239) (-1.512, -0.222) (-1.509, -0.205) (-1.505, -0.187) (-1.502, -0.168) (-1.498, -0.149) (-1.494, -0.129) (-1.49, -0.108) (-1.486, -0.087) (-1.481, -0.066) (-1.477, -0.044) (-1.472, -0.023) (-1.467, -0.001) (-1.461, 0.021) (-1.456, 0.042) (-1.45, 0.063) (-1.444, 0.085) (-1.438, 0.105) (-1.432, 0.126) (-1.425, 0.145) (-1.418, 0.164) (-1.411, 0.183) (-1.403, 0.201) (-1.396, 0.217) (-1.387, 0.233) (-1.379, 0.248) (-1.37, 0.262) (-1.361, 0.275) (-1.351, 0.287) (-1.341, 0.297) (-1.331, 0.307) (-1.32, 0.315) (-1.309, 0.322) (-1.297, 0.327) (-1.285, 0.331) (-1.272, 0.334) (-1.259, 0.335) (-1.245, 0.335) (-1.23, 0.334) (-1.215, 0.331) (-1.2, 0.327) (-1.184, 0.322) (-1.167, 0.315) (-1.149, 0.307) (-1.131, 0.297) (-1.112, 0.287) (-1.092, 0.275) (-1.071, 0.262) (-1.05, 0.248) (-1.027, 0.233) (-1.004, 0.217) (-0.979, 0.201) (-0.954, 0.183) (-0.928, 0.164) (-0.9, 0.145) (-0.872, 0.126) (-0.842, 0.105) (-0.811, 0.085) (-0.779, 0.063) (-0.745, 0.042) (-0.71, 0.021) (-0.673, -0.001) (-0.636, -0.023) (-0.596, -0.044) (-0.555, -0.066) (-0.512, -0.087) (-0.467, -0.108) (-0.421, -0.129) (-0.373, -0.149) (-0.322, -0.168) (-0.27, -0.187) (-0.215, -0.205) (-0.158, -0.222) (-0.099, -0.239) (-0.038, -0.254) (0.027, -0.269) (0.093, -0.282) (0.163, -0.294) (0.235, -0.305) (0.311, -0.315) (0.389, -0.324) (0.471, -0.332) (0.556, -0.338) (0.645, -0.342) (0.737, -0.346) (0.833, -0.348) (0.933, -0.349) 
        };
        \addlegendentry{$\widehat{\mathbf{u}}(x)$};

    	\addplot[color=black, style={thick}, mark=*, only marks, mark options={scale=0.5}] coordinates {
    	    (-0.9668, -0.01)};

        \draw[black,thin,->,>=angle 45] (-1.5, -1) -- (-0.9668, -0.01) node[midway,left,yshift=-2.5ex,xshift=-1.5ex]{$\overline{\mathbf{u}}$};

    	\addplot[color=black, style={thin}, mark=o, only marks, mark options={scale=0.5}] coordinates {
    	    (4.589, -1.0)};

        \draw[black,thin,->,>=angle 45] (4, -1.5) -- (4.589, -1.0) node[midway,left,yshift=-2.5ex,xshift=-1.5ex]{$g^*$};

	\end{axis}
\end{tikzpicture}}
        \caption{\label{fig:gstar_limiting} Schematic for an example system and constraints showing that limiting based on the extremum of the constraint functional ($g^*$) does \emph{not} guarantee that the limited solution is continuously bounds-preserving.}
    \end{figure}
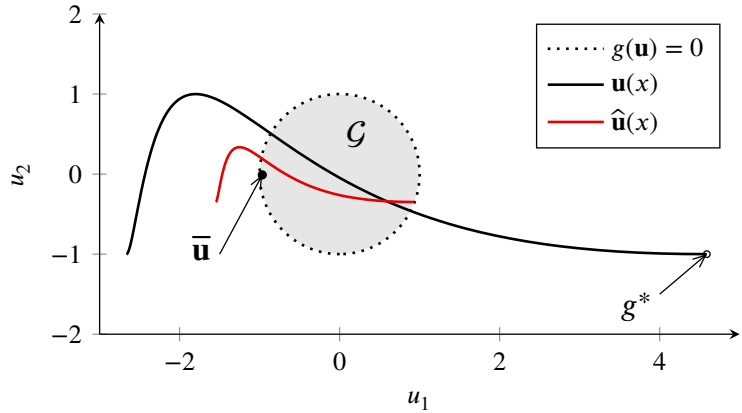

A simple example of this can be given for a vector-valued solution $\mathbf{u}(x) = [u_1(x), u_2(x)]^T$, where $u_1(x)$ and $u_2(x)$ refer to two components of an arbitrary system of equations, on the domain $\Omega = [-1, 1]$, with the solution and constraint functional taken as
\begin{subequations}
\begin{align}
    u_1(x) &= \exp (2x) - 2.8,\\
    u_2(x) &= \cos (\pi x),\\
    g(\mathbf{u}) &= 1 - u_1^2 - u_2^2.
\end{align}
\end{subequations}
The solution profiles and the constraint functional are arbitrarily chosen to illustrate the concept. \cref{fig:gstar_limiting} shows that for this example system, limiting based on the extremum of the constraint functional does in fact guarantee that the limited solution is bounds-preserving at that point, but it does \emph{not} guarantee that the solution is continuously bounds-preserving as $\exists x : \widehat{\mathbf{u}}(x) \notin \mathcal G$. In this scenario, the limiting procedure would have to be iterated until the constraints are satisfied continuously which would incur a large computational cost as the optimization procedure would have to be repeated for each iteration.

The primary novelty of this work is to introduce a \emph{modified} formulation for the constraint functionals such that the limiting approach can ensure a continuously bounds-preserving solution with only a single spatial minimization problem per element. To this end, we define these modified constraint functionals as 
\begin{equation}\label{eq:h}
    h(\mathbf{u}) = 
    \begin{cases}
    h^+(\mathbf{u}), \quad \text{if } \gu \geq 0, \\
    h^-(\mathbf{u}), \quad \text{else},
    \end{cases}
\end{equation}
where the sub-functionals are defined as 
\begin{equation}\label{eq:hplus}
    h^+(\mathbf{u}) = \frac{\gu}{\gbar}
\end{equation}
and
\begin{equation}\label{eq:hminus}
    h^-(\mathbf{u}) = \frac{\gu}{\gbar - \gu}.
\end{equation}
The purpose of these modifications as well as this piecewise formulation will be explained momentarily. We claim that if one instead optimizes over the modified functional, yielding
\begin{equation}\label{eq:hstar}
    h^* = \underset{\mathbf{x}}{\min}\ h(\mathbf{u}(\mathbf{x})), \quad \mathbf{x}^* = \underset{\mathbf{x}}{\arg \min}\ h(\mathbf{u}(\mathbf{x})),
\end{equation}
then by setting the limiting factor as 
\begin{equation}\label{eq:alpha}
    \alpha = \max(0, -h^*),
\end{equation}
the limited solution will be continuously bounds-preserving with only one optimization pass (i.e., only one calculation of the spatial minimization problem). From this formulation, it can be seen that $\alpha \in [0,1]$ as $\gbar$ is strictly non-negative.

The justification for this formulation is as follows. The sub-functional of interest for a bounds-violating solution is $h^-(\mathbf{u})$ as $h(\mathbf{u}) = h^-(\mathbf{u})$ if $g(\mathbf{u}) < 0$. It can be seen that this sub-functional mimics \cref{eq:discrete_alpha}. In fact, optimizing over this functional is effectively optimizing over the limiting factor (or at least a linear approximation thereof) necessary to ensure that the limited solution is bounds-preserving at that given spatial location. This transforms the optimization problem from finding where the constraints are violated the most to finding where the most limiting is required. As will be shown, finding the minimum of $h^-(\mathbf{u})$ and limiting based off of it is sufficient to ensure that the limited solution is continuously bounds-preserving. However, this sub-functional is not necessarily well-behaved due to its rational nature as its value may become undefined --- for scalar systems/constraints, this behavior is guaranteed in at least one location per element. To remedy this, $h^-(\mathbf{u})$ is augmented in a piecewise manner with $h^+(\mathbf{u})$ to yield a modified constraint functional $h(\mathbf{u})$ with the following properties:
\begin{itemize}
    \item Well-defined for all $\mathbf{x}$ to yield numerically stable behavior.
    \item Of the same sign as $g(\mathbf{u})$ for consistency (i.e., $h(\mathbf{u}) < 0$ for a bounds-violating solution and $h(\mathbf{u}) \geq 0$ for a bounds-preserving solution).
    \item $C^1$ continuous with respect to $\mathbf{x}$ to be amenable to optimization algorithms. 
\end{itemize}
These desired properties motivate the definition of $h^+(\mathbf{u})$ in \cref{eq:hplus}, such that the modified formulation of $h(\mathbf{u})$ in \cref{eq:h} achieves much more desirable numerical behavior, namely that it and its derivative are always well-defined (if $\gu$ and its derivative are well-defined).

With this proposed formulation in mind, the proofs of certain key properties of the proposed scheme are now presented. The following proofs rely on two key assumptions, the quasiconcavity of the constraint functional and the bounds-preserving property of the element-wise mean. Furthermore, we assume a strict inequality on the element-wise mean (i.e., $g(\overline{\mathbf{u}}) > 0$) for brevity, with the limiting procedure for the simple edge-case of $g(\overline{\mathbf{u}}) = 0$ shown later.

\begin{theorem}[Bounds preservation]\label{thm:cbp}
    The limited solution as given by \cref{eq:limiter} is continuously bounds-preserving (i.e., $\ghat \geq 0 \ \forall \ \mathbf{x}$) if the limiting factor is computed as \cref{eq:alpha}.
\end{theorem}

\newproof{pot_cbp}{Proof of \cref{thm:cbp}}
\begin{pot_cbp}
Per \cref{eq:alpha}, the limiting factor may either be $\alpha = 0$ or $\alpha = -h^*$. For the former case, $h^* \geq 0$, which gives the equivalency $h^* = h^+(\mathbf{u}(\mathbf{x}^*))$. Given the positivity of the denominator in \cref{eq:hplus}, it can be seen that 
\begin{equation}
    \gux \geq \gstar \geq 0 \ \forall \ \mathbf{x}.
\end{equation}
It then follows that for $\alpha = 0$, 
\begin{equation}
    \ghat = \gux \geq 0 \ \forall \ \mathbf{x}.
\end{equation}

For the latter case, $h^* \leq 0$, which gives the equivalency $h^* = h^-(\mathbf{u}(\mathbf{x}^*))$. Therefore,
\begin{equation}
    \uhat =  \ux - \frac{\gstar}{\gbar - \gstar} \left (\ubar - \ux \right) = (1 - \alpha) \ux + \alpha \ubar ,
\end{equation}
where
\begin{equation}
    \alpha = -\frac{\gstar}{\gbar - \gstar}.
\end{equation}
Furthermore, for $h^* \leq 0$, it can be seen that $\gstar \leq 0$. Due to Jensen's inequality, 
\begin{equation}
    \ghat = g \left ( (1 - \alpha) \ux + \alpha \ubar \right) \geq (1 - \alpha) \gux + \alpha \gbar  =  \gux + \alpha \left ( \gbar - \gux \right).
\end{equation}
As such, it is sufficient to show that 
\begin{equation}\label{eq:gineq}
    \gux + \alpha \left ( \gbar - \gux \right) = \gux - \frac{\gstar}{\gbar - \gstar} \left (\gbar - \gux\right)\geq 0\  \forall\  \mathbf{x}.
\end{equation}
Note that only the case of $\gux < 0$ has to be proven for as $\ghat\geq 0$ for any $\alpha \in [0,1]$ if $\gux \geq 0$ (i.e., any convex combination of bounds-preserving states is bounds-preserving due to the convexity of $\mathcal G$). In this case, the quantity $\gbar - \gux$ is positive, such that \cref{eq:gineq} can be rearranged to
\begin{equation}
    \frac{\gux}{\gbar - \gux} - \frac{\gstar}{\gbar - \gstar} \geq 0 \  \forall\  \mathbf{x}.
\end{equation}
This is simply the condition
\begin{equation}\label{eq:hstarcond}
    \hux - \hstar \geq 0\ \forall \ \mathbf{x},
\end{equation}
which is satisfied as, by definition, $\hux \geq \hstar$ for all $\mathbf{x}$. Consequently, 
\begin{equation}
    \ghat \geq 0  \ \forall \ \mathbf{x}.
\end{equation}
\end{pot_cbp}

\begin{corollary}\label{thm:nolimiting}
From \cref{thm:cbp}, the limiting procedure does not modify the high-order solution $\mathbf{u}(\mathbf{x})$ if the high-order solution is continuously bounds-preserving. 
\end{corollary}

\begin{theorem}[Continuity]\label{thm:cont}
    The modified constraint functional in \cref{eq:h} is $C^1$ continuous with respect to $\mathbf{x}$ if $\gux$ is $C^1$ continuous with respect to $\mathbf{x}$.
\end{theorem}
\newproof{pot_cont}{Proof of \cref{thm:cont}}
\begin{pot_cont}
To show $C^1$ continuity for the modified constraint functional in \cref{eq:h}, it is sufficient to show that $h^+(\ux) \in C^1$ for $\gux \geq 0$, $h^-(\ux) \in C^1$ for $\gux \leq 0$, and $h^{+'}(\ux) = h^{-'}(\ux)$ for $\gux = 0$. Since $\gbar$ is constant across an element, $g'(\ubar) = 0$. Therefore,
\begin{equation}
    h^{+'}(\ux) = \frac{g'(\ux)}{\gbar} \quad \text{and} \quad h^{-'}(\ux) = \frac{\gbar g'(\ux)}{(\gbar - \gux)^2}.
\end{equation}
Due to the positivity of $\gbar$, it can be seen that if $\gux \in C^1$, then $h^+(\ux)$ and $h^-(\ux)$ are also $C^1$ continuous across their respective domains. Furthermore, for $\gux = 0$, 
\begin{equation}
    h^{+'}(\ux) = h^{-'}(\ux).
\end{equation}

\end{pot_cont}

\begin{lemma}\label{thm:lowerbound}
The limiting factor required for continuous bounds-preservation (computed by \cref{eq:alpha}) is bounded from below by the limiting factor required for discrete bounds-preservation (computed by \cref{eq:discrete_alpha}).
\end{lemma}

\bigskip
\bigskip

With these key properties of the proposed scheme shown, some further details are presented concerning edge cases, higher-order limiting, and multiple constraint functionals. First, the edge case of $g(\ubar) \to 0$ can result in some numerical difficulties, particularly as $\gu \to 0$. In this scenario, the limiting factor is  explicitly set to $\alpha = 1$ for $g(\ubar) < \epsilon$, where $\epsilon = 10^{-12}$ is some small positive number on the order of machine precision, as it can be seen that $\alpha \to 1$ as $g(\ubar) \to 0$. 

Second, the denominator in \cref{eq:hminus} effectively acts as a linear approximation of the convex limiting factor necessary to ensure that the limited solution is continuously bounds-preserving, approximating \emph{from above} the solution of 
\begin{equation*}
    g \left ( (1 - \alpha) \ux + \alpha \ubar \right) = 0.
\end{equation*}
For constraint functionals that are linear with respect to the solution variables, the resulting limiting factor is exact, i.e., it is the smallest necessary limiting factor. However, for arbitrary nonlinear quasiconcave constraint functionals, the resulting limiting factor is only sufficient (i.e., the above equality becomes a positive inequality), and there may be a smaller limiting factor for which the limited solution is still continuously bounds-preserving. It is possible to formulate the modified constraint functionals such that the resulting limiting factor is exact for higher-order constraint functionals (e.g., quadratic constraints such as pressure/internal energy in gas dynamics) or to incorporate a nonlinear limiting procedure in the spatial optimization problem to better approximate the minimum necessary limiting factor for arbitrary nonlinear constraints. This may be achieved by replacing the subfunctional $h^-(\mathbf{u})$ with the exact (or sufficiently approximate) solution of the intersection problem,
\begin{equation*}
    h^-(\mathbf{u}) = -\alpha^* : g \left ( (1 - \alpha^*) \mathbf{u} + \alpha^* \ubar \right) = 0,
\end{equation*}
which enforces an equality on the constraint functional of the limited solution (i.e., a necessary condition) instead of an inequality (i.e., a sufficient condition). However, this is left as a topic of future work. We remark here that this effect is only present for solutions where the nonlinear constraints are violated -- per \cref{thm:nolimiting}, for a solution that is already continuously bounds-preserving, no limiting whatsoever is applied. 

Finally, the use of multiple constraint functionals provides a choice in how to approach the limiting procedure. The first method is to individually compute $\alpha$ for each of the $N_c$ constraints and take the maximum value. The second method is to sequentially perform the limiting procedure, where for each constraint, the limited solution as computed from the previous constraint is used as the ``unlimited'' solution for the next constraint. We posit here that the first method will likely yield more performant behavior (i.e., less numerical dissipation) for general nonlinear constraints whereas the two methods would be identical for linear constraints, the reasoning being closely related to the previous discussion on the sufficient versus necessary limiting factor. For nonlinear constraints, computing $\alpha$ by the proposed approach may yield an overprediction in the necessary limiting factor. As such, computing $\alpha$ sequentially may compound this overprediction, causing a higher degree of limiting than computing $\alpha$ individually. However, in some scenarios it may be necessary to compute the limiting sequentially, such as when the convexity of $\mathcal G$ or the well-definedness of some constraint functionals rely on other constraints being satisfied. An example of this can be seen in the gas dynamics equations, where the pressure/internal energy functional is quasiconcave only if the density is positive and the entropy functional is typically only well-defined if the density and pressure are positive. In this case, it is necessary to perform limiting sequentially, first on density, then pressure, then entropy.

\subsection{Optimization procedure}
The proposed limiting approach relies on finding the continuous minimum of the modified constraint functional within the element. For general applications, this minimum does not have a closed-form expression. As such, it must be found numerically using standard optimization techniques. Note that irrespective of the properties of the constraint functional itself, this is generally a global (but bounded) non-convex optimization problem as the solution itself can be and often is spatially non-convex. However, certain properties of the limiting procedure as well as judicious choices for the numerical formulation reduce this complexity. 

\paragraph{Bounded optimization.} While the optimization problem itself is global, it is bounded to a fixed domain (i.e., the element). As there is no inherent requirement to perform the search for the minimum of the constraint functional in physical space, the optimization can instead be performed in reference space. This greatly simplifies the implementation of bounded optimization schemes as the elements reduce to much simpler (convex) representations in reference space, particularly so for curved elements. The enforcement of search bounds then becomes easier to implement as it is then trivial to evaluate if the search coordinate is within the element bounds (e.g., the reference quadrilateral is simply $[-1, 1]^2$), and bounded optimization methods such as Lagrange multiplier-type approaches can be readily used as the element boundaries in reference space reduce to planar faces with fixed normals. Furthermore, this allows for the user to implement the optimization procedure only for a given reference element without having to consider the intricacies of performing bounded optimization on arbitrarily shaped/curved elements. 

\paragraph{Continuous evaluation.} Since the optimization problem requires repeatedly evaluating the solution at arbitrary locations, the standard nodal representation of the solution is not well-suited. Instead, the solution may be transformed from its nodal representation to a modal representation, given as 
\begin{equation}\label{eq:modal}
        \mathbf{u}(\hat{\mathbf{x}}) = \sum_{i=1}^{n_s} \mathbf{u}_{i}\phi_i(\hat{\mathbf{x}}) = \sum_{i=1}^{n_s} \widetilde{\mathbf{u}}_i\psi_i(\hat{\mathbf{x}}),
\end{equation}
where $\widetilde{\mathbf{u}}_i$ are a set of modal basis coefficients and $\psi_i(\hat{\mathbf{x}})$ are their respective basis functions. Note that $\hat{\mathbf{x}}$ is used instead of $\mathbf{x}$ to represent that the basis is with respect to the reference space instead of the physical space. This formulation is typically much more amenable to continuous evaluation. In particular, we use the monomial (power) basis to represent the solution, given as  
\begin{equation}\label{eq:monomial}
    \psi_i(\hat{\mathbf{x}}) = \hat{\mathbf{x}}^{\mathbf{q}_i} = \prod_{j=1}^{d} \hat{\mathbf{x}}_{j}^{\mathbf{q}_{ij}},
\end{equation}
which makes it trivial to evaluate the solution anywhere within the element. While the transformation to the monomial basis may be ill-conditioned in some scenarios, this is typically avoided if the solution is not of very high order. However, in such cases, one may use other bases such as orthogonal polynomials which may possess more robust numerical behavior. 

\paragraph{Initialization.}
The choice of initial guess for the minimization problem can play a large role in non-convex optimization methods. Even if the optimization landscape has many local minima, convex optimization techniques can still converge to the global minimum if the initial guess on a bounded domain is sampled sufficiently enough such that the initial guess falls within a locally convex region near the global minimum. Fortunately, for DG schemes, the solution is interpolated to the volume and surface quadrature nodes at each temporal step. As such, these quadrature nodes form a suitable set of locations to initially evaluate solution at. In fact, it can be shown through B\'ezout's theorem that for simpler cases such as linear constraint functionals the number of strict local minima is substantially fewer than the number of quadrature points~\citep{Durfee1993}.

\bigskip

With these points in mind, the proposed optimization procedure is presented, first with an introduction and explanation of the various features of the procedure followed by an overview of the overall optimization method. We remark here that since the proof in \cref{thm:cbp} is contingent on being able to find the global minimum of the modified constraint functional in the element (or at least a lower bound thereof), if one wants to truly \emph{guarantee} that limited solution is continuously bounds-preserving, global non-convex optimization techniques on bounded convex domains should be used. For this, branch-and-bound optimization methods~\citep{Lawler1966} would likely be a good option as the monomial basis allows for simple --- albeit somewhat loose --- estimates of upper/lower bounds of the solution on convex subdomains which, by the quasiconcavity of $\gu$, allows for estimates on the lower bound on the constraint functional necessary for branch-and-bound methods. Stricter bounds can be placed by instead using the convex hull property of the Bernstein basis and are often used to provide bounds for nonlinear functionals of the solution on convex subdomains of the element~\citep{Johnen2013}. Alternatively, one may use simulated annealing-type methods~\citep{Bertsimas1993} as it will be shown later that finding an approximate global minimum instead of a precise local minimum is sufficient. 

However, while the computational cost of these bounded global optimization methods is not intractable, particularly so on GPU computing architectures for which the compute-heavy, element-local operations of the methods are well-suited, their algorithmic complexity motivates the use of more straightforward optimization methods. To this end, we instead utilize convex optimization techniques, namely Newton--Raphson and gradient descent methods, and rely on the fact that the solution is sampled sufficiently finely on the element such that it is very likely that one of the sampled points lies within a locally convex region near the global minimum. As this approach does not necessarily ensure that it will find the \emph{global} minimum, it does however lose the mathematical property that the solution is \emph{guaranteed} to be continuously bounds-preserving for general nonlinear constraint functionals, but in the numerical experiments soon to be shown, this approach always yielded a continuously bounds-preserving limited solution. 

First, an iterative method must be chosen for optimization. The primary underlying optimization method in this work is the Newton--Raphson scheme, chosen due to its high rate of convergence, with the gradient descent method used as a secondary optimization method under edge cases where the Newton--Raphson scheme is not well-suited. Using the solution interpolated to the surface and volume quadrature nodes, denoted as $\mathbf{x}^q$, the initial guess is taken as 
\begin{equation}
    \mathbf{x}_0 = \underset{\mathbf{x}_i \in \mathbf{x}^q}{\arg \min}\ h(\mathbf{u}(\mathbf{x}^i)).
\end{equation}
Each iteration of the Newton--Raphson scheme can be given as
\begin{equation}
    \mathbf{x}_n = \mathbf{x}_{n-1} - \mathbf{H}_{n-1}^{-1} \mathbf{J}_{n-1},
\end{equation}
where $\mathbf{H}_{n-1}$ and $\mathbf{J}_{n-1}$ represent the Hessian and Jacobian, respectively, of the cost functional $h(\mathbf{u}(\mathbf{x}))$ evaluated at iteration $n-1$. While these quantities may be evaluated analytically using the chain rule, for memory bandwidth-bound computing architectures such as GPUs (which are primarily used in this work), it is generally much more efficient to evaluate them numerically. As such, a centered difference approximation is used, e.g. for a two-dimensional domain, the Jacobian and Hessian can be computed as
\begin{equation}\label{eq:jacobian}
    \mathbf{J}_n = \left[\frac{h(\mathbf{u}(x_n+\delta x, y_n) - h(\mathbf{u}(x_n-\delta x, y_n)}{2\delta x}, \frac{h(\mathbf{u}(x_n, y_n+\delta x) - h(\mathbf{u}(x_n, y_n-\delta x)}{2\delta x}\right]^T
\end{equation}
and
\begin{equation}\label{eq:hessian}
    \mathbf{H} = \begin{bmatrix}
                \mathbf{H}_{11} & \mathbf{H}_{12} \\
                \mathbf{H}_{21} & \mathbf{H}_{22} \\
                \end{bmatrix},
\end{equation}
where
\begin{subequations}
\begin{align}
    \mathbf{H}_{11} & =
                    \frac{h(\mathbf{u}(x_n+\delta x, y_n) - 2 h(\mathbf{u}(x_n, y_n) + h(\mathbf{u}(x_n-\delta x, y_n)}{\delta x^2} \\
    \mathbf{H}_{12} &= \mathbf{H}_{21} = \frac{h(\mathbf{u}(x_n+\delta x, y_n+\delta x) - h(\mathbf{u}(x_n-\delta x, y_n+\delta x) - h(\mathbf{u}(x_n+\delta x, y_n-\delta x) + h(\mathbf{u}(x_n-\delta x, y_n-\delta x)}{4\delta x^2} \\
    \mathbf{H}_{22} &=  \frac{h(\mathbf{u}(x_n, y_n+\delta x) - 2 h(\mathbf{u}(x_n, y_n) + h(\mathbf{u}(x_n, y_n-\delta x)}{\delta x^2}.
\end{align}
\end{subequations}
The step size is taken as $\delta x = 10^{-4}$ to yield a minimum denominator of approximately the square root of the FP64 machine precision error~\citep{Dennis1996}. It will later be shown that the possible convergence issues due to using the approximate Jacobian/Hessian are avoided as it will not necessary to compute the exact minimum. 

For the gradient descent method, each iteration is simply computed as
\begin{equation}
    \mathbf{x}_n = \mathbf{x}_{n-1} - \beta \frac{\mathbf{J}_{n-1}}{\|\mathbf{J}_{n-1}\|_2},
\end{equation}
where $\beta$ is a free parameter dictating the step size of the iteration. This step size is computed adaptively using backtracking line search with the initial step length taken as $\beta_0 = 2/(p+1)$ as this is the approximate average distance between nodal quadrature points in reference space. At each inner iteration of the line search, corresponding to one outer iteration of the gradient descent method, the step size is halved (i.e., $\beta_{i+1} = \beta_i/2$), and the process is repeated until either a maximum of five line search iterations have completed or the Armijo--Goldstein~\citep{Armijo1966} condition is satisfied, i.e.,
\begin{equation}
    h(\mathbf{u}(\mathbf{x}_n)) \leq h(\mathbf{u}(\mathbf{x}_{n-1})) - c\beta\|\mathbf{J}_{n-1}\|_2,
\end{equation}
where the control parameter is taken as $c = 0.5$ as per \citet{Armijo1966}. This line search also has the benefit of helping prevent the optimizer from getting stuck in local minima. The gradient descent scheme is used as a fallback for when the Hessian is not invertible or if the functional is not locally convex. Therefore, the Newton--Raphson method is used when both $det(\mathbf{H}) > \epsilon$ and $\mathbf{H}$ is positive semi-definite, where $\epsilon$ is the same tolerance used for the edge case of $\gbar \to 0$. Otherwise, the gradient descent method is used. 

Next, since the Newton--Raphson and gradient descent methods do not naturally ensure that the spatial location at each iteration will remain within the element, it is necessary to apply some techniques for bounded optimization. The overarching theoretical framework for this approach derives from Lagrange multiplier methods --- or, more aptly, the Karush--Kuhn--Tucker conditions --- but the simplicity of the bounding domain allows for its straightforward implementation as simple linear projections. We assume here that the reference element is convex and its boundary can be represented by a finite set of straight edges/planar faces. The iterative update $\mathbf{x}_n \to \mathbf{x}_{x+1}$, equivalently represented in terms of a search direction $\Delta \mathbf{x}_n$ as $\mathbf{x}_{x+1} = \mathbf{x}_n + \Delta \mathbf{x}_n$, is modified according to one of three conditions.
\begin{enumerate}
    \item If $\mathbf{x}_{n+1} \in \Omega_k$, i.e., the iterative update moves strictly within the element, leave $\mathbf{x}_{n+1}$ as is. 
    \item If $\mathbf{x}_{n+1} \notin \Omega_k$ and $\mathbf{x}_{n} \notin \partial \Omega_k$, i.e., the iterative update moves from within the element (but not on the boundary) to outside of the element, project the update \emph{distance} onto the element boundary \emph{along the search direction} by solving for $\eta \in [0,1]$ such that $\mathbf{x}_{n+1} = \mathbf{x}_{n} + \eta \Delta \mathbf{x}_n \in \partial \Omega_k$.
    \item If $\mathbf{x}_{n+1} \notin \Omega_k$ and $\mathbf{x}_{n} \in \partial \Omega_k$, i.e., the iterative update moves from the element boundary to outside of the element, project the update \emph{direction} onto the element boundary \emph{parallel to the boundary} as $\mathbf{x}_{n+1} = \mathbf{x}_{n} + (1 - \mathbf{n})\odot\Delta \mathbf{x}_n \in \partial \Omega_k$, where $\mathbf{n}$ is the normal direction of the particular boundary edge/surface that $\mathbf{x}_{n}$ is currently on. While this ensures that the iterative update does not intersect the particular boundary edge/surface that $\mathbf{x}_{n}$ is currently on, it does not ensure that it does not intersect any of the other boundary edges/surfaces. As such, project as per (2) if necessary.
\end{enumerate}
It can be seen here that there are one of two possibilities in this optimization process. The optimizer will either converge towards a local minimum within the element, or the optimizer will converge towards a local minimum along the boundary of the element where the gradient of the constraint functional is parallel to the normal direction of the boundary (i.e., the Lagrange multiplier condition). Due to the fact that the boundary of reference elements can be represented as a small number of straight edges/planar faces, the computations involved in these projections at most simply reduce to computing the intersection of a line with other lines or planes. For the line search in the gradient descent method, this projection is performed only at the initial step $\Delta \mathbf{x} = -\beta_0 \mathbf{J}$ as every subsequent step is guaranteed to be within the element bounds. Algorithmic implementation details of these projection steps are presented in \cref{app:algs}.

Finally, the use of iterative methods requires applying some small convergence tolerance to the optimization process. However, this can cause issues in the limiting method as the bounds-preserving properties shown in \cref{thm:cbp} are contingent on truly finding the minimum of $\hu$ as per \cref{eq:hstarcond}. If one instead limits based on $h^* + \delta h$, where $\delta h$ is some positive value smaller than the convergence tolerance, it is possible that the limited solution may violate the bounds on the order of $\delta h$. In most cases, this is not acceptable as these bounds must be strictly enforced to ensure the robustness of the numerical scheme, such that one would have to converge to effectively machine precision levels to ensure a continuously bounds-preserving scheme which is computationally intractable. Instead, we propose a final step in the optimization process to alleviate this constraint and allow for reasonable values of the convergence tolerance. Note that by \cref{eq:hstarcond}, one does not necessarily need to find the minimum of $\hu$ per se, simply a lower bound for its minimum. This corresponds to the observation that one may not limit marginally less than what is necessary to ensure a continuously bounds-preserving scheme, but limiting marginally more than necessary is perfectly acceptable. 
As such, the final step in the optimization approach is to extrapolate a lower bound for $\hu$ given the information derived from the iterative method. If we denote $\hstar$ as the value of the true minimum of $\hu$, we instead seek some value $h^{**} \leq \hstar$ that will instead be used for limiting as per \cref{eq:alpha}. 

For this, we utilize some convergence properties of Newton--Raphson/backtracking gradient descent methods, where we seek some correction $\Delta h < 0$, such that 
\begin{equation}
    h^n - \Delta h \leq h^*. 
\end{equation}
If we assume that the optimization process after $n$ iterations is within some locally convex region of $h(\mathbf{u}(\mathbf{x})$ in the interior (or convex subregion on the face/edge), then the magnitude of the Jacobian is strictly decreasing across the line segment $[\mathbf{x}^n, \mathbf{x}^*]$, with $\mathbf{J}(\mathbf{x}^*) = \mathbf{0}$, i.e.,
\begin{equation}
    \left \|\mathbf{J}^n\right \|_2 \geq \underset{x \in [x^n, x^*]}{\max}\ \left \|\mathbf{J}(\mathbf{x})\right\|_2. 
\end{equation}
Note here that if the $\mathbf{x}^n$ resides on the face/edge of the element, the Jacobian here refers to the projected Jacobian as per step (3) above. Furthermore, error estimates for the chosen optimization algorithms under the given assumptions yield~\citep{More1982, Dunn1980}
\begin{equation}
    \left \|\mathbf{x}^* - \mathbf{x}^n\right\|_2 \leq \left \|\mathbf{x}^n - \mathbf{x}^{n-1}\right\|_2 = \left \|\Delta \mathbf{x}_{n-1}\right\|_2.
\end{equation}
The combination of these two properties yields the relation
\begin{equation}
    \left | h^n - h^* \right | \leq \left \|\mathbf{J}^n\right\|_2 \left \|\Delta \mathbf{x}_{n-1}\right\|_2.
\end{equation}
Therefore, we set
\begin{equation}
    \Delta h = \left \|\mathbf{J}^n\right\|_2 \left \|\Delta \mathbf{x}_{n-1}\right\|_2
\end{equation}
and 
\begin{equation}
    h^{**} = \max \left(-1, h^{n} - \Delta h \right),
\end{equation}
which satisfies the condition $h^{**} \leq \hstar$. A maximum is taken to ensure that $\alpha \in [0,1]$. A simple schematic representation of this tolerance correction is shown in \cref{fig:opt_ext}.

    \begin{figure}[htbp!]
        \centering
        \adjustbox{width=0.4\linewidth,valign=b}{\begin{tikzpicture}[spy using outlines={rectangle, height=3cm,width=2.3cm, magnification=3, connect spies}]
	\begin{axis}[name=plot1,
		axis line style={latex-latex},
	    axis x line=left,
        axis y line=left,
		xlabel={$x$},
    	xmin=-0.5, xmax=1,
    	ylabel={$h(\mathbf{u})$},
    	ymin=-0.5, ymax=1,
    	ticks = none,
        clip mode=individual,
        axis on top,
    	legend style={at={(0.03, 0.97)},anchor=north west},
    	legend cell align={left},
    	style={font=\normalsize}]

        \addplot [domain=-1:1, samples=101,unbounded coords=jump, color=black, style={thick}]{1.5*x^3 + x^2 -  x};
        \addlegendentry{$g(\mathbf{u})$};
        
    	\addplot[color=black, style={thick}, mark=*, only marks, mark options={scale=0.7}] coordinates {
    	    (0.3, -0.17)};
         
    	\addplot[color=black, style={thick}, mark=*, only marks, mark options={scale=0.7}] coordinates {
    	    (0.1, -0.0885)};
         
    	\addplot[color=black, style={thick}, mark=*, only marks, mark options={scale=0.7}] coordinates {
    	    (-0.3, 0.35)};
         
    	\addplot[color=red!90!black, style={thick}] coordinates {
    	    (0.1, -0.0885)
    	    (0.5, -0.35)};
    	\addplot[color=red!90!black, style={thick}, mark=*, only marks, mark options={scale=0.7}] coordinates {
    	    (0.5, -0.35)};

        \dimline[extension start length=1ex, extension end length=1ex] {(-0.3, -.23)}{(0.1, -.23)}{\small$\Delta x$};
        
        \dimline[extension start length=1ex, extension end length=1ex] {(0.5, -.4)}{(0.1, -.4)}{\small$\Delta x$};
        
        \dimline[extension start length=0ex, extension end length=0ex] {(0.7, -0.35)}{(0.7, -0.085)}{\small$\Delta h$};
         
    	\addplot[color=red!90!black, style={thick, dotted}] coordinates {
    	    (0.1, -0.0885)
    	    (0.1, -0.35)
    	    (0.8, -0.35)};
         
    	\addplot[color=gray, style={thick, dotted}] coordinates {
    	    (-0.3, 0.35)
    	    (-0.3, -0.5)};
    	\addplot[color=gray, style={thick, dotted}] coordinates {
    	    (0.1, -0.0885)
    	    (0.1, -0.5)};
    	\addplot[color=gray, style={thick, dotted}] coordinates {
    	    (0.1, -0.0885)
    	    (0.8, -0.0885)};

        \draw[black,thin,->,>=angle 45] (-0.3, 0.5) -- (-0.3, 0.37) node[midway,midway,yshift=3.2ex,xshift=1.7ex]{$x^{n-1}$};
        
        \draw[black,thin,->,>=angle 45] (0.1, 0.04) -- (0.1, -.07) node[midway,midway,yshift=3ex,xshift=0.65ex]{$x^{n}$};
        
        \draw[black,thin,->,>=angle 45] (0.3, -0.02) -- (0.3, -0.15) node[midway,midway,yshift=3ex,xshift=0.65ex]{$x^{*}$};
        
        \node at (0.85, -0.0885) {$h^n$};
        \node at (0.87, -0.35) {\textcolor{red!90!black}{$h^{**}$}};
        \node at (0.25, -0.27) {\textcolor{red!90!black}{$J^n$}};
                
	\end{axis}
\end{tikzpicture}}
        \caption{\label{fig:opt_ext} Schematic of the tolerance correction $\Delta h$ computed at the end of the optimization process.}
    \end{figure}

With the proposed approach introduced, algorithmic implementation details are presented in \cref{app:algs}, including an overview of the limiting approach and the bounded optimization process. For brevity, we present the bounded optimization approach in terms of a single boundary face/normal vector, but the approach readily extends to the multiple faces required for multi-dimensional elements.

\subsection{Constraints}
For generality, the proposed approach was introduced with respect to an arbitrary quasiconcave functional $g(\mathbf{u})$. Here, we present some examples of common constraints encountered in the simulation of convection dominated flows, with more in-depth details of constraints relating to particular systems later presented in \cref{sec:results}. The first example is a minimum principle, written with respect to some functional of the solution $\mu(\mathbf{u})$ and some minimum bound $c$ as
\begin{equation}
    \mu(\mathbf{u}) \geq c.
\end{equation}
It is common to enforce positivity of some functional of the solution (e.g., density and pressure in gas dynamics), which can be expressed by setting $c = 0$. Alternatively, a more general minimal principle can be enforced by setting $c = \mu_{\min}$, where $\mu_{\min}$ is some minimum bound on the functional (e.g., entropy). This can be enforced by simply setting the constraint functional as 
\begin{equation}
    g(\mathbf{u}) = \mu(\mathbf{u}) - c.
\end{equation}

In other cases, the chosen functional must be bounded from above as well as from below, i.e.,
\begin{equation}
    \mu(\mathbf{u}) \in [a, b],
\end{equation}
where $a$ and $b$ ($b > a$) are some minimum/maximum bounds, respectively. These bounds are common in the simulation of scalar systems which obey a strict minimum/maximum principle and can be enforced through a pair of constraint functionals as
\begin{equation}
    g_1(\mathbf{u}) = \mu(\mathbf{u}) - a, \quad  g_2(\mathbf{u}) = b - \mu(\mathbf{u}).
\end{equation}
Depending on the constraints and system in question, the order which these constraints are enforced may be of importance (as discussed in \cref{ssec:formulation}). 

As the constraints are enforced independently on an element-wise basis, the proposed scheme does not differentiate whether the given bounds (e.g., the values of $a$, $b$, and $c$) are global or local. In some cases, the bounds may depend on the values of the solution in a local domain of influence~\citep{Dzanic2022}, or they may be fixed throughout the entire domain. In this work, we consider only constraints with global bounds, but the extension to arbitrary local bounds is trivial as long as the bounds are constant across the element. 

\section{Implementation}\label{sec:implementation}
The proposed approach was implemented using the DG scheme recovered by the flux reconstruction (FR) method of \citet{Huynh2007} within PyFR~\citep{Witherden2014}, a high-order unstructured FR solver that can target massively-parallel CPU and GPU computing architectures. Experiments were performed on scalar transport equations and the gas dynamics equations in one and two dimensions on both CPUs and GPUs. At the interfaces, upwinding was performed for scalar equations whereas the Rusanov~\citep{Rusanov1962} Riemann solver with the Davis~\citep{Davis1988} wavespeed estimate was used for the gas dynamics equations. The FR solution/flux points, corresponding to collocated solution/quadrature nodes in DG, were distributed along the Gauss--Lobatto points for tensor-product elements and $\alpha$-optimized points~\citep{Hesthaven2008DG} for simplex elements. We note here that the any set of solution nodes (within reason) can be used without affecting the bounds-preserving properties of the proposed scheme as the approach can enforce constraints across the entire solution polynomial. 

Temporal integration was performed using a three-stage, third-order strong stability preserving Runge--Kutta scheme (SSP-RK3~\citep{Gottlieb2001}). Limiting was applied after each substage of the SSP-RK3 scheme using three iterations of the Newton--Raphson/gradient descent method. Additionally, limiting was applied at the initial time $t = 0$ after interpolation to the solution nodes to ensure that the solution was initially continuously bounds-preserving. In some experiments, comparisons were performed between the proposed continuously bounds-preserving approach and a standard discrete limiting approach, which coincides to limiting based on the solution value at the discrete solution nodes. This discretely bounds-preserving approach could be recovered by simply taking zero iterations of the optimization problem as the initial guess is taken over the discrete solution nodes. 

To evaluate the bounds-preserving properties of the scheme in a continuous sense, the constraints were sampled at a large number of points as a post-processing operation to yield a sufficiently accurate approximation of the continuous minima of the constraint functionals. In this work, we utilize 100 equispaced nodes per spatial dimension. Depending on the experiment, this evaluation was either performed at a fixed time (i.e., a spatial minimum) or continuously in time (i.e., a spatio-temporal minimum). The latter was approximated by evaluating the spatial over-sampling at each discrete temporal integration step. For problems with multiple constraint functionals, the minimum of the respective minima of the constraint functionals was used. 

As the proposed approach only enforces given constraints on the solution, it may not guarantee that the solution will be well-behaved around strong discontinuities if the constraints are not restrictive enough. For robustness around discontinuities, it is generally necessary to combine bounds-preserving limiters with additional numerical stabilization methods (e.g., slope limiting~\citep{Zhang2011b}) or enforce more restrictive constraints (e.g., minimum entropy constraints~\citep{Dzanic2022,Ching2023}). In this work, we utilize a simple subcell-type slope limiter based on a low-order subcell approximation used in the works of \citet{Pazner2021}, \citet{RuedaRamrez2022}, and \citet{Lin2023}, the details of which are given in \cref{app:subcell}. We note here that the choice of numerical stabilization method for discontinuity capturing is orthogonal to the purpose of bounds-preserving limiters (see \citet{Zhang2011b}). As long as it does not affect the assumptions that the proposed scheme is built upon (e.g., conservation, bounds-preservation on the element-wise mean, etc.), the combination with slope limiting is largely irrelevant with regards to the benefits and novelties of the proposed scheme, and we do not make any claims about the efficacy of the chosen slope limiter in comparison to other approaches. One may instead use one of a variety of other stabilization approaches of varying degrees of robustness and accuracy.

\section{Results}\label{sec:results}
\subsection{Linear advection}
The proposed approach was first applied to scalar linear transport problems through the linear advection equation, given in conservative form as
\begin{equation}
    \partial_t u + \boldsymbol{\nabla}{\cdot}\left(\mathbf{c}(\mathbf{x}) u\right) = 0,
\end{equation}
where $\mathbf{c}(\mathbf{x})$ is some constant or spatially-varying advection field. As with many scalar systems, entropy solutions of the linear advection equation satisfy a maximum principle, i.e.,
\begin{equation}
    u(\mathbf{x}, t) \in \left[u_{0,\min}, u_{0,\max}\right] \ \ \forall \ \mathbf{x}, t,
\end{equation}
where
\begin{equation}
    u_{0,\min} = \underset{\mathbf{x}}{\min}\ u_0(\mathbf{x}) \quad \text{and} \quad u_{0,\max} =\underset{\mathbf{x}}{\max}\ u_0(\mathbf{x}).
\end{equation}
This was simply enforced in terms of a pair of constraint functionals,
\begin{equation}
    g_1(u) = u - u_{0,\min} \quad \text{and} \quad g_2(u) = u_{0,\max} - u.
\end{equation}

\subsubsection{Advecting waveforms}
The bounds-preserving and error convergence properties of the proposed approach were first evaluated through the advecting waveforms problem. The domain is taken as $\Omega = [0, 1]$ with periodic boundary conditions, and the initial conditions are given as
\begin{equation}
    u_0(x) = \begin{cases}
        \exp\left( -300(2x - 0.3)^2 \right), &\mbox{if } |2x - 0.3| \leq 0.25, \\
        1, &\mbox{if } |2x - 0.9| \leq 0.2, \\
        \sqrt{1 - \left(\frac{2x-1.6}{0.2}\right)^2}, &\mbox{if } |2x - 1.6| \leq 0.2, \\
        0, &\mbox{else}.
    \end{cases} 
\end{equation}
These conditions yield waveforms with varying degrees of continuity, which allows for the evaluation of the approach for both smooth and discontinuous solutions. The advection velocity was simply set as $c = 1$, and global maximum principle bounds were enforced, i.e., $u \in [0,1]$. 

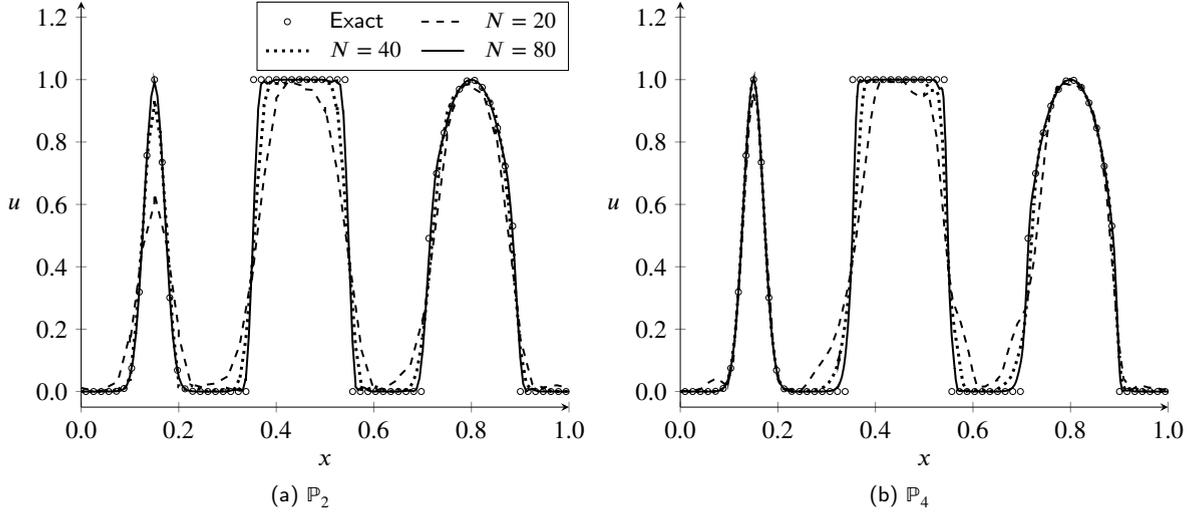
\begin{figure}[htbp!]
        \centering
        \subfloat[$\mathbb P_2$]{\adjustbox{width=0.48\linewidth,valign=b}{    \begin{tikzpicture}[spy using outlines={rectangle, height=3cm,width=2.3cm, magnification=3, connect spies}]
		\begin{axis}[name=plot1,
		    axis x line=left,
            axis y line=left,
		    xlabel={$x$},
		    xtick={0,0.2,0.4,0.6,0.8,1},
    		xmin=0,
    		xmax=1,
    		x tick label style={
        		/pgf/number format/.cd,
            	fixed,
            	fixed zerofill,
            	precision=1,
        	    /tikz/.cd},
    		ylabel={$u$},
    		ylabel style={rotate=-90},
		    ytick={0,0.2,0.4,0.6,0.8,1,1.2},
    		ymin=-0.05,
    		ymax=1.25,
    		y tick label style={
        		/pgf/number format/.cd,
            	fixed,
            	fixed zerofill,
            	precision=1,
        	    /tikz/.cd},
    		legend style={at={(1,1)},anchor=north east ,font=\small, column sep=0.2cm},
    		legend cell align={left},
    		legend columns=2,
    		style={font=\normalsize}]
    		
			\addplot[color=black, style={ultra thin}, only marks, mark=o, mark options={scale=0.6}, mark repeat = 8, mark phase = 6]
				table[x=x,y=u,col sep=comma,unbounded coords=jump]{./figs/data/shapes_exact.csv};
    		\addlegendentry{Exact}    		

			\addplot[color={black}, style={thick, dashed}]
				table[x=x,y=u,col sep=comma,unbounded coords=jump]{./figs/data/shapes_p2_n20_continuous_solution.csv};
    		\addlegendentry{$N = 20$}
      
			\addplot[color={black}, style={very thick, dotted}]
				table[x=x,y=u,col sep=comma,unbounded coords=jump]{./figs/data/shapes_p2_n40_continuous_solution.csv};
    		\addlegendentry{$N = 40$}
      
			\addplot[color={black}, style={thick}]
				table[x=x,y=u,col sep=comma,unbounded coords=jump]{./figs/data/shapes_p2_n80_continuous_solution.csv};
    		\addlegendentry{$N = 80$}

		\end{axis}

	\end{tikzpicture}}}
        \subfloat[$\mathbb P_4$]{\adjustbox{width=0.48\linewidth,valign=b}{    \begin{tikzpicture}[spy using outlines={rectangle, height=3cm,width=2.3cm, magnification=3, connect spies}]
		\begin{axis}[name=plot1,
		    axis x line=left,
            axis y line=left,
		    xlabel={$x$},
		    xtick={0,0.2,0.4,0.6,0.8,1},
    		xmin=0,
    		xmax=1,
    		x tick label style={
        		/pgf/number format/.cd,
            	fixed,
            	fixed zerofill,
            	precision=1,
        	    /tikz/.cd},
    		ylabel={$u$},
    		ylabel style={rotate=-90},
		    ytick={0,0.2,0.4,0.6,0.8,1,1.2},
    		ymin=-0.05,
    		ymax=1.25,
    		y tick label style={
        		/pgf/number format/.cd,
            	fixed,
            	fixed zerofill,
            	precision=1,
        	    /tikz/.cd},
    		legend style={at={(1,1)},anchor=north east ,font=\small, column sep=0.2cm},
    		legend cell align={left},
    		legend columns=2,
    		style={font=\normalsize}]
    		
			\addplot[color=black, style={ultra thin}, only marks, mark=o, mark options={scale=0.6}, mark repeat = 8, mark phase = 6]
				table[x=x,y=u,col sep=comma,unbounded coords=jump]{./figs/data/shapes_exact.csv};

			\addplot[color={black}, style={thick, dashed}]
				table[x=x,y=u,col sep=comma,unbounded coords=jump]{./figs/data/shapes_p4_n20_continuous_solution.csv};
      
			\addplot[color={black}, style={very thick, dotted}]
				table[x=x,y=u,col sep=comma,unbounded coords=jump]{./figs/data/shapes_p4_n40_continuous_solution.csv};
      
			\addplot[color={black}, style={thick}]
				table[x=x,y=u,col sep=comma,unbounded coords=jump]{./figs/data/shapes_p4_n80_continuous_solution.csv};

		\end{axis}

	\end{tikzpicture}}}
        \caption{\label{fig:advection_profiles} Solution profiles for the advecting waveforms problem at $t=1$ computed using a $\mathbb P_2$ (left) and $\mathbb P_4$ (right) approximation at varying mesh resolution with \emph{continuously} bounds-preserving limiting and global max principle bounds ($u \in [0,1]$).  }
    \end{figure}

The solutions profiles after one advection cycle ($t = 1$) through the domain as computed by a $\mathbb P_2$ and $\mathbb P_4$ scheme with varying mesh resolution are shown in 
\cref{fig:advection_profiles}. It can be seen that the approach shows clear qualitative convergence to the analytic solution with increasing resolution, both with the $\mathbb P_2$ and $\mathbb P_4$ approximation. Furthermore, the enforcement of the global maximum principle bounds ensured that the solution remained well-behaved in the vicinity of discontinuities. 

    \setlength{\extrarowheight}{.2em}

A more quantitative evaluation of the proposed approach was performed by analyzing the error convergence rates. The $L^1$ norm of the solution error after one advection cycle was computed via quadrature at the solution nodes at varying approximation orders using both the proposed continuous limiting approach and a discrete limiting approach, the results of which are shown in \cref{tab:advection_cbp_error} and \cref{tab:advection_dbp_error}, respectively. Both approaches showed very similar error behavior, with nearly identical convergence rates at all approximation orders. Furthermore, the magnitude of the error was also similar, with the continuous limiting approach showing marginally less error than the discrete limiting approach at all approximation orders and mesh resolution.

    \begin{figure}[htbp!] 
        \centering
        \begin{tabular}{|r | cccc |}
        \hline
        $N$ & $\mathbb P_2$ & $\mathbb P_3$ & $\mathbb P_4$ & $\mathbb P_5$ \\ 
        \hline
        20 & $7.52 \times 10^{-2}$ & $7.63 \times 10^{-2}$ & $7.77 \times 10^{-2}$ & $7.69 \times 10^{-2}$ \\
        40 & $3.50 \times 10^{-2}$ & $3.24 \times 10^{-2}$ & $3.39 \times 10^{-2}$ & $3.59 \times 10^{-2}$ \\
        60 & $2.16 \times 10^{-2}$ & $2.03 \times 10^{-2}$ & $2.16 \times 10^{-2}$ & $2.30 \times 10^{-2}$ \\
        80 & $1.55 \times 10^{-2}$ & $1.48 \times 10^{-2}$ & $1.57 \times 10^{-2}$ & $1.68 \times 10^{-2}$ \\
        100 & $1.22 \times 10^{-2}$ & $1.16 \times 10^{-2}$ & $1.23 \times 10^{-2}$ & $1.32 \times 10^{-2}$ \\
        120 & $1.01 \times 10^{-2}$ & $9.52 \times 10^{-3}$ & $1.01 \times 10^{-2}$ & $1.09 \times 10^{-2}$ \\
        \hline
        \textbf{RoC} & $1.13$ & $1.16$ & $1.14$ & $1.09$\\
        \hline
        \end{tabular}
        \captionof{table}{\label{tab:advection_cbp_error} Convergence in the $L^1$ norm for the advecting waveforms problem at varying approximation order with \emph{continuously} bounds-preserving limiting and global max principle bounds ($u \in [0,1]$).}
    \end{figure}
    \begin{figure}[htbp!] 
        \centering
        \begin{tabular}{|r | cccc |}
        \hline
        $N$ & $\mathbb P_2$ & $\mathbb P_3$ & $\mathbb P_4$ & $\mathbb P_5$ \\ 
        \hline
        20 & $8.68 \times 10^{-2}$ & $8.35 \times 10^{-2}$ & $8.40 \times 10^{-2}$ & $8.34 \times 10^{-2}$ \\
        40 & $4.05 \times 10^{-2}$ & $3.66 \times 10^{-2}$ & $3.69 \times 10^{-2}$ & $3.91 \times 10^{-2}$ \\
        60 & $2.52 \times 10^{-2}$ & $2.31 \times 10^{-2}$ & $2.36 \times 10^{-2}$ & $2.50 \times 10^{-2}$ \\
        80 & $1.81 \times 10^{-2}$ & $1.69 \times 10^{-2}$ & $1.72 \times 10^{-2}$ & $1.83 \times 10^{-2}$ \\
        100 & $1.43 \times 10^{-2}$ & $1.33 \times 10^{-2}$ & $1.35 \times 10^{-2}$ & $1.44 \times 10^{-2}$ \\
        120 & $1.18 \times 10^{-2}$ & $1.10 \times 10^{-2}$ & $1.11 \times 10^{-2}$ & $1.19 \times 10^{-2}$ \\
        \hline
        \textbf{RoC} & $1.12$ & $1.13$ & $1.13$ & $1.09$\\
        \hline
        \end{tabular}
        \captionof{table}{\label{tab:advection_dbp_error} Convergence in the $L^1$ norm for the advecting waveforms problem at varying approximation order with \emph{discretely} bounds-preserving limiting and global max principle bounds ($u \in [0,1]$).}
    \end{figure}

The bounds-preserving properties of the proposed continuous limiting approach was then verified by evaluating the spatio-temporal minimum of the constraint functionals over the domain and simulation time, i.e., $\min \ (g_1(u), g_2(u))$ sampled at 100 points per element at each temporal integration step. These minima are tabulated in \cref{tab:advection_cbp_gmin} for the continuous limiting approach. It can be seen that the proposed approach effectively enforces constraints \emph{continuously} down to machine precision levels, even with coarse mesh resolution. In comparison, the results of the discrete limiting approach, tabulated in \cref{tab:advection_dbp_gmin}, showed significant violations of the bounds which did not reduce with increasing resolution. In fact, this bounds violation stayed constant with increasing resolution, which can be attributed to the initial interpolation error of a discontinuous field onto a piece-wise continuous solution space. As such, even with highly-resolved meshes, discrete limiting approaches would likely be inappropriate for situations where constraints must be enforced continuously whereas the proposed approach would be effective regardless of the resolution. 

    \begin{figure}[htbp!] 
        \centering
        \begin{tabular}{|r | cccc |}
        \hline
        $N$ & $\mathbb P_2$ & $\mathbb P_3$ & $\mathbb P_4$ & $\mathbb P_5$ \\ 
        \hline
        20 & ${0.0}$ & $-1.33 \times 10^{-15}$ & ${0.0}$ & $-4.44 \times 10^{-16}$ \\
        40 & ${0.0}$ & $-1.33 \times 10^{-15}$ & ${0.0}$ & $-4.44 \times 10^{-16}$ \\
        60 & $-5.55 \times 10^{-17}$ & $-1.33 \times 10^{-15}$ & ${0.0}$ & $-4.44 \times 10^{-16}$ \\
        80 & ${0.0}$ & $-1.33 \times 10^{-15}$ & ${0.0}$ & $-4.44 \times 10^{-16}$ \\
        100 & ${0.0}$ & $-1.33 \times 10^{-15}$ & ${0.0}$ & $-4.44 \times 10^{-16}$ \\
        120 & ${0.0}$ & $-1.33 \times 10^{-15}$ & ${0.0}$ & $-4.44 \times 10^{-16}$ \\
        \hline
        \end{tabular}
        \captionof{table}{\label{tab:advection_cbp_gmin} Spatio-temporal minimum constraint functional value encountered during simulation for the advecting waveforms problem at varying approximation order and mesh resolution with \emph{continuously} bounds-preserving limiting and global max principle bounds ($u \in [0,1]$).}
    \end{figure}
    \begin{figure}[htbp!] 
        \centering
        \begin{tabular}{|r | cccc |}
        \hline
        $N$ & $\mathbb P_2$ & $\mathbb P_3$ & $\mathbb P_4$ & $\mathbb P_5$ \\ 
        \hline
        20 & $-1.25 \times 10^{-1}$ & $-1.31 \times 10^{-1}$ & $-1.32 \times 10^{-1}$ & $-1.32 \times 10^{-1}$ \\
        40 & $-1.25 \times 10^{-1}$ & $-1.31 \times 10^{-1}$ & $-1.32 \times 10^{-1}$ & $-1.32 \times 10^{-1}$ \\
        60 & $-1.25 \times 10^{-1}$ & $-1.31 \times 10^{-1}$ & $-1.32 \times 10^{-1}$ & $-1.32 \times 10^{-1}$ \\
        80 & $-1.25 \times 10^{-1}$ & $-1.31 \times 10^{-1}$ & $-1.32 \times 10^{-1}$ & $-1.32 \times 10^{-1}$ \\
        100 & $-1.25 \times 10^{-1}$ & $-1.31 \times 10^{-1}$ & $-1.32 \times 10^{-1}$ & $-1.32 \times 10^{-1}$ \\
        120 & $-1.25 \times 10^{-1}$ & $-1.31 \times 10^{-1}$ & $-1.32 \times 10^{-1}$ & $-1.32 \times 10^{-1}$ \\
        \hline
        \end{tabular}
        \captionof{table}{\label{tab:advection_dbp_gmin} Spatio-temporal minimum constraint functional value encountered during simulation for the advecting waveforms problem at varying approximation order and mesh resolution with \emph{discretely} bounds-preserving limiting and global max principle bounds ($u \in [0,1]$).}
    \end{figure}

\subsubsection{Solid body rotation}
The extension to two-dimensional problems and unstructured meshes was then performed through the solid body rotation problem of \citet{LeVeque1996}. In this problem, the domain is set as $\Omega = [0, 1]^2$ with periodic boundary conditions, and the initial conditions are given
\begin{equation}
    u_0(x) = \begin{cases}
        1, &\text{if } (x - 0.5)^2 + (y - 0.75)^2 \leq 0.15^2 \\
        & \quad \text{ and } x,y \notin [0.475,0.525]\times[0.6, 0.85], \\
        0.25\left (1 + \cos \left (\frac{\pi}{0.15} \sqrt{(x - 0.25)^2 + (y - 0.5)^2}\right ) \right ), &\text{if } (x - 0.25)^2 + (y - 0.5)^2 \leq 0.15^2, \\
        1 - \sqrt{(x - 0.5)^2 + (y - 0.25)^2}/0.15, &\text{if } (x - 0.5)^2 + (y - 0.25)^2 \leq 0.15^2, \\
        0, &\text{else}.
    \end{cases} 
\end{equation}
The advection field is set as  
\begin{equation}
    c(\mathbf{x}) = [-2 \pi (y - 0.5),\ 2 \pi (x - 0.5)]^T,
\end{equation}
which results in counter-clockwise rotation of constant angular velocity about the center of the domain with a unit rotation period. Similarly to the advecting waveforms problem, the initial conditions yield solution profiles of varying degrees of continuity, including a $C^\infty$ cosinusoidal hump, a $C^0$ sharp cone, and a discontinuous notched cylinder. These profiles can present a challenge to high-resolution numerical schemes, particularly the notched cylinder which is sensitive to excessive diffusion in numerical stabilization methods.

\begin{figure}[htbp!]
        \centering
        \subfloat[$N = 32^2$]{
        \adjustbox{width=0.33\linewidth,valign=b}{\includegraphics{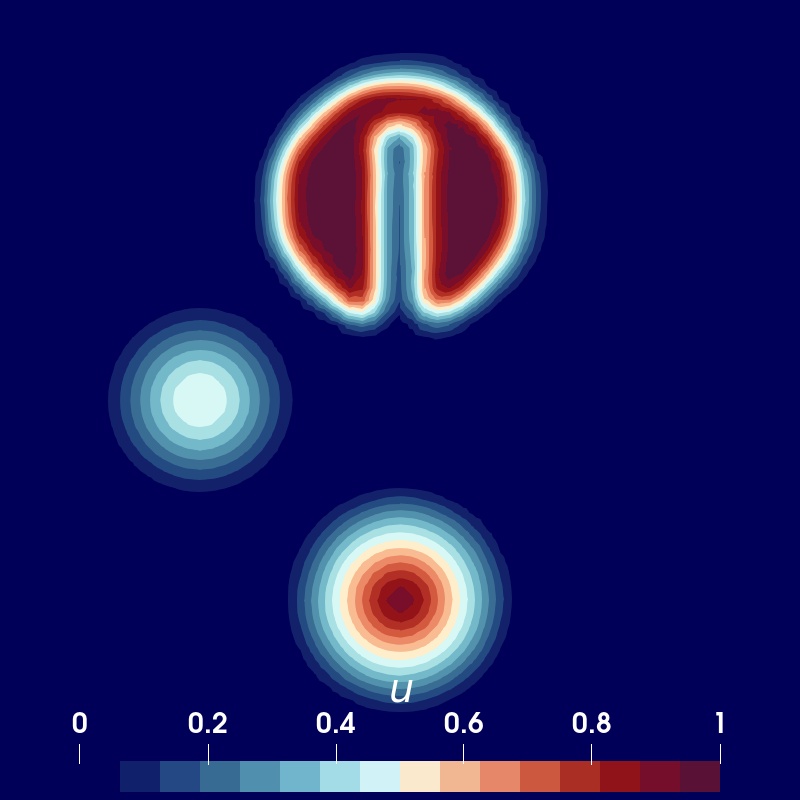}}}
        \subfloat[$N = 64^2$]{
        \adjustbox{width=0.33\linewidth,valign=b}{\includegraphics{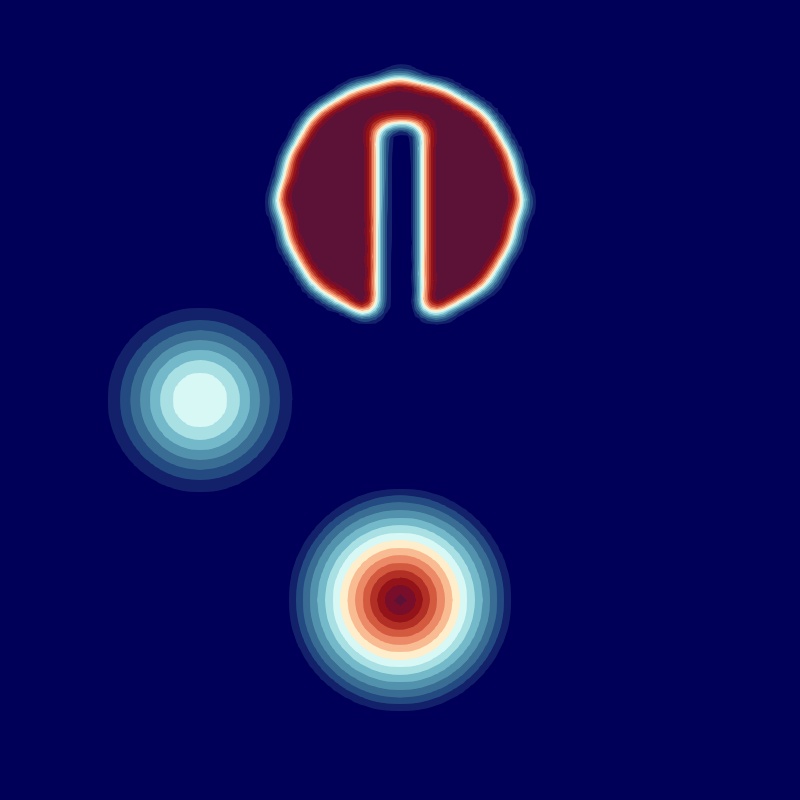}}}
        \subfloat[$N = 128^2$]{
        \adjustbox{width=0.33\linewidth,valign=b}{\includegraphics{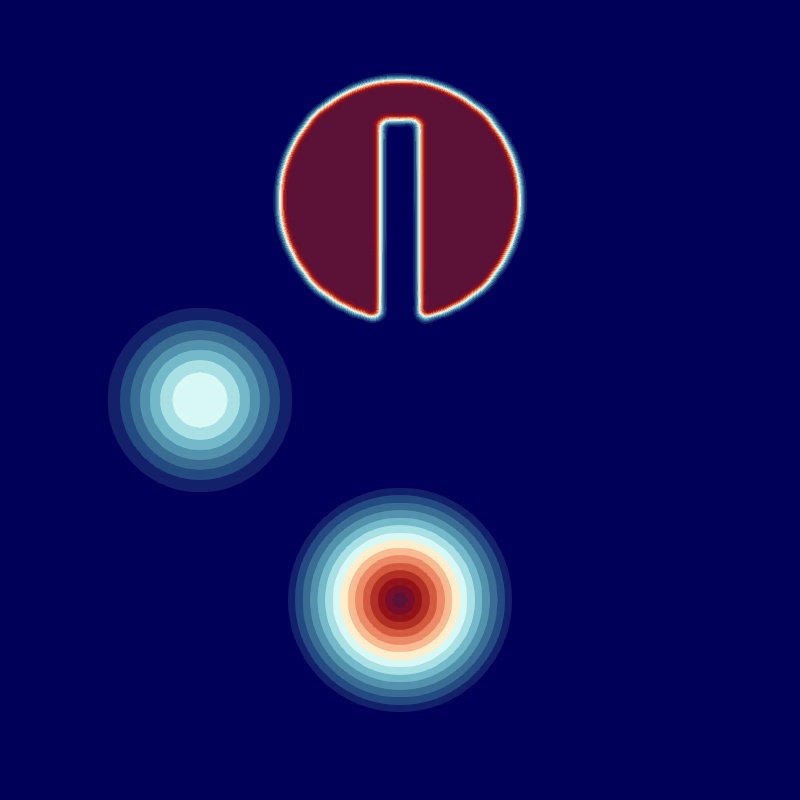}}}
        \caption{\label{fig:sbr_p2_contours} Solution contours for the solid body rotation problem at $t = 1$ computed using a $\mathbb P_2$ approximation with \emph{continuously} bounds-preserving limiting and global max principle bounds ($u \in [0,1]$) on a structured quadrilateral mesh with varying mesh resolution.} 
    \end{figure}

The problem was solved using a $\mathbb P_2$ approximation with global max principle bounds ($u \in [0,1]$) using both structured and unstructured meshes. The solution contours after one rotation period as computed on structured quadrilateral meshes of varying resolution are shown in \cref{fig:sbr_p2_contours}. It can be seen that for linear problems, the imposition of global max principle bounds effectively stabilized the solution around discontinuities. With increasing resolution, the numerical diffusion around the notched cylinder decreased, such that the initial profile was well-recovered at the highest resolution. The minimum value of the constraint functionals at the final time, sampled at $100^2$ points per element, was $g_{\min} = -2.5{\cdot}10^{-18}$, $-1.6{\cdot}10^{-18}$, and $-1.4{\cdot}10^{-19}$ for the $N = 32^2$, $64^2$, and $128^2$ mesh, respectively. The proposed limiting approach could ensure continuous boundedness to machine precision levels. 

\begin{figure}[htbp!]
        \centering
        \subfloat[$N = 2{\cdot}32^2$]{
        \adjustbox{width=0.33\linewidth,valign=b}{\includegraphics{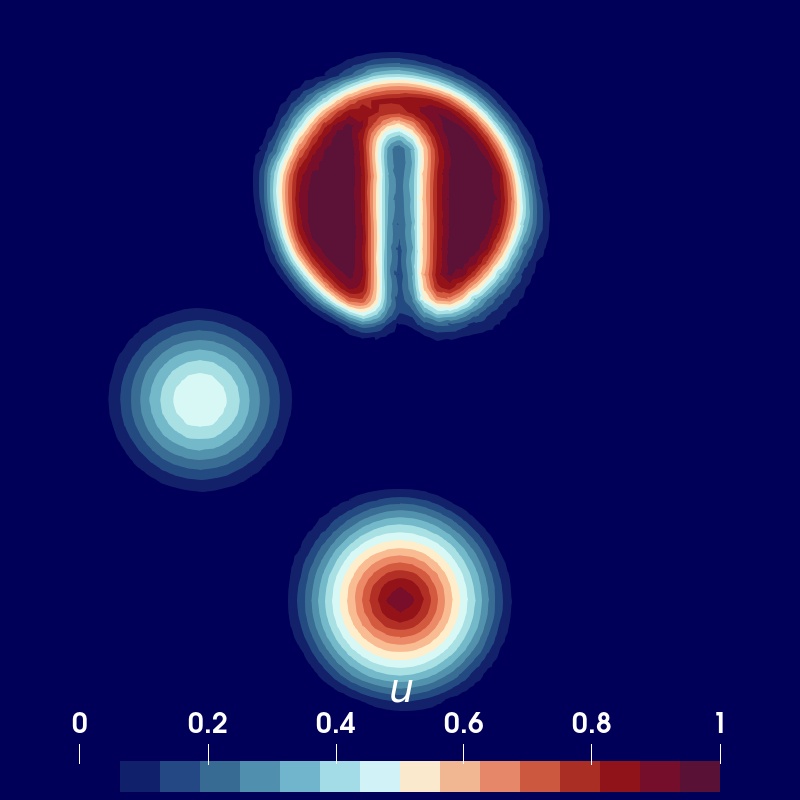}}}
        \subfloat[$N = 2{\cdot}64^2$]{
        \adjustbox{width=0.33\linewidth,valign=b}{\includegraphics{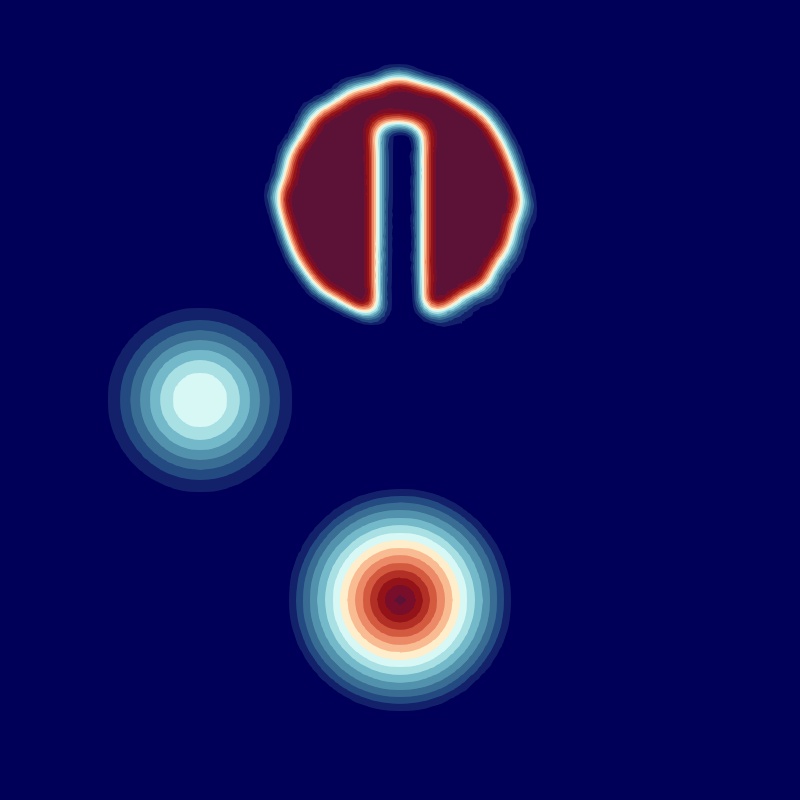}}}
        \subfloat[$N = 2{\cdot}128^2$]{
        \adjustbox{width=0.33\linewidth,valign=b}{\includegraphics{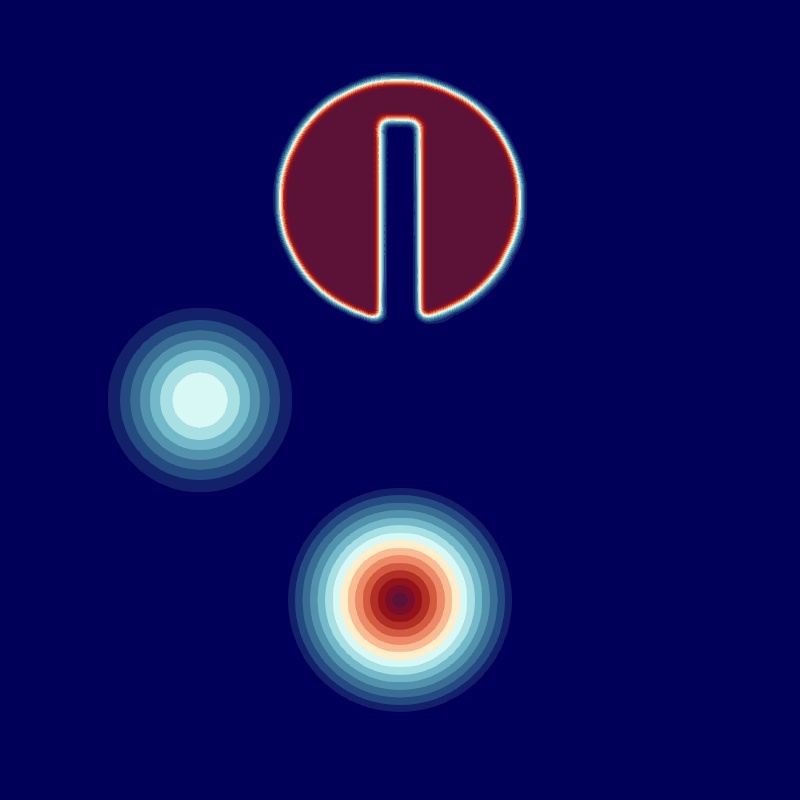}}}
        \caption{\label{fig:sbr_p2_contours_uns} Solution contours for the solid body rotation problem at $t = 1$ computed using a $\mathbb P_2$ approximation with \emph{continuously} bounds-preserving limiting and global max principle bounds ($u \in [0,1]$) on an unstructured triangular mesh with varying mesh resolution.} 
    \end{figure}

To verify the proposed approach for unstructured meshes, the solid body rotation problem was similarly computed on triangular meshes of equivalent resolution. The meshes were generated by uniformly triangulating the previous quadrilateral meshes, such that the number of triangular elements was double the number of quadrilateral elements. The solution contours after one rotation period as computed on the unstructured triangular meshes of varying resolution are shown in \cref{fig:sbr_p2_contours_uns}. It can be seen that nearly identical results are recovered with the unstructured meshes, indicating the proposed approach is nearly independent of the mesh discretization. Similarly as with the quadrilateral meshes, the minimum value of the constraint functional at the final time, sampled at $100^2$ points per element, was $g_{\min} = -4{\cdot}10^{-19}$, $-1.8{\cdot}10^{-19}$, and $-2.5{\cdot}10^{-18}$ for the $N = 2{\cdot}32^2$, $2{\cdot}64^2$, and $2{\cdot}128^2$ mesh, respectively.

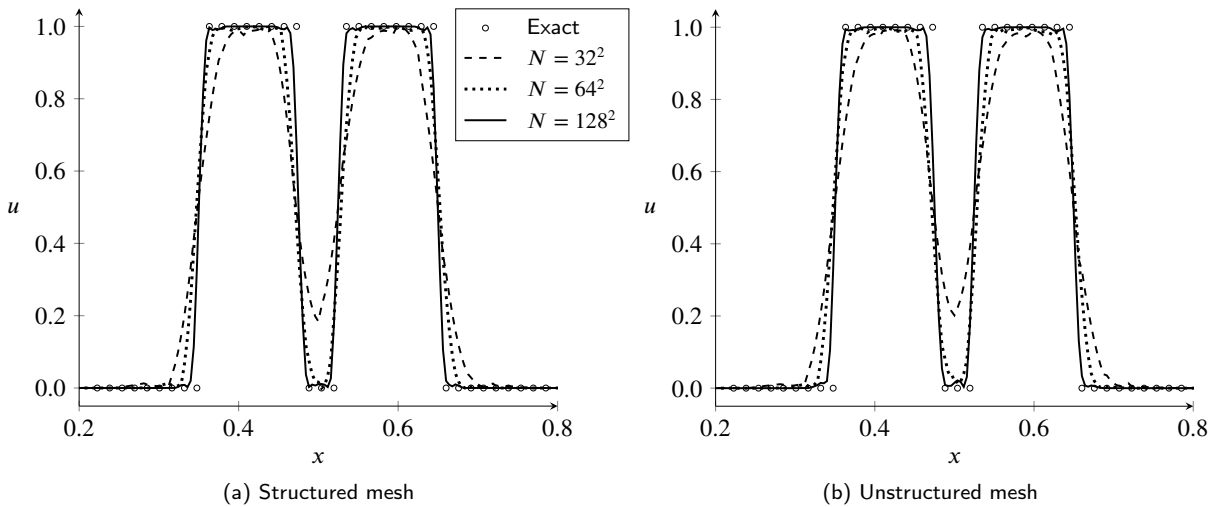
\begin{figure}[htbp!]
        \centering
        \subfloat[Structured mesh]{\adjustbox{width=0.51\linewidth,valign=b}{    \begin{tikzpicture}[spy using outlines={rectangle, height=3cm,width=2.3cm, magnification=3, connect spies}]
		\begin{axis}[name=plot1,
		    axis x line=left,
            axis y line=left,
		    xlabel={$x$},
		    xtick={0.2,0.4,0.6,0.8},
    		xmin=0.2,
    		xmax=0.8,
    		x tick label style={
        		/pgf/number format/.cd,
            	fixed,
            	fixed zerofill,
            	precision=1,
        	    /tikz/.cd},
    		ylabel={$u$},
    		ylabel style={rotate=-90},
		    ytick={0,0.2,0.4,0.6,0.8,1},
    		ymin=-0.05,
    		ymax=1.05,
    		y tick label style={
        		/pgf/number format/.cd,
            	fixed,
            	fixed zerofill,
            	precision=1,
        	    /tikz/.cd},
    		legend style={at={(1.15,1)},anchor=north east ,font=\small, column sep=0.2cm},
    		legend cell align={left},
    		style={font=\normalsize}]
    		
			\addplot[color=black, style={ultra thin}, only marks, mark=o, mark options={scale=0.6}, mark repeat = 4, mark phase = 6]
				table[x=x,y=u,col sep=comma,unbounded coords=jump]{./figs/data/sbr_exact.csv};
    		\addlegendentry{Exact}    		

			\addplot[color={black}, style={thick, dashed}]
				table[x=x,y=u,col sep=comma,unbounded coords=jump]{./figs/data/sbr_p2_32.csv};
    		\addlegendentry{$N = 32^2$}
      
			\addplot[color={black}, style={very thick, dotted}]
				table[x=x,y=u,col sep=comma,unbounded coords=jump]{./figs/data/sbr_p2_64.csv};
    		\addlegendentry{$N = 64^2$}
      
			\addplot[color={black}, style={thick}]
				table[x=x,y=u,col sep=comma,unbounded coords=jump]{./figs/data/sbr_p2_128.csv};
    		\addlegendentry{$N = 128^2$}

		\end{axis}

	\end{tikzpicture}}}
        \subfloat[Unstructured mesh]{\adjustbox{width=0.47\linewidth,valign=b}{    \begin{tikzpicture}[spy using outlines={rectangle, height=3cm,width=2.3cm, magnification=3, connect spies}]
		\begin{axis}[name=plot1,
		    axis x line=left,
            axis y line=left,
		    xlabel={$x$},
		    xtick={0.2,0.4,0.6,0.8},
    		xmin=0.2,
    		xmax=0.8,
    		x tick label style={
        		/pgf/number format/.cd,
            	fixed,
            	fixed zerofill,
            	precision=1,
        	    /tikz/.cd},
    		ylabel={$u$},
    		ylabel style={rotate=-90},
		    ytick={0,0.2,0.4,0.6,0.8,1},
    		ymin=-0.05,
    		ymax=1.05,
    		y tick label style={
        		/pgf/number format/.cd,
            	fixed,
            	fixed zerofill,
            	precision=1,
        	    /tikz/.cd},
    		legend style={at={(1.15,1)},anchor=north east ,font=\small, column sep=0.2cm},
    		legend cell align={left},
    		style={font=\normalsize}]
    		
			\addplot[color=black, style={ultra thin}, only marks, mark=o, mark options={scale=0.6}, mark repeat = 4, mark phase = 6]
				table[x=x,y=u,col sep=comma,unbounded coords=jump]{./figs/data/sbr_exact.csv};

			\addplot[color={black}, style={thick, dashed}]
				table[x=x,y=u,col sep=comma,unbounded coords=jump]{./figs/data/sbr_p2_32uns.csv};
      
			\addplot[color={black}, style={very thick, dotted}]
				table[x=x,y=u,col sep=comma,unbounded coords=jump]{./figs/data/sbr_p2_64uns.csv};
      
			\addplot[color={black}, style={thick}]
				table[x=x,y=u,col sep=comma,unbounded coords=jump]{./figs/data/sbr_p2_128uns.csv};

		\end{axis}

	\end{tikzpicture}}}
        \caption{\label{fig:sbr_profiles} Solution profile on the cross-section $y = 0.75$ for the solid body rotation problem at $t = 1$ computed using a $\mathbb P_2$ approximation with \emph{continuously} bounds-preserving limiting and global max principle bounds ($u \in [0,1]$). Solution showed for structured quadrilateral meshes (left) and unstructured triangular meshes (right) with varying mesh resolution.}
    \end{figure}

An additional comparison between the structured and unstructured meshes was compared by observing the cross-sectional profiles of the notched cylinder after one rotation period. The solution profile on the cross-section $y = 0.75$ as computed by both mesh types is shown in \cref{fig:sbr_profiles}. Both mesh types resulted in effectively identical solution profiles, not just in terms of stability around discontinuities but also dissipation of the solution features. These results further showcase the generalizability of the proposed approach to unstructured meshes. 
    
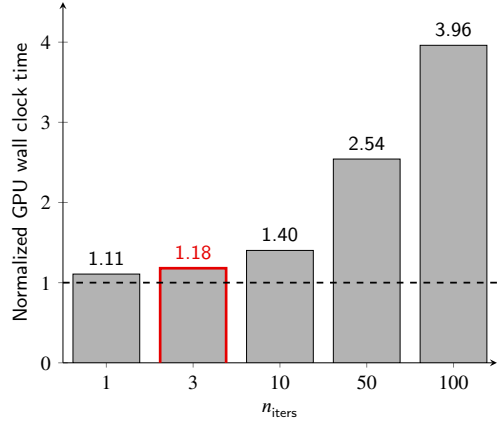
\begin{figure}[htbp!]
        \centering
        \adjustbox{width=0.4\linewidth,valign=b}{\begin{tikzpicture}[spy using outlines={rectangle, height=3cm,width=2.3cm, magnification=3, connect spies}]
\begin{axis} [ybar,
		axis line style={latex-latex},
	    axis x line=left,
        axis y line=left,
        bar width=30pt,
    	xmin=0.25, xmax=2.75,
    	ymin=0, ymax=4.5,
    	ylabel={Normalized GPU wall clock time},
        xtick={0.5, 1, 1.5, 2, 2.5},
    	xlabel={$n_{\text{iters}}$},
        xticklabels={$1$, $3$, $10$, $50$, $100$},
        clip mode=individual,
    	legend style={at={(0.03, 0.97)},anchor=north west},
    	legend cell align={left}]
     
\addplot[draw=black,fill=black!30]
    coordinates {
    	(0.5, 1.107) 
    	(1.0, 1.180) 
    	(1.5, 1.402) 
    	(2.0, 2.542) 
    	(2.5, 3.961) 
    };
\draw [red!90!black, very thick] (0.81, 0) -- (0.81, 1.18) -- (1.19, 1.18) -- (1.19, 0);
\draw [black, thick, dashed] (0.25, 1) -- (2.75, 1);

\node at (0.5, 1.307) {1.11};
\node at (1.0, 1.380) {\textcolor{red!90!black}{1.18}};
\node at (1.5, 1.602) {1.40};
\node at (2.0, 2.742) {2.54};
\node at (2.5, 4.161) {3.96};

\end{axis}

\end{tikzpicture}}
        \caption{\label{fig:cost} Computational cost of the continuously bounds-preserving limiting approach (normalized by the cost of the discrete limiting approach) with respect to the number of optimization iterations for the solid body rotation problem computed with a $\mathbb P_2$ approximation and $N = 128^2$ elements. Number of iterations used in this work $n_{\text{iters}} = 3$ highlighted in red. Unit value (dashed line) denotes no additional cost. }
    \end{figure}

An evaluation of the computational cost of the proposed approach was also performed for the two-dimensional experiments. The relative computational cost of the limiting was calculated by observing the elapsed wall clock time of the the structured $N=128^2$ experiment for a varying number of optimization iterations ($n_{\text{iters}}$). This cost was normalized by the computational cost of the discrete limiting approach, with a value of one corresponding to no additional computational cost. Comparisons were performed on one NVIDIA V100 GPU. The normalized cost is shown in \cref{fig:cost} for 1, 3, 10, 50, and 100 optimization iterations, with the standard number of iterations used in this work ($n_{\text{iters}} = 3$) highlighted in red. It can be seen that the cost of the proposed approach is actually quite minor, with 1 iteration requiring only an $11\%$ increase in the computational cost, largely due to an overhead cost associated with transforming the solution basis into its modal form, and 3 iterations requiring only an $18\%$ increase. It was not until the iteration was increased on the order of $n_{\text{iters}} = 30$ that the cost of the optimization approach became more than the underlying solver itself, with $n_{\text{iters}} = 50$ causing a $154\%$ cost increase and $n_{\text{iters}} = 100$ causing a $296\%$ increase. However, this is many more iterations than what is used in this work, which does not appreciably impact the overall computational cost. In fact, this indicates that the further computational effort can be expended without notably degrading the performance of the underlying solver, such as for sampling even more initial points for optimization or performing multiple optimization passes using separate initial guesses.

\subsection{Burgers' equation}
The proposed approach was extended to nonlinear governing equations through Burgers' equation, given in the form of a hyperbolic conservation law as
\begin{equation}
    \partial_t u + \partial_x \left(\frac{1}{2} u^2 \right) = 0.
\end{equation}
Similarly to the advection equation, the entropy solution of Burgers' equation satisfies a maximum principle. Due to the nonlinear nature of the equation, further numerical stabilization (e.g., slope limiting, entropy constraints) may be required to achieve a well-behaved solution as a global maximum principle may not be restrictive enough. However, as the purpose of this experiment is simply to compare the bounds-preserving properties of the proposed approach to other methods, we only enforce global max principle constraints and observe the behavior of the solution near constraint boundaries. 

\subsubsection{Compression/expansion wave}
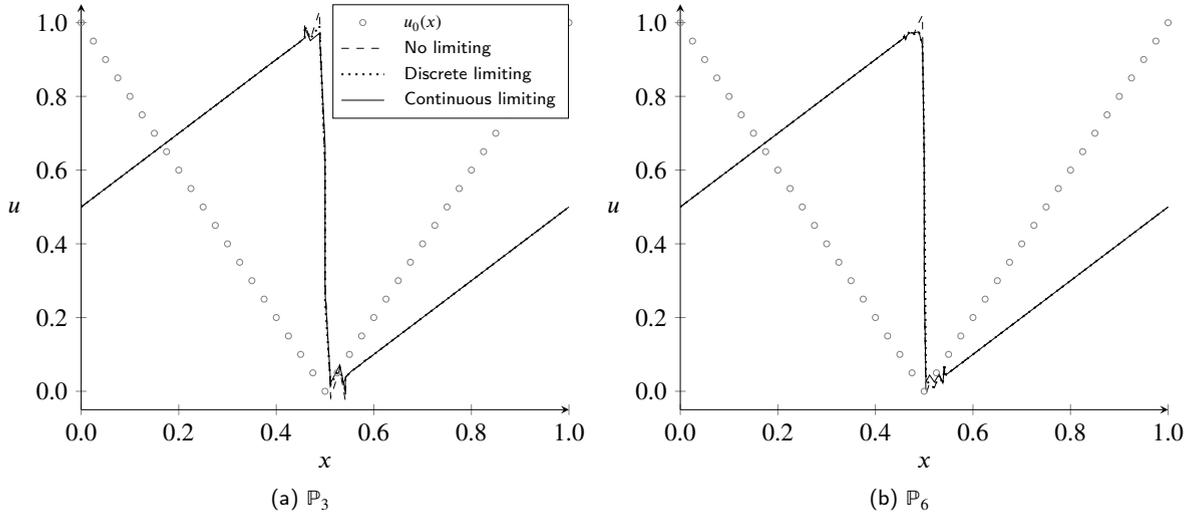
\begin{figure}[htbp!]
    \centering
    \subfloat[$\mathbb P_3$]{\adjustbox{width=0.48\linewidth,valign=b}{    \begin{tikzpicture}[spy using outlines={rectangle, height=6cm,width=2.3cm, magnification=3, connect spies}]
		\begin{axis}[name=plot1,
		    axis x line=left,
            axis y line=left,
		    xlabel={$x$},
		    xtick={0,0.2,0.4,0.6,0.8,1},
    		xmin=0,
    		xmax=1,
    		x tick label style={
        		/pgf/number format/.cd,
            	fixed,
            	fixed zerofill,
            	precision=1,
        	    /tikz/.cd},
    		ylabel={$u$},
    		ylabel style={rotate=-90},
		    ytick={0,0.2,0.4,0.6,0.8,1},
    		ymin=-0.05,
    		ymax=1.05,
    		y tick label style={
        		/pgf/number format/.cd,
            	fixed,
            	fixed zerofill,
            	precision=1,
        	    /tikz/.cd},
    		legend style={at={(1,1)},anchor=north east ,font=\scriptsize, column sep=0.2cm},
    		legend cell align={left},
    		style={font=\normalsize}]
    		
   
            \addplot [color=black!50, style={ultra thin}, only marks, mark=o, mark options={scale=0.6}] coordinates{
            (0.0, 1.)
            (0.025, 0.95)
            (0.05, 0.9)
            (0.075, 0.85)
            (0.1, 0.8)
            (0.125, 0.75)
            (0.15, 0.7)
            (0.175, 0.65)
            (0.2, 0.6)
            (0.225, 0.55)
            (0.25, 0.5)
            (0.275, 0.45)
            (0.3, 0.4)
            (0.325, 0.35)
            (0.35, 0.3)
            (0.375, 0.25)
            (0.4, 0.2)
            (0.425, 0.15)
            (0.45, 0.1)
            (0.475, 0.05)
            (0.5, 0.0) 
            (0.525, 0.05)
            (0.55, 0.1)
            (0.575, 0.15)
            (0.6, 0.2) 
            (0.625, 0.25)
            (0.65, 0.3)
            (0.675, 0.35)
            (0.7, 0.4) 
            (0.725, 0.45)
            (0.75, 0.5)
            (0.775, 0.55)
            (0.8, 0.6) 
            (0.825, 0.65)
            (0.85, 0.7)
            (0.875, 0.75)
            (0.9, 0.8) 
            (0.925, 0.85)
            (0.95, 0.9) 
            (0.975, 0.95)
            (1,1)};
    		\addlegendentry{$u_0(x)$}

			\addplot[color={black}, style={ dashed}]
				table[x=x,y=u,col sep=comma,unbounded coords=jump]{./figs/data/burgers_triangle_p3_n24_none_solution.csv};
    		\addlegendentry{No limiting}
      
			\addplot[color={black}, style={thick, dotted}]
				table[x=x,y=u,col sep=comma,unbounded coords=jump]{./figs/data/burgers_triangle_p3_n24_discrete_solution.csv};
    		\addlegendentry{Discrete limiting}
      
			\addplot[color={black}, style={}]
				table[x=x,y=u,col sep=comma,unbounded coords=jump]{./figs/data/burgers_triangle_p3_n24_continuous_solution.csv};
    		\addlegendentry{Continuous limiting}

		\end{axis} 		    
		
	\end{tikzpicture}}}
    \subfloat[$\mathbb P_6$]{\adjustbox{width=0.48\linewidth,valign=b}{    \begin{tikzpicture}[spy using outlines={rectangle, height=6cm,width=2.3cm, magnification=3, connect spies}]
		\begin{axis}[name=plot1,
		    axis x line=left,
            axis y line=left,
		    xlabel={$x$},
		    xtick={0,0.2,0.4,0.6,0.8,1},
    		xmin=0,
    		xmax=1,
    		x tick label style={
        		/pgf/number format/.cd,
            	fixed,
            	fixed zerofill,
            	precision=1,
        	    /tikz/.cd},
    		ylabel={$u$},
    		ylabel style={rotate=-90},
		    ytick={0,0.2,0.4,0.6,0.8,1},
    		ymin=-0.05,
    		ymax=1.05,
    		y tick label style={
        		/pgf/number format/.cd,
            	fixed,
            	fixed zerofill,
            	precision=1,
        	    /tikz/.cd},
    		legend style={at={(1,1)},anchor=north east ,font=\small, column sep=0.2cm},
    		legend cell align={left},
    		style={font=\normalsize}]
    		
   
            \addplot [color=black!50, style={ultra thin}, only marks, mark=o, mark options={scale=0.6}] coordinates{
            (0.0, 1.)
            (0.025, 0.95)
            (0.05, 0.9)
            (0.075, 0.85)
            (0.1, 0.8)
            (0.125, 0.75)
            (0.15, 0.7)
            (0.175, 0.65)
            (0.2, 0.6)
            (0.225, 0.55)
            (0.25, 0.5)
            (0.275, 0.45)
            (0.3, 0.4)
            (0.325, 0.35)
            (0.35, 0.3)
            (0.375, 0.25)
            (0.4, 0.2)
            (0.425, 0.15)
            (0.45, 0.1)
            (0.475, 0.05)
            (0.5, 0.0) 
            (0.525, 0.05)
            (0.55, 0.1)
            (0.575, 0.15)
            (0.6, 0.2) 
            (0.625, 0.25)
            (0.65, 0.3)
            (0.675, 0.35)
            (0.7, 0.4) 
            (0.725, 0.45)
            (0.75, 0.5)
            (0.775, 0.55)
            (0.8, 0.6) 
            (0.825, 0.65)
            (0.85, 0.7)
            (0.875, 0.75)
            (0.9, 0.8) 
            (0.925, 0.85)
            (0.95, 0.9) 
            (0.975, 0.95)
            (1,1)};

			\addplot[color={black}, style={dashed}]
				table[x=x,y=u,col sep=comma,unbounded coords=jump]{./figs/data/burgers_triangle_p6_n24_none_solution.csv};
      
			\addplot[color={black}, style={thick, dotted}]
				table[x=x,y=u,col sep=comma,unbounded coords=jump]{./figs/data/burgers_triangle_p6_n24_discrete_solution.csv};
      
			\addplot[color={black}, style={}]
				table[x=x,y=u,col sep=comma,unbounded coords=jump]{./figs/data/burgers_triangle_p6_n24_continuous_solution.csv};

		\end{axis} 		    
		
	\end{tikzpicture}}}
    \caption{\label{fig:burgers_full}Solution profiles for the Burgers' equation compression/expansion wave problem at $t=0.5$ computed using a $\mathbb P_3$ (left) and $\mathbb P_6$ (right) approximation and $N = 24$ elements with various limiting approaches and global max principle bounds ($u \in [0,1]$). }
\end{figure}

We consider the problem of a compression and expansion wave around a static point. The problem is solved on the periodic domain $\Omega = [0,1]$ with the initial conditions
\begin{equation}
    u_0(x) = \max \left( 1 - 2x, 2x - 1 \right).
\end{equation}
These initial conditions consist of a compression front decelerating to zero velocity followed by an expansion front, which results in a shock forming around the static point. The problem was solved with a $\mathbb P_3$ and $\mathbb P_6$ approximation with $N = 24$ elements using global max principle bounds ($u \in [0, 1]$). A comparison of the solution as computed using no limiting, discrete limiting, and the proposed continuous limiting method is shown in \cref{fig:burgers_full} at shock formation time ($t = 0.5$). As the constraints imposed are quite relaxed, the three approaches showed similar qualitative behavior, with the continuous limiting approach showing a marginally better behaved solution around the discontinuity, particularly at higher order. However, the efficacy of the proposed approach became much more pronounced when observing the continuous solution near the shock at a closer scale.

\begin{figure}[htbp!]
    \centering
    \subfloat[No limiting]{\adjustbox{width=0.33\linewidth,valign=b}{    \begin{tikzpicture}[spy using outlines={rectangle, height=6cm,width=2.3cm, magnification=3, connect spies}]
		\begin{axis}[name=plot1,
		    axis x line=left,
            axis y line=left,
		    xlabel={$x$},
		    xtick={0.45, 0.5, 0.55},
    		xmin=0.45,
    		xmax=0.55,
    		x tick label style={
        		/pgf/number format/.cd,
            	fixed,
            	fixed zerofill,
            	precision=2,
        	    /tikz/.cd},
    		ylabel={$u$},
    		ylabel style={rotate=-90},
		    ytick={0,0.2,0.4,0.6,0.8,1},
    		ymin=-0.1,
    		ymax=1.1,
    		y tick label style={
        		/pgf/number format/.cd,
            	fixed,
            	fixed zerofill,
            	precision=1,
        	    /tikz/.cd},
    		legend style={at={(1,1)},anchor=north east ,font=\small, column sep=0.2cm},
    		legend cell align={left},
            height=1.5*\axisdefaultheight,
            width=\axisdefaultheight,
    		style={font=\normalsize}]
    		
            \foreach \x in {10,...,13}{
			\addplot[color={black}, style={thick}]
				table[x=x,y=u,col sep=comma,unbounded coords=jump]{./figs/data/burgers_triangle_p3_n24_none_solution_e\x.csv};
            }

			\addplot[color={black}, style={}, only marks, mark=o, mark options={scale=0.8}]
				table[x=x,y=u,col sep=comma,unbounded coords=jump]{./figs/data/burgers_triangle_p3_n24_none_solution.csv};
    		    		
            \addplot[mark=none, color={black!80}, style={dotted, thick}] coordinates {(0,0) (1,0)};
            \addplot[mark=none, color={black!80}, style={dotted, thick}] coordinates {(0,1) (1,1)};
		\end{axis} 		    
		
	\end{tikzpicture}}}
    \subfloat[Discrete limiting]{\adjustbox{width=0.33\linewidth,valign=b}{    \begin{tikzpicture}[spy using outlines={rectangle, height=6cm,width=2.3cm, magnification=3, connect spies}]
		\begin{axis}[name=plot1,
		    axis x line=left,
            axis y line=left,
		    xlabel={$x$},
		    xtick={0.45, 0.5, 0.55},
    		xmin=0.45,
    		xmax=0.55,
    		x tick label style={
        		/pgf/number format/.cd,
            	fixed,
            	fixed zerofill,
            	precision=2,
        	    /tikz/.cd},
    		ylabel={$u$},
    		ylabel style={rotate=-90},
		    ytick={0,0.2,0.4,0.6,0.8,1},
    		ymin=-0.1,
    		ymax=1.1,
    		y tick label style={
        		/pgf/number format/.cd,
            	fixed,
            	fixed zerofill,
            	precision=1,
        	    /tikz/.cd},
    		legend style={at={(1,1)},anchor=north east ,font=\small, column sep=0.2cm},
    		legend cell align={left},
            height=1.5*\axisdefaultheight,
            width=\axisdefaultheight,
    		style={font=\normalsize}]
    		
            \foreach \x in {10,...,13}{
			\addplot[color={black}, style={thick}]
				table[x=x,y=u,col sep=comma,unbounded coords=jump]{./figs/data/burgers_triangle_p3_n24_discrete_solution_e\x.csv};
            }

			\addplot[color={black}, style={}, only marks, mark=o, mark options={scale=0.8}]
				table[x=x,y=u,col sep=comma,unbounded coords=jump]{./figs/data/burgers_triangle_p3_n24_discrete_solution.csv};
    		    		
            \addplot[mark=none, color={black!80}, style={dotted, thick}] coordinates {(0,0) (1,0)};
            \addplot[mark=none, color={black!80}, style={dotted, thick}] coordinates {(0,1) (1,1)};
		\end{axis} 		    
		
	\end{tikzpicture}}}
    \subfloat[Continuous limiting]{\adjustbox{width=0.33\linewidth,valign=b}{    \begin{tikzpicture}[spy using outlines={rectangle, height=6cm,width=2.3cm, magnification=3, connect spies}]
		\begin{axis}[name=plot1,
		    axis x line=left,
            axis y line=left,
		    xlabel={$x$},
		    xtick={0.45, 0.5, 0.55},
    		xmin=0.45,
    		xmax=0.55,
    		x tick label style={
        		/pgf/number format/.cd,
            	fixed,
            	fixed zerofill,
            	precision=2,
        	    /tikz/.cd},
    		ylabel={$u$},
    		ylabel style={rotate=-90},
		    ytick={0,0.2,0.4,0.6,0.8,1},
    		ymin=-0.1,
    		ymax=1.1,
    		y tick label style={
        		/pgf/number format/.cd,
            	fixed,
            	fixed zerofill,
            	precision=1,
        	    /tikz/.cd},
    		legend style={at={(1,1)},anchor=north east ,font=\small, column sep=0.2cm},
    		legend cell align={left},
            height=1.5*\axisdefaultheight,
            width=\axisdefaultheight,
    		style={font=\normalsize}]
    		
            \foreach \x in {10,...,13}{
			\addplot[color={black}, style={thick}]
				table[x=x,y=u,col sep=comma,unbounded coords=jump]{./figs/data/burgers_triangle_p3_n24_continuous_solution_e\x.csv};
            }

			\addplot[color={black}, style={}, only marks, mark=o, mark options={scale=0.8}]
				table[x=x,y=u,col sep=comma,unbounded coords=jump]{./figs/data/burgers_triangle_p3_n24_continuous_solution.csv};
    		    		
            \addplot[mark=none, color={black!80}, style={dotted, thick}] coordinates {(0,0) (1,0)};
            \addplot[mark=none, color={black!80}, style={dotted, thick}] coordinates {(0,1) (1,1)};
		\end{axis} 		    
		
	\end{tikzpicture}}}
    \caption{\label{fig:burgers_p3_zoom} Enlarged view of the solution profiles near the shock front for the Burgers' equation compression/expansion wave problem at $t=0.5$ computed using a $\mathbb P_3$ approximation and $N = 24$ elements with various limiting approaches and global max principle bounds. Circle markers denote discrete solution nodes. Bounds (isocontour of $g_1(u) = g_2(u) = 0$) shown by dotted lines.}
\end{figure}
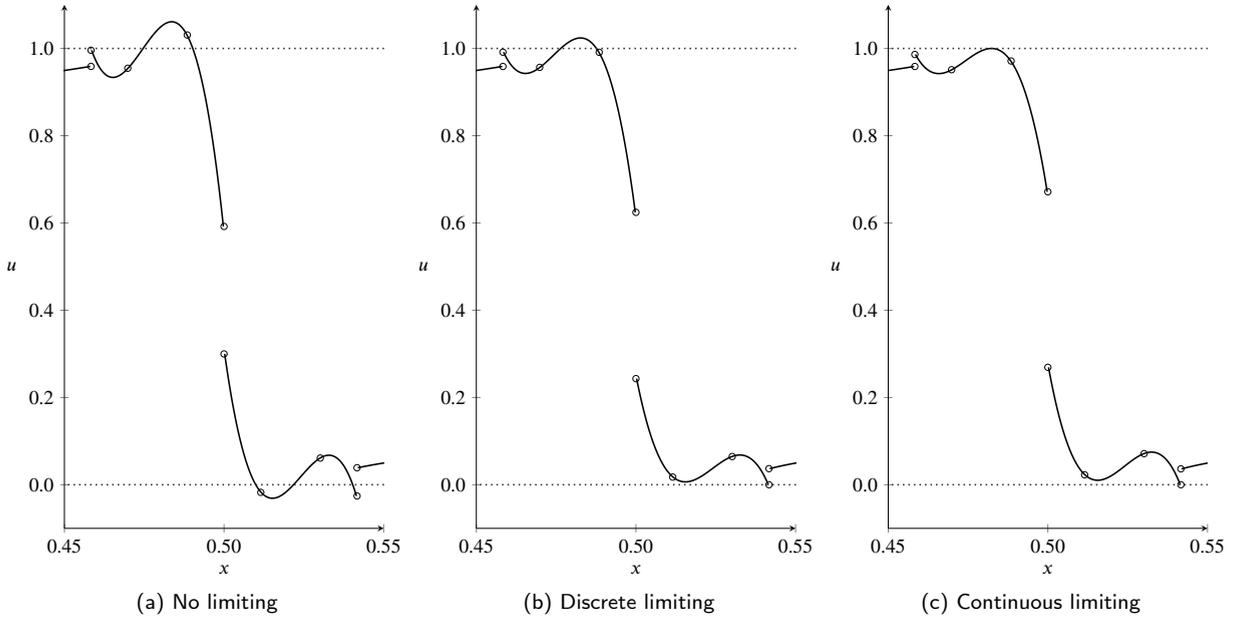

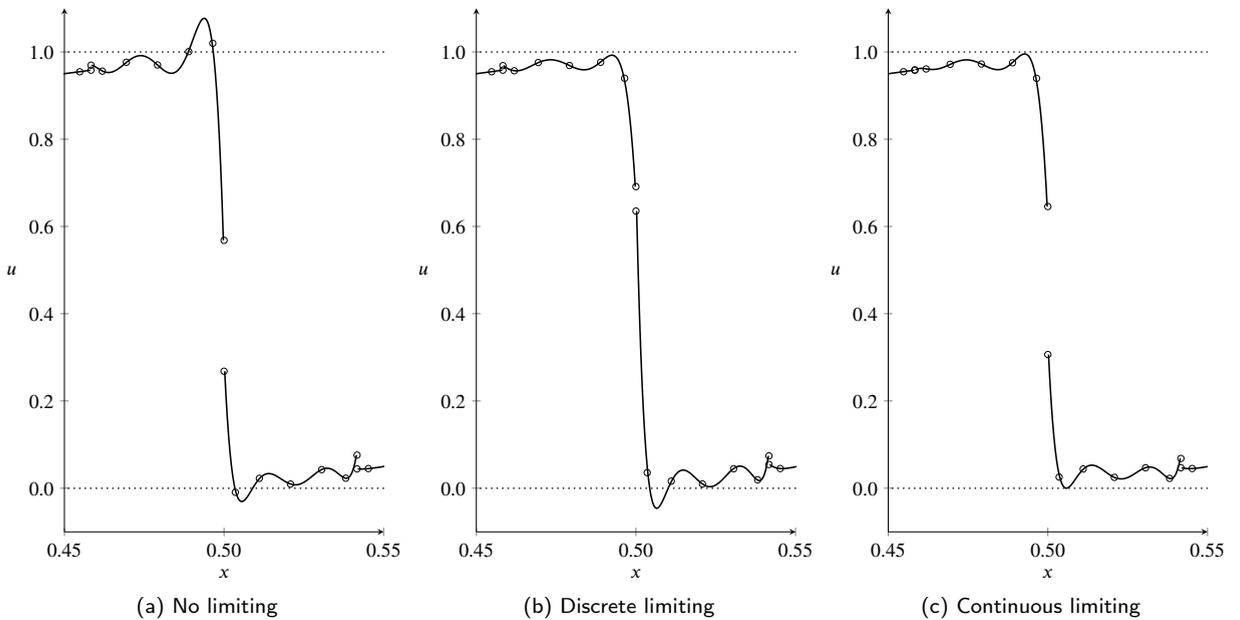
\begin{figure}[htbp!]
    \centering
    \subfloat[No limiting]{\adjustbox{width=0.33\linewidth,valign=b}{    \begin{tikzpicture}[spy using outlines={rectangle, height=6cm,width=2.3cm, magnification=3, connect spies}]
		\begin{axis}[name=plot1,
		    axis x line=left,
            axis y line=left,
		    xlabel={$x$},
		    xtick={0.45, 0.5, 0.55},
    		xmin=0.45,
    		xmax=0.55,
    		x tick label style={
        		/pgf/number format/.cd,
            	fixed,
            	fixed zerofill,
            	precision=2,
        	    /tikz/.cd},
    		ylabel={$u$},
    		ylabel style={rotate=-90},
		    ytick={0,0.2,0.4,0.6,0.8,1},
    		ymin=-0.1,
    		ymax=1.1,
    		y tick label style={
        		/pgf/number format/.cd,
            	fixed,
            	fixed zerofill,
            	precision=1,
        	    /tikz/.cd},
    		legend style={at={(1,1)},anchor=north east ,font=\small, column sep=0.2cm},
    		legend cell align={left},
            height=1.5*\axisdefaultheight,
            width=\axisdefaultheight,
    		style={font=\normalsize}]
    		
            \foreach \x in {10,...,13}{
			\addplot[color={black}, style={thick}]
				table[x=x,y=u,col sep=comma,unbounded coords=jump]{./figs/data/burgers_triangle_p6_n24_none_solution_e\x.csv};
            }

			\addplot[color={black}, style={}, only marks, mark=o, mark options={scale=0.8}]
				table[x=x,y=u,col sep=comma,unbounded coords=jump]{./figs/data/burgers_triangle_p6_n24_none_solution.csv};
    		    		
            \addplot[mark=none, color={black!80}, style={dotted, thick}] coordinates {(0,0) (1,0)};
            \addplot[mark=none, color={black!80}, style={dotted, thick}] coordinates {(0,1) (1,1)};
		\end{axis} 		    
		
	\end{tikzpicture}}}
    \subfloat[Discrete limiting]{\adjustbox{width=0.33\linewidth,valign=b}{    \begin{tikzpicture}[spy using outlines={rectangle, height=6cm,width=2.3cm, magnification=3, connect spies}]
		\begin{axis}[name=plot1,
		    axis x line=left,
            axis y line=left,
		    xlabel={$x$},
		    xtick={0.45, 0.5, 0.55},
    		xmin=0.45,
    		xmax=0.55,
    		x tick label style={
        		/pgf/number format/.cd,
            	fixed,
            	fixed zerofill,
            	precision=2,
        	    /tikz/.cd},
    		ylabel={$u$},
    		ylabel style={rotate=-90},
		    ytick={0,0.2,0.4,0.6,0.8,1},
    		ymin=-0.1,
    		ymax=1.1,
    		y tick label style={
        		/pgf/number format/.cd,
            	fixed,
            	fixed zerofill,
            	precision=1,
        	    /tikz/.cd},
    		legend style={at={(1,1)},anchor=north east ,font=\small, column sep=0.2cm},
    		legend cell align={left},
            height=1.5*\axisdefaultheight,
            width=\axisdefaultheight,
    		style={font=\normalsize}]
    		
            \foreach \x in {10,...,13}{
			\addplot[color={black}, style={thick}]
				table[x=x,y=u,col sep=comma,unbounded coords=jump]{./figs/data/burgers_triangle_p6_n24_discrete_solution_e\x.csv};
            }

			\addplot[color={black}, style={}, only marks, mark=o, mark options={scale=0.8}]
				table[x=x,y=u,col sep=comma,unbounded coords=jump]{./figs/data/burgers_triangle_p6_n24_discrete_solution.csv};
    		    		
            \addplot[mark=none, color={black!80}, style={dotted, thick}] coordinates {(0,0) (1,0)};
            \addplot[mark=none, color={black!80}, style={dotted, thick}] coordinates {(0,1) (1,1)};
		\end{axis} 		    
		
	\end{tikzpicture}}}
    \subfloat[Continuous limiting]{\adjustbox{width=0.33\linewidth,valign=b}{    \begin{tikzpicture}[spy using outlines={rectangle, height=6cm,width=2.3cm, magnification=3, connect spies}]
		\begin{axis}[name=plot1,
		    axis x line=left,
            axis y line=left,
		    xlabel={$x$},
		    xtick={0.45, 0.5, 0.55},
    		xmin=0.45,
    		xmax=0.55,
    		x tick label style={
        		/pgf/number format/.cd,
            	fixed,
            	fixed zerofill,
            	precision=2,
        	    /tikz/.cd},
    		ylabel={$u$},
    		ylabel style={rotate=-90},
		    ytick={0,0.2,0.4,0.6,0.8,1},
    		ymin=-0.1,
    		ymax=1.1,
    		y tick label style={
        		/pgf/number format/.cd,
            	fixed,
            	fixed zerofill,
            	precision=1,
        	    /tikz/.cd},
    		legend style={at={(1,1)},anchor=north east ,font=\small, column sep=0.2cm},
    		legend cell align={left},
            height=1.5*\axisdefaultheight,
            width=\axisdefaultheight,
    		style={font=\normalsize}]
    		
            \foreach \x in {10,...,13}{
			\addplot[color={black}, style={thick}]
				table[x=x,y=u,col sep=comma,unbounded coords=jump]{./figs/data/burgers_triangle_p6_n24_continuous_solution_e\x.csv};
            }

			\addplot[color={black}, style={}, only marks, mark=o, mark options={scale=0.8}]
				table[x=x,y=u,col sep=comma,unbounded coords=jump]{./figs/data/burgers_triangle_p6_n24_continuous_solution.csv};
    		    		
            \addplot[mark=none, color={black!80}, style={dotted, thick}] coordinates {(0,0) (1,0)};
            \addplot[mark=none, color={black!80}, style={dotted, thick}] coordinates {(0,1) (1,1)};
		\end{axis} 		    
		
	\end{tikzpicture}}}
    \caption{\label{fig:burgers_p6_zoom}Enlarged view of the solution profiles near the shock front for the Burgers' equation compression/expansion wave problem at $t=0.5$ computed using a $\mathbb P_6$ approximation and $N = 24$ elements with various limiting approaches and global max principle bounds. Circle markers denote discrete solution nodes. Bounds (isocontour of $g_1(u) = g_2(u) = 0$) shown by dotted lines.}
\end{figure}

A comparison of the proposed approach to the discrete limiting method and the unlimited numerical scheme is shown in the enlarged views in \cref{fig:burgers_p3_zoom} and \cref{fig:burgers_p6_zoom} for the $\mathbb P_3$ and $\mathbb P_6$ approximation, respectively. To highlight the effect of enforcing continuous constraints, the polynomial solution is plotted with the discrete solution nodes (i.e., the nodes where discrete limiting would be performed) shown as individual markers. For the $\mathbb P_3$ approximation, the unlimited numerical scheme violated the global maximum principle bounds, shown by the dotted lines, both in terms of the discrete solution nodes and the continuous solution polynomial. When discrete limiting was applied, it can be seen that the limiting approach ensured that discrete solution nodes resided strictly within the bounds, but the solution polynomial itself showed overshoots which violated the bounds. As such, if any applications required evaluating the solution at different nodal locations, it is likely that the solution would not satisfy the given constraints. In contrast, with the proposed continuous limiting approach, the \emph{entire} solution polynomial resided within the given bounds, which highlights the ability of the approach to continuously enforce constraints.

A similar analysis was performed with the $\mathbb P_6$ approximation, shown in \cref{fig:burgers_p6_zoom}. Due to the higher-order polynomial approximation, the resulting solution polynomial in the vicinity of the discontinuity was significantly more complex, showing more local extrema. As with the $\mathbb P_3$ experiment, the $\mathbb P_6$ results showed that the unlimited scheme and the discrete limiting method both violated the bounds across the solution polynomial, this time with undershoots instead of overshoots in the discrete limiting method. However, the proposed approach still ensured that the entire solution polynomial resided within the bounds, even with the high polynomial order causing a highly-oscillatory optimization landscape. 
These results particularly highlight the differences between the proposed continuous limiting method and simple discrete limiting, showcasing the efficacy of the approach for applications where continuous boundedness may be required.

\subsection{Euler equations}
The robustness of the approach was then evaluated for hyperbolic systems with nonlinear constraints through the Euler equations, given in terms of \cref{eq:gen_hype} as 
\begin{equation}\label{eq:euler}
    \mathbf{u} = \begin{bmatrix}
            \rho \\ \mathbf{m} \\ E
        \end{bmatrix} \quad  \mathrm{and} \quad \mathbf{F}(\mathbf{u}) = \begin{bmatrix}
            \mathbf{m}^T\\
            \mathbf{m}\otimes\mathbf{v} + P\mathbf{I}\\
        (E+P)\mathbf{v}^T
    \end{bmatrix},
\end{equation}
where $\rho$ is the density, $\mathbf{m}$ is the momentum, and $E$ is the total energy. The symbol $\mathbf{I}$ denotes the identity matrix in $\mathbb R^{d \times d}$, $\mathbf{v} = \mathbf{m}/\rho$ denotes the velocity, and $P$ denotes the pressure, computed as 
\begin{equation}
    P = (\gamma - 1)\rho e = (\gamma - 1)\left (E - \frac{1}{2} \mathbf{m}{\cdot}\mathbf{m}/ \rho \right),
\end{equation}
where $\gamma$ is the specific heat ratio and $e$ is the specific internal energy. Unless otherwise stated, the specific heat ratio is set as $\gamma = 1.4$. For brevity, we denote use $\mathbf{Q} = [\rho, \mathbf{v}, P]^T$ to represent a vector of primitive variables. 

Entropy solutions of the Euler equations reside within an invariant set corresponding to solutions with positive density and pressure (i.e., $\mathcal G = \left\{ \mathbf{u}\ |\ \rho \geq 0, P \geq 0 \right\}$) which correspond to physical constraints in gas dynamics. As such, bounds can be imposed in terms of a linear density constraint and a nonlinear (quadratic) pressure constraint as
\begin{equation}
    g_1(\mathbf{u}) = \rho - \rho_{\min} \quad \text{and} \quad g_2(\mathbf{u}) = P - P_{\min},
\end{equation}
with the limiting applied sequentially. While physically the density and pressure may be zero, it is typically necessary to impose some minimum value to avoid numerical issues with the Riemann solver. Therefore, we impose a minimum density/pressure bound as
\begin{equation*}
    \rho_{\min} = P_{\min} = 10^{-11}.
\end{equation*}

\subsubsection{Near-vacuum convecting density pulse}
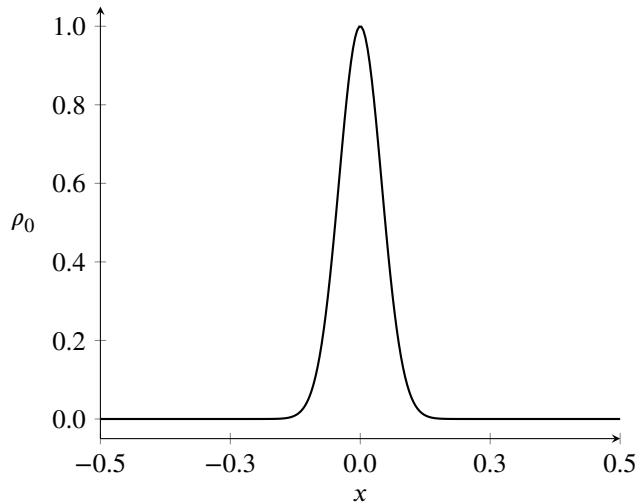
\begin{figure}[htbp!]
    \centering
        \begin{tikzpicture}[spy using outlines={rectangle, height=6cm,width=2.3cm, magnification=3, connect spies}]
		\begin{axis}[name=plot1,
		    axis x line=left,
            axis y line=left,
		    xlabel={$x$},
		    xtick={-0.5, -0.25,  0, 0.25, 0.5},
    		xmin=-0.5,
    		xmax=0.5,
    		x tick label style={
        		/pgf/number format/.cd,
            	fixed,
            	fixed zerofill,
            	precision=1,
        	    /tikz/.cd},
    		ylabel={$\rho_0$},
    		ylabel style={rotate=-90},
		    ytick={0,0.2,0.4,0.6,0.8,1},
    		ymin=-0.05,
    		ymax=1.05,
    		y tick label style={
        		/pgf/number format/.cd,
            	fixed,
            	fixed zerofill,
            	precision=1,
        	    /tikz/.cd},
    		legend style={at={(1,1)},anchor=north east ,font=\scriptsize, column sep=0.2cm},
    		legend cell align={left},
    		style={font=\normalsize}]

            \addplot [color=black, style={thick}, domain=-0.5:0.5, samples=250] {exp(-300*x*x)};

		\end{axis} 		    
		
	\end{tikzpicture}
    \caption{\label{fig:wave_ics} Initial density field for the near-vacuum convecting density pulse. }
\end{figure}

To verify that imposing positivity constraints via the proposed approach retains high-order accuracy at smooth extrema, we consider the problem of a convecting density pulse in a near-vacuum state. The domain is set as 
$\Omega = [-0.5, 0.5]$ with periodic boundary conditions, and the initial conditions are given
\begin{equation}
    \mathbf{Q}_0 (x) = 
    \begin{bmatrix}
        \exp \left(-200 x^2 \right) + 2 \rho_{\min}\\
        1\\
        2 P_{\min}
    \end{bmatrix}.
\end{equation}
The problem consists of a smooth Gaussian density wave, shown in \cref{fig:wave_ics}, convecting uniformly through the domain. The ambient density and pressure are set to $2{\cdot}10^{-11}$ which is very close to the minimum density/pressure bound. At these near-vacuum conditions, an unlimited scheme would quickly diverge as the density and pressure attain negative values, such that it is necessary to utilize positivity-preserving schemes. Furthermore, as the exact solution is only marginally within the admissible set, continuous limiting methods which only provide a \emph{sufficient} degree of limiting (e.g., limiting based on Bernstein polynomials) can degrade the high-order accuracy of the underlying numerical scheme as they may introduce numerical dissipation in regions where it is not necessary. 

In \citet{Zhang2011b}, it was shown that the \emph{discrete} positivity-preserving limiter retains the high-order accuracy of the underlying DG approach. To verify that the proposed continuous limiting approach also exhibits this behavior, the $L^{\infty}$ norm of the density error was evaluated after one pass-through of the domain ($t = 1$). The convergence of this error is shown in \cref{tab:euler_pulse_error} for varying approximation orders and mesh resolution. It can be seen that for the given positivity-preserving constraints, the proposed approach also retains the high-order $p+1$ convergence rates of the underlying scheme, even at smooth extrema, while effectively stabilizing the solution at near-vacuum conditions. This can largely be attributed to the property of the proposed approach in which no limiting is applied to the high-order solution locally if it is already continuously bounds-preserving (\cref{thm:nolimiting}).

    \begin{figure}[htbp!] 
        \centering
        \begin{tabular}{|r | cccc |}
        \hline
        $N$ & $\mathbb P_2$ & $\mathbb P_3$ & $\mathbb P_4$ & $\mathbb P_5$ \\ 
        \hline
        5 & $1.36 \times 10^{-1}$ & $7.01 \times 10^{-2}$ & $5.40 \times 10^{-2}$ & $7.41 \times 10^{-2}$ \\
        10 & $1.62 \times 10^{-2}$ & $5.47 \times 10^{-3}$ & $2.40 \times 10^{-3}$ & $1.90 \times 10^{-4}$ \\
        15 & $9.89 \times 10^{-3}$ & $1.24 \times 10^{-3}$ & $5.12 \times 10^{-5}$ & $1.59 \times 10^{-5}$ \\
        20 & $3.51 \times 10^{-3}$ & $1.85 \times 10^{-4}$ & $1.27 \times 10^{-5}$ & $2.03 \times 10^{-6}$ \\
        25 & $1.41 \times 10^{-3}$ & $5.62 \times 10^{-5}$ & $4.17 \times 10^{-6}$ & $9.41 \times 10^{-7}$ \\
        30 & $7.00 \times 10^{-4}$ & $2.50 \times 10^{-5}$ & $1.72 \times 10^{-6}$ & $5.34 \times 10^{-7}$ \\
        35 & $3.55 \times 10^{-4}$ & $1.30 \times 10^{-5}$ & $9.21 \times 10^{-7}$ & $3.30 \times 10^{-7}$ \\
        40 & $2.10 \times 10^{-4}$ & $7.27 \times 10^{-6}$ & $5.52 \times 10^{-7}$ & $2.21 \times 10^{-7}$ \\
        \hline
        \textbf{RoC} & $3.06$ & $4.56$ & $5.76$ & $6.04$\\
        \hline
        \end{tabular}
        \captionof{table}{\label{tab:euler_pulse_error} Convergence of the density error in the $L^\infty$ norm at $t=1$ for the near-vacuum convecting density pulse at varying approximation order with \emph{continuously} bounds-preserving limiting and positivity-preserving bounds.}
    \end{figure}
    
\subsubsection{Leblanc shock tube}
The extension to strong shocks and extreme flow conditions was performed through simulation of the Leblanc shock tube, a notoriously difficult test case. The domain is taken as 
$\Omega = [0,9]$ with Dirichlet boundary conditions, and the initial conditions are given as
\begin{equation}
    \mathbf{Q}_0 (x) = \begin{cases}
        \mathbf{Q}_L &\mbox{if } x \leq 3, \\
        \mathbf{Q}_R, &\mbox{else},
    \end{cases} 
\end{equation}
where 
\begin{equation}
    \mathbf{Q}_L = 
    \begin{bmatrix}
        1\\
        0\\
        (\gamma -1) 10^{-1}
    \end{bmatrix}
    \quad \text{and} \quad
    \mathbf{Q}_R = 
    \begin{bmatrix}
        10^{-3}\\
        0\\
        (\gamma -1) 10^{-10}
    \end{bmatrix}.
\end{equation}
The specific heat ratio is set as $\gamma = 5/3$. Due to extreme pressure ratios on the order of $10^9$, the problem is extremely challenging for numerical schemes, and a more in-depth description of the numerical difficulties of this case is presented in \citet{Liu2009}. Among them is the tendency for numerical schemes to overpredict the specific internal energy near the contact discontinuity, such that the speed of the shock front is overestimated. This overestimation in the shock speed typically requires an extremely fine resolution to mitigate. Since the purpose of this experiment was to evaluate the bounds-preserving capabilities of the proposed approach, more reasonable resolution levels were used, such that minor overpredictions in the shock speed are to be expected.

\begin{figure}[htbp!]
    \centering
    \subfloat[Density]{\adjustbox{width=0.48\linewidth,valign=b}{    \begin{tikzpicture}[spy using outlines={rectangle, height=6cm,width=2.3cm, magnification=3, connect spies}]
		\begin{semilogyaxis}[name=plot1,
		    axis x line=left,
            axis y line=left,
		    xlabel={$x$},
		    xtick={0,1,2,3,4,5,6,7,8,9},
    		xmin=0,
    		xmax=9,
    		x tick label style={
        		/pgf/number format/.cd,
            	fixed,
            	fixed zerofill,
            	precision=1,
        	    /tikz/.cd},
    		ylabel={$\rho$},
    		ylabel style={rotate=-90},
		    ytick={1e-3, 1e-2, 1e-1, 1e-0},
    		ymin=8e-4,
    		ymax=1.02,
    		legend style={at={(0.03,0.03)},anchor=south west ,font=\normalsize, column sep=0.2cm},
    		legend cell align={left},
    		style={font=\normalsize}
      ]

            \addplot [color=black, style={ultra thin}, only marks, mark=o, mark options={scale=0.5}, mark repeat = 5, mark phase = 0] table[x=x,y=r,col sep=comma,unbounded coords=jump]{./figs/data/leblanc_ref.csv};
    		\addlegendentry{Exact}

			\addplot[color={black}, style={thick, dashed}]
				table[x=x,y=r,col sep=comma,unbounded coords=jump]{./figs/data/euler_leblanc_p3_n1200_discrete_solution.csv};
    		\addlegendentry{Discrete limiting}
   
			\addplot[color={black}, style={thick}]
				table[x=x,y=r,col sep=comma,unbounded coords=jump]{./figs/data/euler_leblanc_p3_n1200_continuous_solution.csv};
    		\addlegendentry{Continuous limiting}

		\end{semilogyaxis} 		    
		
	\end{tikzpicture}}}
    \subfloat[Velocity]{\adjustbox{width=0.48\linewidth,valign=b}{    \begin{tikzpicture}[spy using outlines={rectangle, height=6cm,width=2.3cm, magnification=3, connect spies}]
		\begin{axis}[name=plot1,
		    axis x line=left,
            axis y line=left,
		    xlabel={$x$},
		    xtick={0,1,2,3,4,5,6,7,8,9},
    		xmin=0,
    		xmax=9,
    		x tick label style={
        		/pgf/number format/.cd,
            	fixed,
            	fixed zerofill,
            	precision=1,
        	    /tikz/.cd},
    		ylabel={$v$},
    		ylabel style={rotate=-90},
		    ytick={0,0.2,0.4,0.6,0.8},
    		ymin=-0.05,
    		ymax=0.8,
    		y tick label style={
        		/pgf/number format/.cd,
            	fixed,
            	fixed zerofill,
            	precision=1,
        	    /tikz/.cd},
    		legend style={at={(1,1)},anchor=north east ,font=\scriptsize, column sep=0.2cm},
    		legend cell align={left},
    		style={font=\normalsize}]

            \addplot [color=black, style={ultra thin}, only marks, mark=o, mark options={scale=0.5}, mark repeat = 5, mark phase = 2] table[x=x,y=u,col sep=comma,unbounded coords=jump]{./figs/data/leblanc_ref.csv};

			\addplot[color={black}, style={thick, dashed}]
				table[x=x,y=u,col sep=comma,unbounded coords=jump]{./figs/data/euler_leblanc_p3_n1200_discrete_solution.csv};

			\addplot[color={black}, style={thick}]
				table[x=x,y=u,col sep=comma,unbounded coords=jump]{./figs/data/euler_leblanc_p3_n1200_continuous_solution.csv};

		\end{axis} 		    
		
	\end{tikzpicture}}}
    \newline
    \subfloat[Pressure]{\adjustbox{width=0.48\linewidth,valign=b}{    \begin{tikzpicture}[spy using outlines={rectangle, height=6cm,width=2.3cm, magnification=3, connect spies}]
		\begin{semilogyaxis}[name=plot1,
		    axis x line=left,
            axis y line=left,
		    xlabel={$x$},
		    xtick={0,1,2,3,4,5,6,7,8,9},
    		xmin=0,
    		xmax=9,
    		x tick label style={
        		/pgf/number format/.cd,
            	fixed,
            	fixed zerofill,
            	precision=1,
        	    /tikz/.cd},
    		ylabel={$P$},
    		ylabel style={rotate=-90},
		    ytick={1e-12, 1e-10, 1e-8, 1e-6, 1e-4, 1e-2, 1e-0},
    		ymin=1e-12,
    		ymax=1,
    		legend style={at={(1,1)},anchor=north east ,font=\scriptsize, column sep=0.2cm},
    		legend cell align={left},
    		style={font=\normalsize}]

            \addplot [color=black, style={ultra thin}, only marks, mark=o, mark options={scale=0.5}, mark repeat = 5, mark phase = 0] table[x=x,y=p,col sep=comma,unbounded coords=jump]{./figs/data/leblanc_ref.csv};

			\addplot[color={black}, style={thick, dashed}]
				table[x=x,y=p,col sep=comma,unbounded coords=jump]{./figs/data/euler_leblanc_p3_n1200_discrete_solution.csv};
      
			\addplot[color={black}, style={thick}]
				table[x=x,y=p,col sep=comma,unbounded coords=jump]{./figs/data/euler_leblanc_p3_n1200_continuous_solution.csv};

		\end{semilogyaxis} 		    
		
	\end{tikzpicture}}}
    \subfloat[Specific internal energy]{\adjustbox{width=0.48\linewidth,valign=b}{    \begin{tikzpicture}[spy using outlines={rectangle, height=6cm,width=2.3cm, magnification=3, connect spies}]
		\begin{axis}[name=plot1,
		    axis x line=left,
            axis y line=left,
		    xlabel={$x$},
		    xtick={0,1,2,3,4,5,6,7,8,9},
    		xmin=0,
    		xmax=9,
    		x tick label style={
        		/pgf/number format/.cd,
            	fixed,
            	fixed zerofill,
            	precision=1,
        	    /tikz/.cd},
    		ylabel={$e$},
    		ylabel style={rotate=-90},
		    ytick={0, 0.05, 0.10, 0.15, 0.20},
    		ymin=-0.02,
    		ymax=0.22,
    		y tick label style={
        		/pgf/number format/.cd,
            	fixed,
            	fixed zerofill,
            	precision=2,
        	    /tikz/.cd},
    		legend style={at={(1,1)},anchor=north east ,font=\scriptsize, column sep=0.2cm},
    		legend cell align={left},
    		style={font=\normalsize}]

            \addplot [color=black, style={ultra thin}, only marks, mark=o, mark options={scale=0.5}, mark repeat = 5, mark phase = 2] table[x=x,y=e,col sep=comma,unbounded coords=jump]{./figs/data/leblanc_ref.csv};

			\addplot[color={black}, style={thick, dashed}]
				table[x=x,y expr={1.5*\thisrow{p}/\thisrow{r}},col sep=comma,unbounded coords=jump]{./figs/data/euler_leblanc_p3_n1200_discrete_solution.csv};
    
			\addplot[color={black}, style={thick}]
				table[x=x,y expr={1.5*\thisrow{p}/\thisrow{r}},col sep=comma,unbounded coords=jump]{./figs/data/euler_leblanc_p3_n1200_continuous_solution.csv};

		\end{axis} 		    
		
	\end{tikzpicture}}}
    \caption{\label{fig:leblanc_full} Profiles of density (top left), velocity (top right), pressure (bottom left), and specific internal energy (bottom right) for the Leblanc shock tube problem at $t = 6$ computed with a $\mathbb P_3$ approximation and $N = 1200$ elements. Results shown for positivity-preserving bounds with continuously/discretely bounds-preserving limiting augmented with a subcell slope limiter.}
\end{figure}
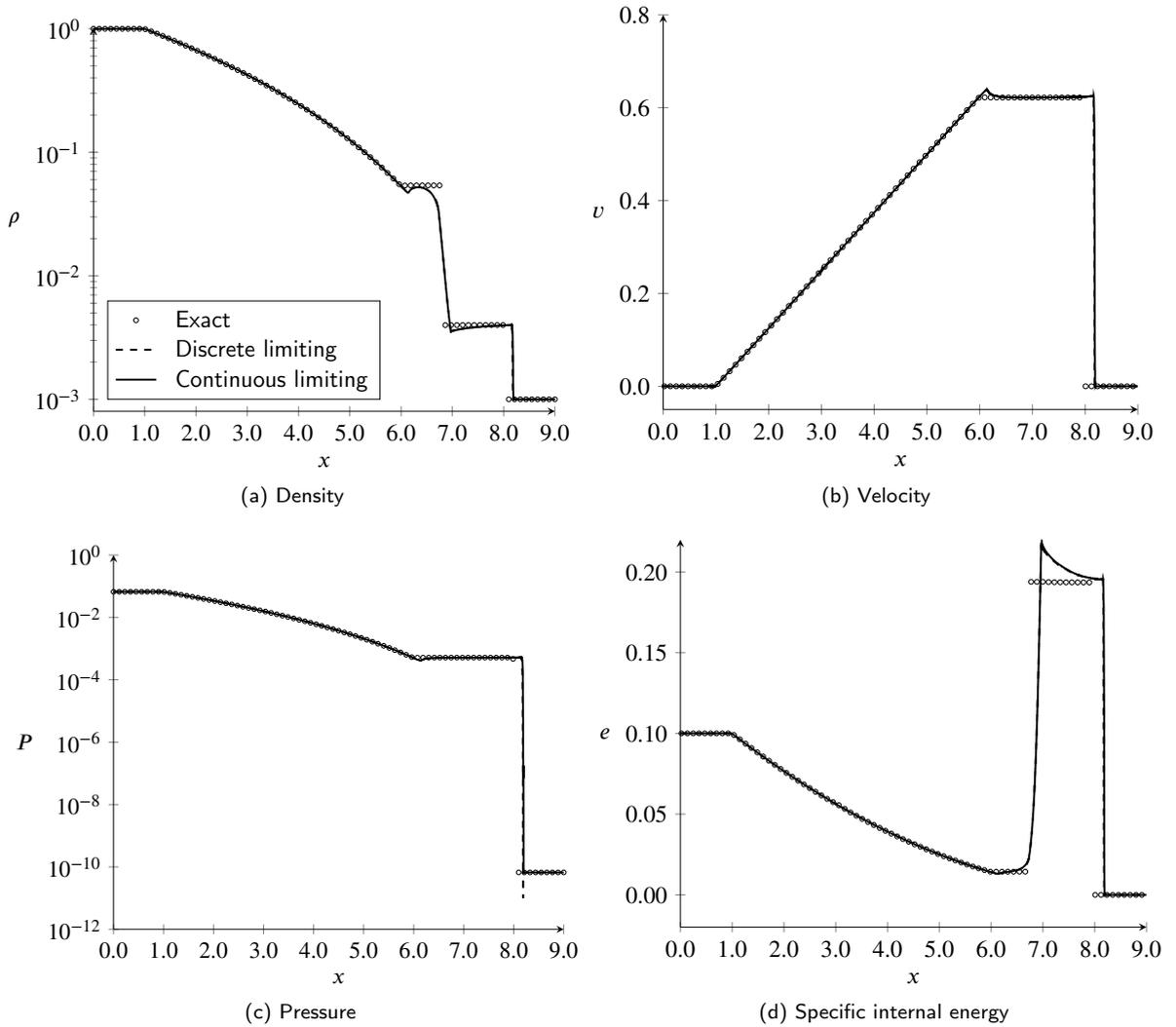
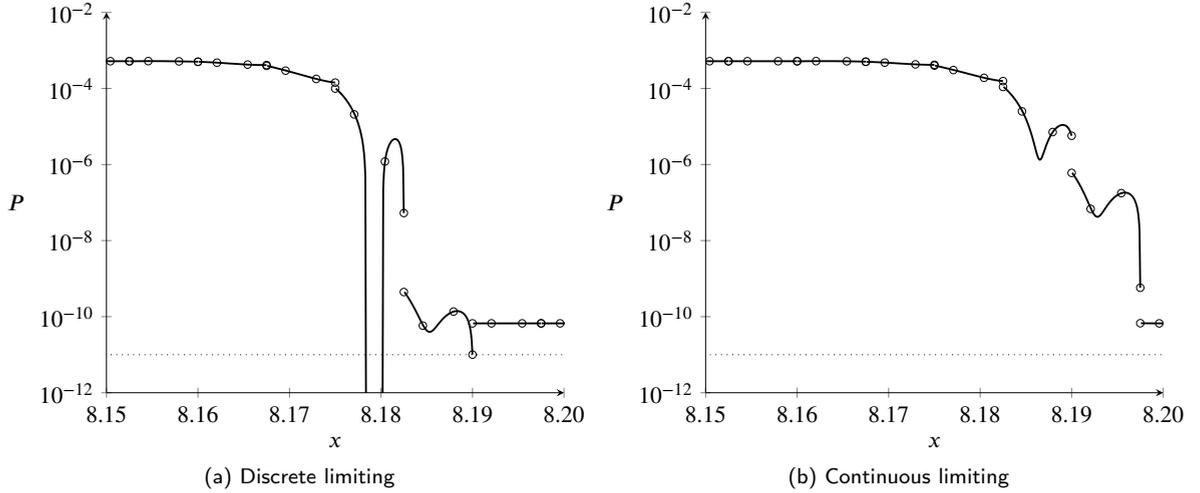
\begin{figure}[htbp!]
    \centering
    \subfloat[Discrete limiting]{\adjustbox{width=0.48\linewidth,valign=b}{    \begin{tikzpicture}[spy using outlines={rectangle, height=6cm,width=2.3cm, magnification=3, connect spies}]
		\begin{semilogyaxis}[name=plot1,
		    axis x line=left,
            axis y line=left,
		    xlabel={$x$},
		    xtick={8.15, 8.16, 8.17, 8.18, 8.19, 8.2},
    		xmin=8.15,
    		xmax=8.2,
    		x tick label style={
        		/pgf/number format/.cd,
            	fixed,
            	fixed zerofill,
            	precision=2,
        	    /tikz/.cd},
    		ylabel={$P$},
    		ylabel style={rotate=-90},
		    ytick={1e-12, 1e-10, 1e-8, 1e-6, 1e-4, 1e-2},
    		ymin=1e-12,
    		ymax=1e-2,
    		legend style={at={(1,1)},anchor=north east ,font=\scriptsize, column sep=0.2cm},
    		legend cell align={left},
    		style={font=\normalsize}]

            \foreach \x in {1086,...,1093}{
			\addplot[color={black}, style={thick}]
				table[x=x,y expr={max(1e-13, \thisrow{p})},col sep=comma,unbounded coords=jump]{./figs/data/euler_leblanc_p3_n1200_discrete_solution_e\x.csv};
            }
            
			\addplot[color={black}, style={}, only marks, mark=o, mark options={scale=0.8}, restrict expr to domain={x}{8.1:8.25}]
				table[x=x,y=p,col sep=comma,unbounded coords=jump]{./figs/data/euler_leblanc_p3_n1200_discrete_solution.csv};

            \draw [dotted] (8.1,1e-11) -- (8.25,1e-11);

		\end{semilogyaxis} 		    
		
	\end{tikzpicture}}}
    \subfloat[Continuous limiting]{\adjustbox{width=0.48\linewidth,valign=b}{    \begin{tikzpicture}[spy using outlines={rectangle, height=6cm,width=2.3cm, magnification=3, connect spies}]
		\begin{semilogyaxis}[name=plot1,
		    axis x line=left,
            axis y line=left,
		    xlabel={$x$},
		    xtick={8.15, 8.16, 8.17, 8.18, 8.19, 8.2},
    		xmin=8.15,
    		xmax=8.2,
    		x tick label style={
        		/pgf/number format/.cd,
            	fixed,
            	fixed zerofill,
            	precision=2,
        	    /tikz/.cd},
    		ylabel={$P$},
    		ylabel style={rotate=-90},
		    ytick={1e-12, 1e-10, 1e-8, 1e-6, 1e-4, 1e-2},
    		ymin=1e-12,
    		ymax=1e-2,
    		legend style={at={(1,1)},anchor=north east ,font=\scriptsize, column sep=0.2cm},
    		legend cell align={left},
    		style={font=\normalsize}]

            \foreach \x in {1086,...,1093}{
			\addplot[color={black}, style={thick}]
				table[x=x,y expr={max(1e-13, \thisrow{p})},col sep=comma,unbounded coords=jump]{./figs/data/euler_leblanc_p3_n1200_continuous_solution_e\x.csv};
            }
            
			\addplot[color={black}, style={}, only marks, mark=o, mark options={scale=0.8}, restrict expr to domain={x}{8.1:8.25}]
				table[x=x,y=p,col sep=comma,unbounded coords=jump]{./figs/data/euler_leblanc_p3_n1200_continuous_solution.csv};

            \draw [dotted] (8.1,1e-11) -- (8.25,1e-11);

		\end{semilogyaxis} 		    
		
	\end{tikzpicture}}}
    \caption{\label{fig:leblanc_zoom}
    Enlarged view of the pressure profiles near the shock front for the Leblanc shock tube problem at $t = 6$ computed with a $\mathbb P_3$ approximation and $N = 1200$ elements. Results shown for positivity-preserving bounds with continuously/discretely bounds-preserving limiting augmented with a subcell slope limiter. Circle markers denote discrete solution nodes. Minimum pressure bound (isocontour of $g_2(\mathbf{u}) = P_{\min}$) shown by dotted line.}
\end{figure}

The problem was computed using a $\mathbb P_3$ scheme with $N = 1200$ elements augmented with a subcell slope limiter. For this case, simply utilizing the subcell slope limiter without positivity-preserving limiting still quickly resulted in the divergence of the numerical scheme. A comparison of the discrete limiting method with the proposed continuous limiting method is shown in the density, velocity, pressure, and specific internal energy profiles at $t = 6$ in \cref{fig:leblanc_full}. Both approaches quantitatively achieved similar result, with the only visually discernible difference being an underprediction in the pressure field around the shock front by the discrete limiting method. 

The distinction in the proposed continuous limiting approach became evident when observing the solution polynomial in the vicinity of the shock. An enlarged view of the pressure profiles near the shock front is shown in \cref{fig:leblanc_zoom} for both the discrete limiting and continuous limiting approach, with the discrete solution nodes shown by markers and the full solution polynomial shown by lines. It can be seen that while the discrete limiting method ensured that the pressure at the discrete solution nodes was positive, it did not prevent undershoots in the polynomial approximation from causing negative pressure values. In fact, the minimum continuous pressure for the discrete limiting method at $t=6$, evaluated at 100 points per element, was $-3.2{\cdot}10^{-6}$. In contrast, the continuous limiting approach showed a positive pressure profile across the entire solution polynomial, even across the strong shock front. The minimum continuous pressure for the proposed continuous limiting approach at $t=6$ was $6.6{\cdot}10^{-11}$. This value aligns with the minimum initial pressure (and the minimum pressure in the exact solution) which indicates that the additional dissipation in the proposed approach is better suited at mitigating spurious undershoots even at the discrete level.

\subsubsection{Sedov blast}
For the final numerical experiment and a validation of the proposed approach for two-dimensional gas dynamics problems, the Sedov blast wave~\citep{Sedov1993} was considered. In this problem, an energy point-source in an ambient gas drives an outward-running radial explosion. The problem setup is taken identically to the work of \citet{Maire2009}, although the entire domain is simulated instead of a single quadrant. The domain is set as $\Omega = [-1.2, 1.2]^2$ with periodic boundary conditions. Due to the finite rate of propagation of the shock wave, the boundary conditions do not impact the solution over the given simulation time. The initial conditions consist of an ambient gas with a constant density $\rho_0 = 1$ and velocity $\mathbf{v} = \mathbf{0}$. Furthermore, the pressure in the domain is set to an ambient value of $P_a = 10^{-6}$ everywhere except in the center-most element (i.e., the element centered at $x = y = 0$), where an overpressure value is set in terms of the internal energy $e_0$ and cell volume $V_0$ as
\begin{equation}
    P_0 = 4(\gamma-1) \rho_0 e_0/V_0.
\end{equation}
The element volume is simply computed as $V_0 = \Delta x \Delta y = 2.4^2/N$, and the initial specific internal energy is set as $\epsilon_0 = 0.244816$ which yields a shock front with a maximum density of $\rho = 6$ and a radius of 1 at $t = 1$. The factor of 4 in the initial pressure is used to account for the fact that the entire domain is being simulated instead of one quadrant.

\begin{figure}[htbp!]
        \centering
        \subfloat[Density]{\adjustbox{width=0.33\linewidth,valign=b}{\includegraphics{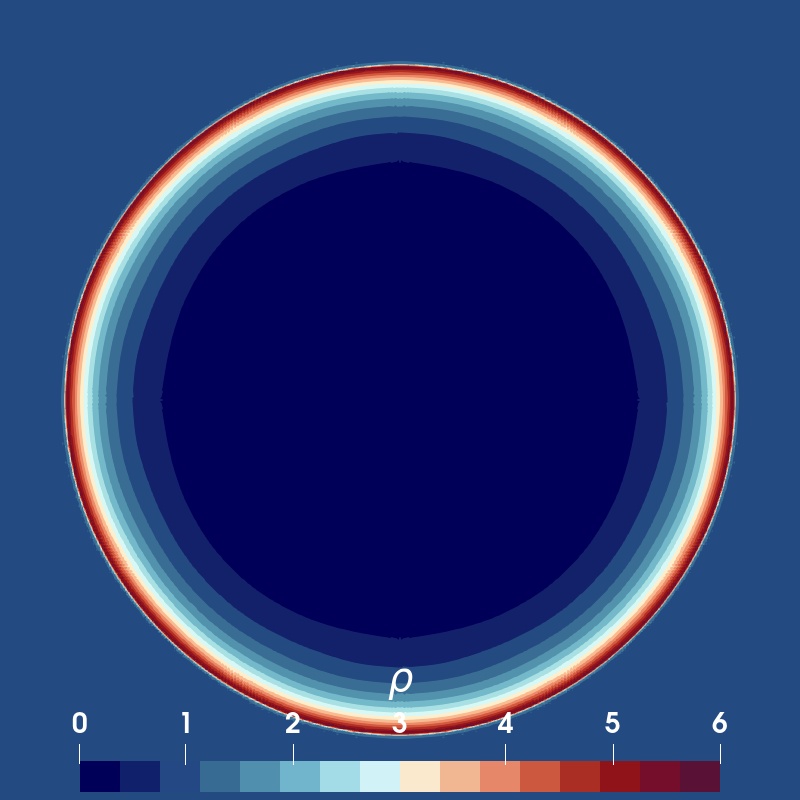}}}
        \subfloat[Radial velocity]{\adjustbox{width=0.33\linewidth,valign=b}{\includegraphics{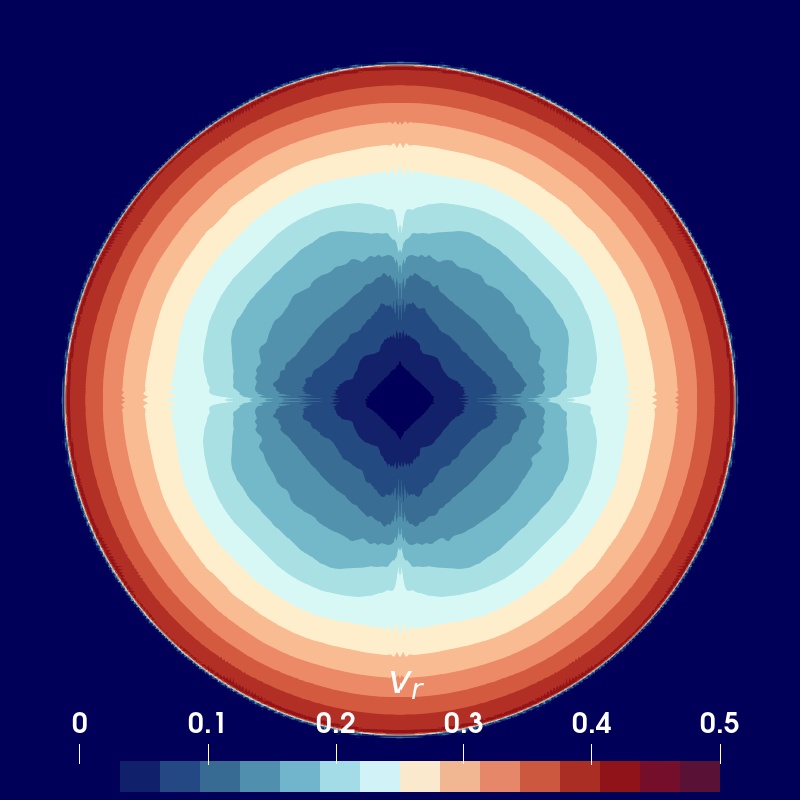}}}
        \subfloat[Pressure]{\adjustbox{width=0.33\linewidth,valign=b}{\includegraphics{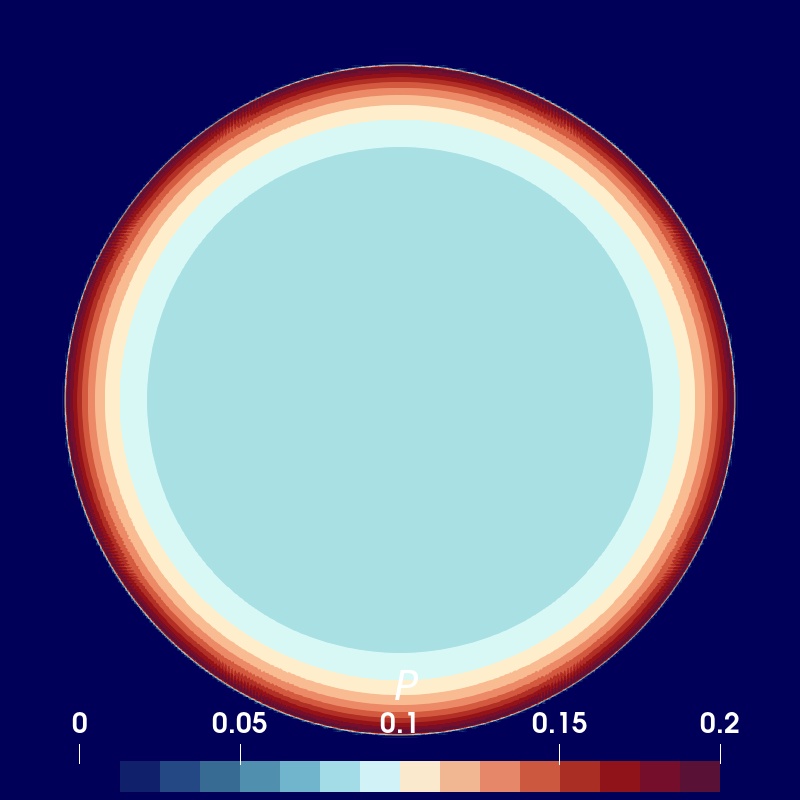}}}
        \caption{\label{fig:sedov_contours} Contours of density (left), radial velocity (middle), and pressure (right) for the two-dimensional Sedov blast problem at $t = 1$ computed with a $\mathbb P_2$ approximation, $N = 257^2$ elements, positivity-preserving bounds, and continuously bounds-preserving limiting augmented with a subcell slope limiter.} 
    \end{figure}

\begin{figure}[htbp!]
    \centering
    \subfloat[Density]{\adjustbox{width=0.33\linewidth,valign=b}{    \begin{tikzpicture}[spy using outlines={rectangle, height=6cm,width=2.3cm, magnification=3, connect spies}]
		\begin{axis}[name=plot1,
		    axis x line=left,
            axis y line=left,
		    xlabel={$r$},
		    xtick={0,0.2, 0.4, 0.6, 0.8, 1.0, 1.2, 1.4},
    		xmin=0,
    		xmax=1.4,
    		x tick label style={
        		/pgf/number format/.cd,
            	fixed,
            	fixed zerofill,
            	precision=1,
        	    /tikz/.cd},
    		ylabel={$\rho$},
    		ylabel style={rotate=-90},
		    ytick={0,1,2,3,4,5,6},
    		ymin=-0.1,
    		ymax=6,
    		legend style={at={(0.03,0.97)},anchor=north west ,font=\normalsize, column sep=0.2cm},
    		legend cell align={left},
    		style={font=\normalsize}
      ]

            \addplot [color=black, style={ultra thin}, only marks, mark=o, mark options={scale=0.5}, mark repeat = 2, mark phase = 0] table[x=r,y=d,col sep=comma,unbounded coords=jump]{./figs/data/sedov_exact2.csv};
    		\addlegendentry{Exact}

			\addplot[color={black}, style={thick}] table[x=r,y=d,col sep=comma,unbounded coords=jump]{./figs/data/euler_sedov_p2_n257_continuous.csv};
    		\addlegendentry{Present work}
      
            \addplot[color=black, style={ultra thin}, only marks, mark=o, mark options={scale=0.5}, samples at={1, 1.02, ..., 1.4}]{1};

		\end{axis} 		    
		
	\end{tikzpicture}}}
    \subfloat[Radial velocity]{\adjustbox{width=0.33\linewidth,valign=b}{    \begin{tikzpicture}[spy using outlines={rectangle, height=6cm,width=2.3cm, magnification=3, connect spies}]
		\begin{axis}[name=plot1,
		    axis x line=left,
            axis y line=left,
		    xlabel={$r$},
		    xtick={0,0.2, 0.4, 0.6, 0.8, 1.0, 1.2, 1.4},
    		xmin=0,
    		xmax=1.4,
    		x tick label style={
        		/pgf/number format/.cd,
            	fixed,
            	fixed zerofill,
            	precision=1,
        	    /tikz/.cd},
    		ylabel={$v_r$},
    		ylabel style={rotate=-90},
		    ytick={0,0.1,0.2,0.3,0.4},
    		ymin=-0.02,
    		ymax=0.45,
    		legend style={at={(0.03,0.03)},anchor=south west ,font=\normalsize, column sep=0.2cm},
    		legend cell align={left},
    		style={font=\normalsize}
      ]
    		
            \addplot [color=black, style={ultra thin}, only marks, mark=o, mark options={scale=0.5}, mark repeat = 2, mark phase = 0] table[x=r,y=v,col sep=comma,unbounded coords=jump]{./figs/data/sedov_exact2.csv};

			\addplot[color={black}, style={thick}] table[x=r,y=v,col sep=comma,unbounded coords=jump]{./figs/data/euler_sedov_p2_n257_continuous.csv};
      
            \addplot[color=black, style={ultra thin}, only marks, mark=o, mark options={scale=0.5}, samples at={1, 1.02, ..., 1.4}]{0};

		\end{axis} 		    
		
	\end{tikzpicture}}}
    \subfloat[Pressure]{\adjustbox{width=0.33\linewidth,valign=b}{    \begin{tikzpicture}[spy using outlines={rectangle, height=6cm,width=2.3cm, magnification=3, connect spies}]
		\begin{axis}[name=plot1,
		    axis x line=left,
            axis y line=left,
		    xlabel={$r$},
		    xtick={0,0.2, 0.4, 0.6, 0.8, 1.0, 1.2, 1.4},
    		xmin=0,
    		xmax=1.4,
    		x tick label style={
        		/pgf/number format/.cd,
            	fixed,
            	fixed zerofill,
            	precision=1,
        	    /tikz/.cd},
    		ylabel={$P$},
    		ylabel style={rotate=-90},
		    ytick={0,0.05,0.1,0.15,0.2},
    		ymin=-0.01,
    		ymax=0.2,
    		y tick label style={
        		/pgf/number format/.cd,
            	fixed,
            	precision=2,
        	    /tikz/.cd},
    		legend style={at={(0.03,0.03)},anchor=south west ,font=\normalsize, column sep=0.2cm},
    		legend cell align={left},
    		style={font=\normalsize}
      ]
    		
            \addplot [color=black, style={ultra thin}, only marks, mark=o, mark options={scale=0.5}, mark repeat = 2, mark phase = 0] table[x=r,y=p,col sep=comma,unbounded coords=jump]{./figs/data/sedov_exact2.csv};

			\addplot[color={black}, style={thick}] table[x=r,y=p,col sep=comma,unbounded coords=jump]{./figs/data/euler_sedov_p2_n257_continuous.csv};
      
            \addplot[color=black, style={ultra thin}, only marks, mark=o, mark options={scale=0.5}, samples at={1, 1.02, ..., 1.4}]{2e-6};

		\end{axis} 		    
		
	\end{tikzpicture}}}
    
    \caption{\label{fig:sedov_profiles} 
    Cross-sectional profiles of density (left), radial velocity (middle), and pressure (right) for the two-dimensional Sedov blast problem at $t = 1$ along the positive $x=y$ diagonal computed with a $\mathbb P_2$ approximation, $N = 257^2$ elements, positivity-preserving bounds, and continuously bounds-preserving limiting augmented with a subcell slope limiter.
    }
\end{figure}
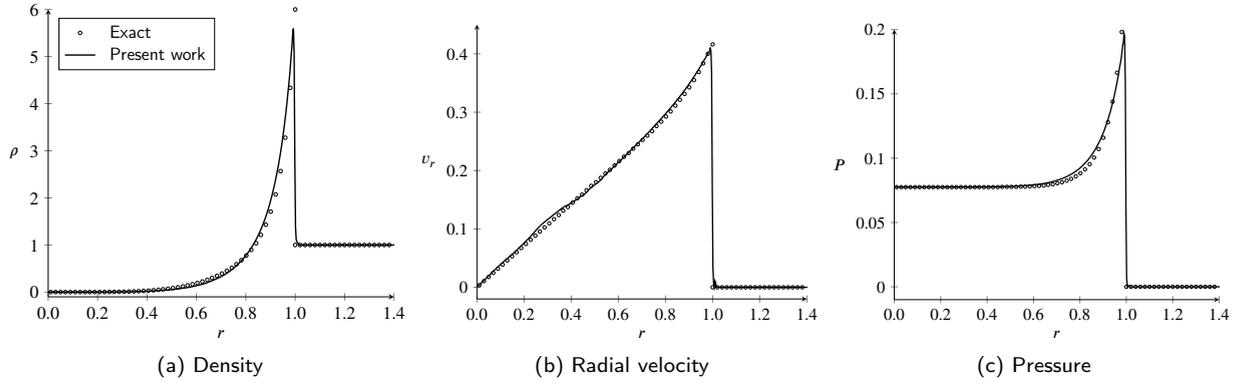

The problem was computed with a $\mathbb P_2$ approximation and $N = 257^2$ elements using the proposed continuously bounds-preserving limiter augmented with a subcell slope limiter. Like with the Leblanc shock tube case, simply using the subcell slope limiter without enforcing positivity constraints resulted in the divergence of the numerical scheme. The contours of density, radial velocity, and pressure at $t=1$ are shown in \cref{fig:sedov_contours}. The typical blast wave structure was observed, with a sharp shock front in the density and pressure field and a smoothly varying velocity field in the interior. Minor mesh imprinting was observed in the velocity field. Additionally, the solution profiles on the positive diagonal cross-section $x = y$ ($x, y \geq 0$) are shown in \cref{fig:sedov_profiles} in comparison to the analytic solution. Good agreement with the analytic solution was obtained with the proposed approach. To evaluate the continuously bounds-preserving property of the proposed scheme, the minimum density and pressure in the solution at $t = 1$ were calculated across $100^2$ points per element. At the final time, the minimum density was $1.5{\cdot}10^{-3}$ and the minimum pressure was $5.9{\cdot}10^{-11}$, well above the minimum density/pressure bounds imposed. 

\section{Conclusion}\label{sec:conclusion}
In this work, we presented a novel limiting approach for ensuring that DG approximations of hyperbolic conservation laws can enforce constraints continuously across the entire solution polynomial. The proposed method relied on introducing a novel modification to the constraint functionals which allow for applying bounded optimization techniques to find a limiting factor that ensures the limited solution is continuously bounds-preserving. An efficient optimization approach was then introduced for finding this limiting factor. The proposed continuously bounds-preserving limiting method was applied to a variety of hyperbolic conservation laws on structured and unstructured meshes, ranging from linear scalar transport to compressible gas dynamics with strong shocks, showcasing the ability of the approach in ensuring constraints are satisfied continuously. 

The proposed approach has the potential for applications outside of the ones shown in this work as well as further possible developments. Other applications include moving/deforming meshes such as in arbitrary Lagrangian--Eulerian codes and enforcing constraints outside of the element for reconstruction-based methods. Improvements can also be made in the optimization approach to ensure that the given solution is the true global minimum, such as combining the optimization approach with a bounding method~\citep{Johnen2013} on the constraint functional or even the Hessian itself. Further improvements in the accuracy can be made by utilizing higher-order approximations of the limiting factor (or iterative limiting procedures) for nonlinear constraints as well as modifying the formulation to allow for enforcing bounds which vary across the element.

\section*{Acknowledgements}
\label{sec:ack}

This work was performed under the auspices of the U.S. Department of Energy by Lawrence Livermore National Laboratory under contract DE--AC52--07NA27344 and the LLNL-LDRD Program under Project tracking No.\ 24--ERD--050. Release number LLNL--JRNL--858838--DRAFT.

\bibliographystyle{unsrtnat}
\bibliography{reference}



\newpage
\appendix
\section{Algorithmic implementation}\label{app:algs}

\begin{algorithm}
\caption{Algorithm for projecting optimizer step to element bounds }
{INPUT}: $\mathbf{x}_1$, $\Delta \mathbf{x}$, $\mathbf{n}$ //  Initial spatial location, optimizer step, and boundary normal vector \newline
{OUTPUT}: $\Delta \mathbf{x}'$ // Projected optimizer step \newline

Project($\mathbf{x}_1$, $\Delta \mathbf{x}$, $\mathbf{n}$) :
\begin{algorithmic} \label{alg:proj}

    \item[] // Compute next candidate spatial location.
    \item[]$\mathbf{x}_2 = \mathbf{x}_1 + \Delta \mathbf{x}$
    \item[] // Apply projection if next spatial location is out of element bounds.
    \IF{$\mathbf{x}_2\notin \Omega_k$}
        \item[] // If initial spatial location is on element boundary, project search direction along boundary.
        \item[] // Else project along search direction to boundary. 
        \IF{$\mathbf{x}_1 \in \partial \Omega_k$}
            \item[] $\Delta \mathbf{x}' = (1 - \mathbf{n})\odot\Delta \mathbf{x}$
        \ELSE
            \item[] $\Delta \mathbf{x}' = \eta \Delta \mathbf{x}$ s.t. $\mathbf{x}_1 + \eta \Delta \mathbf{x} \in \partial \Omega_k$
        \ENDIF
    \ELSE
        \item[] $\Delta \mathbf{x}' = \Delta \mathbf{x}$
    \ENDIF
\end{algorithmic}
\end{algorithm}

\begin{algorithm}
\caption{Optimization and limiting algorithm}
{INPUT}: $\mathbf{u}$, $\mathbf{n}$ // Solution, element boundary normal vector \newline
{OUTPUT}: $\mathbf{u}$ // Limited solution \newline

Limit($\mathbf{u}$, $\mathbf{n}$) :
\begin{algorithmic} \label{alg:opt}
\item[]$\mathbf{u} \leftarrow$  Convert nodal basis to monomial basis
\item[]$\mathbf{x}_{\min}, h_{\min}\leftarrow$ Compute discrete minimum of $h$ across solution/quadrature nodes
\item[] // Set minimum/initial values
\item[]$h^* = h_{\min}$
\item[]$\Delta \mathbf{x} = 0$

\item[] // Iterate optimizer.
\FOR{$n \in \lbrace 1, ..., \text{n}_{\text{iters}} \rbrace$}
\item[] $\mathbf{J}$, $\mathbf{H} \leftarrow$ Compute $\mathbf{J}$ and $\mathbf{H}$ as per \cref{eq:jacobian} and \cref{eq:hessian}
\item[] // If Hessian is positive semi-definite, compute projected Newton--Raphson step.  
\item[] // Else compute projected gradient descent step with backtracking line search.  
\IF{$det(\mathbf{H}) > \epsilon$ and $\mathbf{H}$ is PSD}
    \item[] $\Delta \mathbf{x} = \mathbf{H}^{-1} \mathbf{J}$
    \item[] $\Delta \mathbf{x} = \text{Project} (\mathbf{x}_{\min}^n, \Delta \mathbf{x}, \mathbf{n})$ // Compute projection (if necessary) as per \cref{alg:proj}.
\ELSE
    \item[] $\beta = \beta_0$
    \item[] $\Delta \mathbf{x} = -\beta \mathbf{J}$
    \item[] $\Delta \mathbf{x} = \text{Project}(\mathbf{x}_{\min}^n, \Delta \mathbf{x}, \mathbf{n})$ // Compute projection (if necessary) as per \cref{alg:proj}.
    \item[] // Iterate backtracking line search.
    \FOR{$i \in \lbrace 1, ..., 5\rbrace$}
        \item[] $\beta = 0.5 \beta$
        \item[] $\Delta \mathbf{x} = -\beta \mathbf{J}$
        \IF{$h(\mathbf{u}(\mathbf{x}_{n-1} + \Delta x)) \leq h_{\min} - 0.5\beta\|\mathbf{J}\|_2$}
            \item[] \textbf{break}
        \ENDIF
    \ENDFOR    
\ENDIF
\item[] // Set next spatial location.
\item[] $\mathbf{x}_{\min}^n = \mathbf{x}_{\min}^{n-1} + \Delta \mathbf{x}$
\item[] $h_{\min} = h(\mathbf{u}(\mathbf{x}_{\min}^n))$
\item[] $h^* = \min(h^*, h_{\min})$
\ENDFOR
\item[] // Compute and apply tolerance correction.
\item[] $\Delta h = \|\mathbf{J}\|_2 \|\Delta \mathbf{x}\|_2$
\item[] $h^* = \max(-1, h^* - \Delta h)$
\item[] // Compute limiting factor and limit solution.
\item[] $\alpha = \max(0, -h^*)$
\item[] $\mathbf{u} = (1-\alpha)\mathbf{u} + \alpha \overline{\mathbf{u}}$
\end{algorithmic}
\end{algorithm}

\section{Subcell slope limiter}\label{app:subcell}
The slope limiter relies on formulating a high-order (unlimited) approximation, a low-order (slope-limited) approximation, and an indicator metric for switching between the two. The high-order candidate approximation, denoted by $\mathbf{u}_{H}$, is simply the temporal update from the underlying high-order DG scheme. The low-order candidate approximation, denoted by $\mathbf{u}_{L}$, is computed using a subcell approximation as described by \citet{RuedaRamrez2022}, Section 2.2. The temporal update is then performed by switching between the candidate solutions based on an \textit{a posteriori} indicator metric $\pi(\mathbf{u}_{L}, \mathbf{u}_{H})$ as 
\begin{equation}
    \mathbf{u}^{n+1} = \begin{cases}
        \mathbf{u}_{L}^{n+1},\quad\, \text{ if } \pi(\mathbf{u}_{L}^{n+1}, \mathbf{u}_{H}^{n+1}) = 1,\\
        \mathbf{u}_{H}^{n+1},\quad \text{ else.}
    \end{cases}
\end{equation}
The indicator used in this work is a discrete maximum principle indicator similar to the multi-dimensional
optimal order detection (MOOD) indicator in the work of \citet{Dumbser2014} and \citet{Dumbser2016}. Given some indicator variable of the solution $\mu(\mathbf{u})$, the indicator metric is computed as
\begin{equation}
    \pi(\mathbf{u}_{L}, \mathbf{u}_{H}) = \begin{cases}
        0,\quad \text{ if } \mu_{\min} \leq \mu(\mathbf{u}_{H,i}) \leq \mu_{\max} \ \forall \ i \in S,\\
        1,\quad \text{ else,}
    \end{cases}
\end{equation}
where $S = \{1, ..., n_s \}$ represents the discrete solution nodes. The bounds $\mu_{\min}$ and $\mu_{\min}$ are computed through a relaxed discrete maximum principle on the low-order solution, i.e.,
\begin{equation}
    \mu_{\min} = (1 - \epsilon_D) \left(\underset{i \in S}{\min}\ \mu(\mathbf{u}_{L,i})\right)
\end{equation}
and
\begin{equation}
    \mu_{\max} = (1 + \epsilon_D) \left(\underset{i \in S}{\max}\ \mu(\mathbf{u}_{L,i})\right).
\end{equation}
$\epsilon_D$ represents a relaxation coefficient and is taken as $\epsilon_D = 10^{-2}$ in this work. For the gas dynamics equations, the density field is used for the indicator variable $\mu$. 

\end{document}